\RequirePackage{fix-cm}
\documentclass[twocolumn]{svjour3}          % twocolumn
\smartqed  % flush right qed marks, e.g., at end of proof
\usepackage{graphicx}
\usepackage[utf8]{inputenc}
\usepackage{enumerate}
\usepackage[authoryear,round,longnamesfirst]{natbib}
\usepackage{amsmath}
\usepackage{svg}
\usepackage{epsfig}
\usepackage{bm}
\usepackage{amsfonts}
\usepackage{booktabs}
\usepackage{adjustbox}
\usepackage{amsmath}
\usepackage{amsmath,amssymb,algpseudocode,algorithm}

\usepackage{algorithm}
\usepackage{algpseudocode}
\usepackage{soul}
\usepackage{diagbox}
\usepackage{mathtools}
\usepackage{mathrsfs}
\usepackage[mathscr]{eucal}
\usepackage[title]{appendix}
\usepackage[hyphens,spaces,obeyspaces]{url}
\usepackage[inline]{enumitem}
\usepackage{tcolorbox}
\usepackage{lscape}
\usepackage{listings}
\definecolor{codegreen}{rgb}{0,0.6,0}
\definecolor{codegray}{rgb}{0.5,0.5,0.5}
\definecolor{codepurple}{rgb}{0.58,0,0.82}
\definecolor{backcolour}{rgb}{0.95,0.95,0.92}
 
\lstdefinestyle{mystyle}{
    backgroundcolor=\color{backcolour},   
    commentstyle=\color{codegreen},
    keywordstyle=\color{magenta},
    numberstyle=\tiny\color{codegray},
    stringstyle=\color{codepurple},
    breakatwhitespace=false,         
    breaklines=true,                 
    captionpos=b,                    
    keepspaces=true,                 
    numbers=left,                    
    numbersep=10pt,                  
    showspaces=false,                
    showstringspaces=false,
    showtabs=false,                  
    tabsize=4
}
 
\definecolor{aurometalsaurus}{rgb}{0.43, 0.5, 0.5}
\definecolor{arsenic}{rgb}{0.23, 0.27, 0.29}

\newcommand{\norm}[1]{\left\lVert #1 \right\rVert}

\newcommand{\specialcell}[2][c]{%
	\begin{tabular}[#1]{@{}c@{}}#2\end{tabular}}
\begin{document} \sloppy

\title{A partitioned scheme for adjoint shape sensitivity analysis of fluid-structure interactions involving non-matching meshes}


\author{Reza Najian Asl \and Ihar Antonau \and Aditya Ghantasala \and Wulf G. Dettmer \and Roland Wüchner \and Kai-Uwe Bletzinger}

\authorrunning{R. Najian Asl \and I. Antonau \and A. Ghantasala \and W.G. Dettmer \and R. Wüchner \and K.-U. Bletzinger} % if too long for running head

\institute{R. Najian Asl \and I. Antonau \and A. Ghantasala \and R. Wüchner \and K.-U. Bletzinger \at
              Lehrstuhl für Statik, Technische Universität München, Arcisstr. 21, 80333 München, Germany \\
              Tel.: +49-89-28922422\\
              Fax:  +49-89-28922421\\
              \email{reza.najian-asl@tum.de}           %  \\
              \and
              W.G. Dettmer \at
              Zienkiewicz Centre for Computational Engineering, College of Engineering, 			  Swansea University, Fabian Way, Swansea SA1 8EN, Wales, UK \\
              \email{w.g.dettmer@swansea.ac.uk}
}

\date{Received: date / Accepted: date}

\maketitle

\begin{abstract}
This work presents a partitioned solution procedure to compute shape gradients in fluid-structure interaction (FSI) using black-box adjoint solvers. Special attention is paid to project the gradients onto the undeformed configuration. This is due to the mixed Lagrangian-Eulerian formulation of large-displacement FSI in this work. Adjoint FSI problem is partitioned as an assembly of well-known adjoint fluid and structural problems, without requiring expensive cross-derivatives. The sub-adjoint problems are coupled with each other by augmenting the target functions with auxiliary functions, independent of the concrete choice of the underlying adjoint formulations. The auxiliary functions are linear force-based or displacement-based functionals which are readily available in well-established single-disciplinary adjoint solvers. Adjoint structural displacements, adjoint fluid displacements, and domain-based adjoint sensitivities of the fluid are the coupling fields to be exchanged between the adjoint solvers. A reduced formulation is also derived for the case of boundary-based adjoint shape sensitivity analysis for fluids. Numerical studies show that the complete formulation computes accurate shape gradients whereas inaccuracies appear in the reduced gradients, specially in regions of strong flow gradients and near singularities. Nevertheless, reduced gradient formulations are found to be a compromise between computational costs and accuracy. Mapping techniques including nearest element interpolation and the mortar method are studied in computational adjoint FSI. It is numerically shown that the mortar method does not introduce spurious oscillations in primal and sensitivity fields along non-matching interfaces, unlike the nearest element interpolation.

\end{abstract}

\section{Introduction}
Recently, adjoint-based sensitivity analysis in fluid-structure interaction (FSI) problems has been revisited by the research community from the mathematical and, particularly, the implementation point of view. This is mainly due to increases in computational power and the growing interest from industry. Mathematically speaking, numerical methods devised for solving coupled problems can be sorted into two main categories. The first category includes Jacobian-free methods like the classical fixed-point iteration method, whereas the second category needs interdisciplinary Jacobians (cross-derivatives) or matrix-vector products of these Jacobians multiplied by unknown variables. The Jacobian-based algorithms have shown superior accuracy and performance, however, they put a burden on the coupling of black-box solvers in a partitioned procedure. Both categories are very well covered and discussed in the FSI literature (see e.g. \cite{FELIPPA2001,Dettmer2006,degroote2010performance,Sicklinger2014}), but to the authors’ knowledge, the adjoint FSI problem for shape sensitivity analysis has been driven specially by the second category. This paper presents a cross-derivative-free procedure for the adjoint shape sensitivity analysis of steady-state FSI using black-box adjoint solvers on non-matching meshes. Furthermore, the fact that in a partitioned FSI environment, primal and adjoint fluid solvers operate on the deformed fluid domain (due to structural displacements) is carefully taken into account. The spatial coupling of non-matching interfaces is also considered herein. Although this has been routinely done in previous studies, e.g. by \cite{Maute2001}, the accuracy of the sensitivity information obtained by different types of mapping algorithms has not been comparatively assessed yet.

Early attempts in the adjoint-based shape sensitivity analysis for FSI were made by \cite{Maute2001}, \cite{Lund2003} and \cite{Martins2004}. This topic of research has been followed by \citep{Marcelet2008,Mani2009,Martins2013,Jenkins2016,Zhang2018}.
Among recent trends and developments in this area, the following works are notable and addressed here. \cite{Sanchez2017} established an open-source framework for coupled adjoint-based sensitivity analysis, which is based on fixed-point iterations for the adjoint variables of the coupled system using an automatic differentiation (AD) tool. The main benefit of such an approach is that, without sacrificing the gradient accuracy, there is no need to compute and store exact Jacobians to be used in the adjoint problem, especially when higher-order schemes or complex kinematics are involved. However, applicability of AD to existing solvers might be hindered due to distinct software implementation and large memory requirements, unless particular attention is paid \citep{carnarius2010adjoint}. 
\cite{Kiviaho2017} presented a coupling framework for aeroelastic analysis and optimization using discrete adjoint-based gradients. They systematically derived the discrete adjoint corresponding to the steady aeroelastic analysis in a consistent way. Applicability of this approach might be limited in a partitioned adjoint FSI environment due to the lack of availability of the required cross-derivative terms in every software package. 

A literature review of the studies by various authors shows that the coupled-adjoint sensitivity analysis for high-fidelity aero-structural design is divided into two main formulations: a three-field formulation followed by \cite{Kiviaho2017,Zhang2017,Sanchez2017,Barcelos2006} and a two-field formulation followed by \cite{Heners2017,Stavropoulou2016,Fazzolari2007}. The three-field formulation accounts for aerodynamic, structural, and mesh deformation residuals in adjoint-based sensitivity analysis while the two-field formulation either implicitly includes or completely excludes mesh motion in the sensitivity analysis. For example, \cite{kenway2014scalable} derived a two-field-based formulation which incorporates the effect of the structural displacements on the interface forces and fluid residuals through the left and right hand sides of the adjoint structural equation. Although this approach bypasses the adjoint mesh motion problem, the structural Jacobian should be modified for the adjoint problem. Therefore, it is not possible to reuse existing self-adjoint structural solvers. A reduced two-field formulation, which is followed by \cite{Heners2017,Stavropoulou2016,Fazzolari2007}, can be achieved by assuming that the FSI solution is invariant with respect to (w.r.t) the fluid interior mesh. In other words, it is assumed that the interface forces and fluid residuals are only a function of the fluid surface boundary mesh, which yields the so-called boundary-based or reduced gradient formulations \citep{Lozano2017,Kavvadias2015}. 

This paper is structured as follows: In Section \ref{prob_form}, we formulate the stationary fluid-structure interaction problem in a partitioned manner. Section \ref{Aeroelastic_sensitivity_analysis} will focus on the partitioning of the adjoint FSI problem using unique sets of Dirichlet and Neumann-type coupling conditions for multidisciplinary objective functionals. Section \ref{numerical_studies} presents two multiphysics frameworks that are used for the assessment of the adjoint formulations and the well-established mapping algorithms. Finally, in Section \ref{conclusion_sec}, we will give main conclusions of this work.

\section{Fluid-structure interaction }
\label{prob_form}
This section starts with the mathematical description of the stationary fluid-structure interaction problem including a continuous form of the governing equations and an appropriate set of steady coupling conditions at the fluid-structure interface. Without loss of generality, the equations and the interface boundary conditions are then discretized and written in discrete residual form. 
It is important to emphasize that all subsequent derivations are independent of discretization method, e.g., finite-element and finite-volume methods. Lastly, the so-called Dirichlet-Neumann partitioned FSI scheme is presented.

\subsection{Continuous fluid-structure interaction problem}
The system under consideration consists of three main parts: fluid domain $ \Omega^{\mathcal{F}}$, structural domain $ \Omega^{\mathcal{S}}$ and wet fluid–structure interface $ \Gamma_{\mathcal{I}}$. The superscripts $\mathcal{F}$ and $\mathcal{S}$ denote the fluid and structure respectively, and is the convention used throughout the paper. As has been usually done and will be pursued here, Eulerian and total Lagrangian approaches are used to describe fluid and structure motions, respectively. Note that the same descriptions have conventionally been used in single-disciplinary solver implementations. A total Lagrangian approach formulates structural governing equations $\boldsymbol{\mathcal{R}}^{\mathcal{S}}$ with respect to the undeformed configuration $X$ while a Eulerian approach
formulates fluid governing equations $\boldsymbol{\mathcal{R}}^{\mathcal{F}}$ with respect to the deformed configuration $x$. In order to couple the governing equations, we require kinematic continuity as well as the equilibrium of interface traction fields at the fluid-structure interface.
Here, we also describe the motion of the fluid domain by structural/pseudo-structural governing equations $\boldsymbol{\mathcal{R}}^{\mathcal{M}}$. Assuming steady-state conditions, the continuous form of the problem can be written as:
\begin{subequations}
\label{eq:conti_FSI}
\begin{align}
\quad  \quad  \; \; \; \boldsymbol{\mathcal{R}}^{\mathcal{F}}\left(\boldsymbol{w}^{\mathcal{F}},\boldsymbol{x}^{\mathcal{F}}\right)& = \boldsymbol{0} \, \quad \quad  \quad \text{in } ^{x}\Omega^{\mathcal{F}}\label{eq:conti_FSI_f}\\
\boldsymbol{\mathcal{R}}^{\mathcal{S}}\left(\boldsymbol{u}^{\mathcal{S}},\boldsymbol{X}^{\mathcal{S}}\right)& = \boldsymbol{0} \quad  \quad \quad \,  \text{in } ^{X}\Omega^{\mathcal{S}}\label{eq:conti_FSI_s}\\
\boldsymbol{\mathcal{R}}^{\mathcal{M}}\left(\boldsymbol{u}^{\mathcal{F}},\boldsymbol{X}^{\mathcal{F}}\right)& = \boldsymbol{0} \quad  \quad \quad \, \text{in } ^{X}\Omega^{\mathcal{F}}\label{eq:conti_FSI_m}
\end{align}
subject to
\begin{align}
\boldsymbol{v}^{\mathcal{F}}_{\Gamma_{\mathcal{I}}} & =  \boldsymbol{0}  \,  \, \quad  \quad   \quad  \text{on } ^{x}\Gamma^{\mathcal{F}}_{\mathcal{I}}\label{eq:conti_FSI_v}\\
\boldsymbol{\sigma}^{\mathcal{F}}_{\Gamma_{\mathcal{I}}} \cdot \boldsymbol{n}^{\mathcal{F}} + \boldsymbol{\sigma}^{\mathcal{S}}_{\Gamma_{\mathcal{I}}} \cdot \boldsymbol{n}^{\mathcal{S}} & =  \boldsymbol{0}  \, \, \quad  \quad   \quad  \text{on } ^{x}\Gamma^{\mathcal{S}}_{\mathcal{I}}\label{eq:conti_FSI_t}\\
\boldsymbol{u}^{\mathcal{S}}_{\Gamma_{\mathcal{I}}} - \boldsymbol{u}^{\mathcal{F}}_{\Gamma_{\mathcal{I}}} & =  \boldsymbol{0}  \, \, \quad  \quad   \quad  \text{on } ^{X}\Gamma_{\mathcal{I}}^{\mathcal{F}}\label{eq:conti_FSI_sd}\\
\boldsymbol{X}^{\mathcal{F}}_{\Gamma_{\mathcal{I}}}+\boldsymbol{u}^{\mathcal{F}}_{\Gamma_{\mathcal{I}}} -\boldsymbol{x}^{\mathcal{F}}_{\Gamma_{\mathcal{I}}}  & =  \boldsymbol{0}  \, \,   \quad  \quad \quad  \text{on } ^{X}\Gamma^{\mathcal{F}}_{\mathcal{I}}\label{eq:conti_FSI_md_b}\\
\boldsymbol{X}^{\mathcal{F}}_{\Omega}+\boldsymbol{u}^{\mathcal{F}}_{\Omega} - \boldsymbol{x}^{\mathcal{F}}_{\Omega} & = \boldsymbol{0}  \, \, \, \quad  \quad \quad  \text{in } ^{X}\Omega^{\mathcal{F}}\label{eq:conti_FSI_md}
\end{align}
\end{subequations}
where the notation $ ^\alpha \left( \cdot \right)^{\beta}_{\gamma}$ is introduced for the sake of clarity; the left superscript $\alpha \in \{x,X\}$ indicates the configuration for evaluation; the right superscript $\beta \in \{\mathcal{F},\mathcal{S},\mathcal{M}\}$ denotes that the variable belongs to fluid or structure or fluid mesh motion; and the subscript $\gamma \in \{\Omega,\Gamma_{(\cdot)}\}$ indicates 
whether the quantity is evaluated inside the domain or on a boundary. $\boldsymbol{x}$ and $\boldsymbol{X}$ refer to the Cartesian coordinates of deformed and undeformed configurations, respectively. The quantity $ \boldsymbol{w}^{\mathcal{F}}$ denotes the state variables of the fluid, typically velocities $ \boldsymbol{v}^{\mathcal{F}}$ with the pressure $p^{\mathcal{F}}$ or the density and the internal energy. The displacement fields $\boldsymbol{u}^{\mathcal{S}}$ and $\boldsymbol{u}^{\mathcal{F}}$ represent the displacements of structure and fluid, respectively. The vector $\boldsymbol{n}$ is the surface unit normal vector and $\boldsymbol{\sigma}$ is the Cauchy stress tensor (i.e., stress measured in the deformed configuration).  

$\boldsymbol{\mathcal{R}}^{\mathcal{F}}$ represents the continuum equations that govern the fluid flow. We describe the motion of the fluid by the full Navier–Stokes compressible equations, from which all the types of governing flow equations can be derived. Defining a conservative variable $\boldsymbol{w}^{\mathcal{F}} = (\rho^{\mathcal{F}}, \rho^{\mathcal{F}} \boldsymbol{v}^{\mathcal{F}}, \rho^{\mathcal{F}}E^{\mathcal{F}})$, their steady-state formulation for a viscous, compressible, Newtonian flow can be written in the following form:

\begin{equation}
\label{e:COMP_NS_function}
\begin{aligned}
\begin{cases}
 \boldsymbol{\mathcal{R}}^{\mathcal{F}} = \nabla_x \cdot \boldsymbol{F}^{c}-\nabla_x \cdot\boldsymbol{F}^{v} =  \boldsymbol{0}     & \text{in } ^{x}\Omega^{\mathcal{F}}\\
  \boldsymbol{v}^{\mathcal{F}}_{\Gamma_{\mathcal{I}}} =  \boldsymbol{0}  & \text{on } ^{x}\Gamma^{\mathcal{F}}_{\mathcal{I}} \\
\boldsymbol{n}^{\mathcal{F}} \cdot \nabla_x T^{\mathcal{F}}_{\Gamma_{\mathcal{I}}} = 0 & \text{on } ^{x}\Gamma^{\mathcal{F}}_{\mathcal{I}} \\
\boldsymbol{w}^{\mathcal{F}}_{\Gamma_\infty} = \bar{\boldsymbol{w}}^{\mathcal{F}}_{\infty} & \text{on } ^{x}\Gamma^{\mathcal{F}}_{\infty} 
  \end{cases}
\end{aligned}
\end{equation}
where the operator $\nabla_x$ denotes the derivatives with respect to the deformed configuration {\itshape x}, over-bar $(\bar{\cdot})$ indicates prescribed value, and the convective fluxes, viscous fluxes are 

\begin{equation}
\begin{aligned}
& \boldsymbol{F}^{c}  =  \begin{cases}
\rho^{\mathcal{F}}\boldsymbol{v}^{\mathcal{F}} \\
\rho^{\mathcal{F}}\boldsymbol{v}^{\mathcal{F}}\otimes \boldsymbol{v}^{\mathcal{F}} + \boldsymbol{I} p^{\mathcal{F}} \\
\rho^{\mathcal{F}} \boldsymbol{v}^{\mathcal{F}} E^{\mathcal{F}} + p^{\mathcal{F}} \boldsymbol{v}^{\mathcal{F}}
\end{cases}\\
& \boldsymbol{F}^{v}  =  \begin{cases}
0 \\
\boldsymbol{\tau}^{\mathcal{F}} \\
\boldsymbol{\tau}^{\mathcal{F}}\cdot\boldsymbol{v}^{\mathcal{F}} + (\mu^{\mathcal{F}}/Pr) C_{p}\nabla_x T^{\mathcal{F}})
\end{cases}
\end{aligned}
\end{equation}
where $\rho^{\mathcal{F}}$ is the fluid density, $ \boldsymbol{v}^{\mathcal{F}}$ represents the flow velocities in all dimensions, $ p^{\mathcal{F}}$ is the physical pressure, $\boldsymbol{I}$ is the identity matrix, $E^{\mathcal{F}}$ is the total energy of the flow per unit mass, $\mu^{\mathcal{F}}$ is the fluid viscosity, $Pr$ is the Prandtl number, $C_{p}$ is the specific heat, $T^{\mathcal{F}}$ is the temperature, and $\boldsymbol{\tau}^{\mathcal{F}}$ is the viscous stress tensor and defined as 
\begin{equation}
\boldsymbol{\tau}^{\mathcal{F}} = \mu^{\mathcal{F}} \left(\nabla_x \boldsymbol{v}^{\mathcal{F}} + \nabla_x (\boldsymbol{v}^{\mathcal{F}})^T - \frac{2}{3} \boldsymbol{I} \left(\nabla_x\cdot \boldsymbol{v}^{\mathcal{F}}\right) \right).
\end{equation}
After having solved the governing flow equations for a given set of boundary conditions, the fluid Cauchy stress tensor reads
\begin{equation}
\label{e:fluid_cauchy_stress}
\boldsymbol{\sigma}^{\mathcal{F}} ~=~  \boldsymbol{I}  p^{\mathcal{F}} - \boldsymbol{\tau}^{\mathcal{F}}.
\end{equation}
Note that in the case of inviscid flow, only the pressure field contributes to the stress tensor, furthermore, the no-slip condition in Eq. \ref{eq:conti_FSI_v} gets modified to the Euler slip condition (i.e. $\boldsymbol{v}^{\mathcal{F}}\cdot \boldsymbol{n}^{\mathcal{F}}=0$). 

Following a total Lagrangian approach, the static and continuous conservation of momentum written in terms of the second Piola-Kirchhoff stress $\boldsymbol{S}^{\mathcal{S}}$ and the Lagrangian coordinates $\boldsymbol{X}$ of structural domain is

\begin{equation}
\label{e:structure}
\begin{aligned}
\begin{cases}
\boldsymbol{\mathcal{R}}^{\mathcal{S}} = \nabla_{X} \cdot \left(\boldsymbol{F}^{\mathcal{S}}\cdot\boldsymbol{S}^{\mathcal{S}}\right)+ \rho^{\mathcal{S}}\boldsymbol{b}^{\mathcal{S}} = \vec{0}  & \text{in } ^{X}\Omega^{\mathcal{S}}\\
\boldsymbol{\sigma}^{\mathcal{S}}_{\Gamma_{\mathcal{I}}} \cdot \boldsymbol{n}^{\mathcal{S}} + \boldsymbol{\sigma}^{\mathcal{F}}_{\Gamma_{\mathcal{I}}} \cdot \boldsymbol{n}^{\mathcal{F}} = \boldsymbol{0} & \text{on } ^{x}\Gamma^{\mathcal{S}}_{\mathcal{I}} \\
\boldsymbol{S}^{\mathcal{S}}_{\Gamma_{\mathcal{N}}} \cdot \boldsymbol{n}^{\mathcal{S}}  = \bar{\boldsymbol{t}}^{\mathcal{S}}_{\Gamma_{\mathcal{N}}} & \text{on } ^{X}\Gamma^{\mathcal{S}}_{\mathcal{N}} \\
\boldsymbol{u}^{\mathcal{S}}_{\Gamma_{\mathcal{D}}} = \bar{\boldsymbol{u}}^{\mathcal{S}}_{\Gamma_{\mathcal{D}}} & \text{on } ^{X}\Gamma^{\mathcal{S}}_{\mathcal{D}}
\end{cases}
\end{aligned}
\end{equation}
here, $\boldsymbol{F}^{\mathcal{S}}=\nabla_{X}\boldsymbol{x}^{\mathcal{S}}$ represents the deformation gradient;  $\rho^{S}$ is the density of the structural domain; $\boldsymbol{b}^{\mathcal{S}}$ is the volumetric body force. Note that $\nabla_{X}$ indicates the spatial gradient operator acting on the undeformed configuration {\itshape X}. The second Piola-Kirchhoff stress tensor $\boldsymbol{S}^{\mathcal{S}}$ is related to the Green-Lagrangian strains via
\begin{equation}
\label{e:strain_stress_rel}
\boldsymbol{S}^{\mathcal{S}} = \boldsymbol{C : E}^{\mathcal{S}}  \quad \text{with} \quad \boldsymbol{E}^{\mathcal{S}} = \frac{1}{2} \left((\boldsymbol{F}^{\mathcal{S}})^{T}\cdot \boldsymbol{F}^{\mathcal{S}}-\boldsymbol{I}\right)
\end{equation}
where $\boldsymbol{C}$ denotes the material tensor. 
Last but not least, the motion of the fluid domain can then be described by

\begin{equation}
\label{e:mesh_motion}
\begin{aligned}
\begin{cases}
\boldsymbol{\mathcal{R}}^{\mathcal{M}} = \nabla_{X} \cdot \boldsymbol{\sigma}^{\mathcal{M}} = \vec{0}  & \text{in } ^{X}\Omega^{\mathcal{F}}\\
\boldsymbol{u}^{\mathcal{F}}_{\Gamma_{\mathcal{I}}} = \boldsymbol{u}^{\mathcal{S}}_{\Gamma_{\mathcal{I}}}  & \text{on } ^{X}\Gamma^{\mathcal{F}}_{\mathcal{I}} \\
\boldsymbol{u}^{\mathcal{F}}_{\Gamma_{\infty}} = \boldsymbol{0} & \text{on } ^{X}\Gamma^{\mathcal{F}}_{\infty}
\end{cases}
\end{aligned}
\end{equation}
Considering the case of pseudo-linear elasticity, $\boldsymbol{\sigma}^{\mathcal{M}}$ is defined as
\begin{equation}
\label{e:linear_strain_stress_rel}
\boldsymbol{\sigma}^{\mathcal{M}} = \lambda \text{tr}\left( \boldsymbol{\epsilon}\left(\boldsymbol{u}^{\mathcal{F}}\right)\right) \boldsymbol{I} + 2 \mu^{\mathcal{M}} \boldsymbol{\epsilon}\left(\boldsymbol{u}^{\mathcal{F}}\right)
\end{equation}
where $\text{tr()}$ is the trace operator, $\lambda$ and $\mu^{\mathcal{M}}$ are the Lame constants, and $\boldsymbol{\epsilon}$ is the strain
tensor:
\begin{equation}
\label{e:linear_strain}
\boldsymbol{\epsilon}\left(\boldsymbol{u}^{\mathcal{F}}\right) = \frac{1}{2}\left( \nabla_{X}\boldsymbol{u}^{\mathcal{F}}+(\nabla_{X}\boldsymbol{u}^{\mathcal{F}})^T\right).
\end{equation}

\subsection{Discrete fluid-structure interaction problem}
In order to solve the explained coupled problem numerically, spatial discretization of the governing equations and all unknown fields is required. Having arbitrarily discretized the fluid and structure domains with $m^F$ and $m^S$ nodes, respectively, the residual form of the stationary FSI problem in Eq. \ref{eq:conti_FSI} reads independently of the spatial discretization scheme as follows: 

\begin{subequations}
\label{eq:Res_Form_FSI_with_BCs}
\begin{align}
& \boldsymbol{r}^{\mathcal{F}}\left(\boldsymbol{w}^{\mathcal{F}},\boldsymbol{x}^{\mathcal{F}}\right) = \boldsymbol{0} \label{eq:Res_Form_FSI_with_BCs_f}\\
\nonumber \\
& \boldsymbol{r}^{\mathcal{S}}\left(\boldsymbol{u}^{\mathcal{S}},\boldsymbol{X}^{\mathcal{S}},\boldsymbol{w}^{\mathcal{F}},\boldsymbol{x}^{\mathcal{F}},\underline{\boldsymbol{X}^{\mathcal{F}}_{\Gamma_{\mathcal{I}}}}\right) = \nonumber\\ & \boldsymbol{f}^{{\mathcal{S}},ext}\left(\boldsymbol{w}^{\mathcal{F}},\boldsymbol{x}^{\mathcal{F}},\underline{\boldsymbol{X}^{\mathcal{F}}_{\Gamma_{\mathcal{I}}}},\underline{\boldsymbol{X}^{\mathcal{S}}_{\Gamma_{\mathcal{I}}}}\right)-\boldsymbol{f}^{{\mathcal{S}},int}\left(\boldsymbol{u}^{\mathcal{S}},\boldsymbol{X}^{\mathcal{S}}\right)  = \boldsymbol{0} \label{Res_Form_FSI_with_BCs_s}\\
\nonumber\\
 &  \boldsymbol{r}^{\mathcal{M}}\left(\boldsymbol{u}^{\mathcal{F}},\boldsymbol{X}^{\mathcal{F}},\boldsymbol{u}^{\mathcal{S}},\underline{\boldsymbol{X}^{\mathcal{S}}_{\Gamma_{\mathcal{I}}}}\right) = \nonumber\\ & \boldsymbol{f}^{\mathcal{M},ext}-\boldsymbol{f}^{\mathcal{M},int}\left(\boldsymbol{u}^{\mathcal{F}},\boldsymbol{X}^{\mathcal{F}}\right) = \boldsymbol{0} \label{Res_Form_FSI_with_BCs_m}
\end{align}
subject to
\begin{align}
\boldsymbol{v}^{\mathcal{F}}_{\Gamma_{\mathcal{I}}} & =  \left[\boldsymbol{0}\right]_{m^{\mathcal{F}}_{\Gamma_{\mathcal{I}}}\times 1}  \,  \, \quad  \text{on } ^{x}\Gamma^{\mathcal{F}}_{\mathcal{I}}\label{Res_Form_FSI_with_BCs_v}\\
\boldsymbol{r}^{\mathcal{S}}_{\Gamma_{\mathcal{I}}} =  \boldsymbol{H}^{\mathcal{F}} \cdot \boldsymbol{f}^{\mathcal{F}}_{\Gamma_\mathcal{I}} - \boldsymbol{f}^{{\mathcal{S}},int}_{\Gamma_{\mathcal{I}}} & =  \left[\boldsymbol{0}\right]_{m^{\mathcal{S}}_{\Gamma_{\mathcal{I}}}\times 1}  \,  \,    \quad  \text{on } ^{X}\Gamma^{\mathcal{S}}_{\mathcal{I}}\label{Res_Form_FSI_with_BCs_t}\\
\boldsymbol{r}^{\mathcal{M}}_{\Gamma_\mathcal{I}} = \boldsymbol{H}^{\mathcal{S}} \cdot \boldsymbol{u}^{\mathcal{S}}_{\Gamma_\mathcal{I}} - \boldsymbol{u}^{\mathcal{F}}_{\Gamma_\mathcal{I}}  & = \left[\boldsymbol{0}\right]_{m^{\mathcal{F}}_{\Gamma_{\mathcal{I}}}\times 1}  \, \,  \quad  \text{on } ^{X}\Gamma_{\mathcal{I}}^{\mathcal{F}}\label{Res_Form_FSI_with_BCs_sd}\\
\boldsymbol{X}^{\mathcal{F}}_{\Gamma_{\mathcal{I}}}+\boldsymbol{u}^{\mathcal{F}}_{\Gamma_{\mathcal{I}}} -\boldsymbol{x}^{\mathcal{F}}_{\Gamma_{\mathcal{I}}}  & =  \left[\boldsymbol{0}\right]_{m^{\mathcal{F}}_{\Gamma_{\mathcal{I}}}\times 1}  \, \,   \quad   \text{on } ^{X}\Gamma^{\mathcal{F}}_{\mathcal{I}}\label{Res_Form_FSI_with_BCs_md_bc}\\
 \boldsymbol{X}^{\mathcal{F}}_{\Omega}+\boldsymbol{u}^{\mathcal{F}}_{\Omega} -\boldsymbol{x}^{\mathcal{F}}_{\Omega}  & =  \left[\boldsymbol{0}\right]_{m^{\mathcal{F}}_{\Omega}\times 1}\, \, \, \, \, \, \, \, \, \, \, \text{in } ^{X}\Omega^{\mathcal{F}}\label{Res_Form_FSI_with_BCs_md}
\end{align}
\end{subequations}
where $\boldsymbol{r}^{\mathcal{F}}$, $\boldsymbol{r}^{\mathcal{S}}$, and $\boldsymbol{r}^{\mathcal{M}}$ are the full residual vectors (including the internal and the boundary unknowns of the PDEs) of the fluid, the structure, and the mesh motion, respectively. $\boldsymbol{x}^{\mathcal{F}}$ represents the nodal coordinates of the fluid mesh in the deformed configuration, while $\boldsymbol{X}^{\mathcal{S}}$ and $\boldsymbol{X}^{\mathcal{F}}$ are the nodal coordinates of the structural and fluid meshes in the undeformed configuration, respectively. $\boldsymbol{f}^{\mathcal{S},int}$ and $\boldsymbol{f}^{\mathcal{S},ext}$ are the vector of internal forces and the vector of external forces in the structure, respectively, while $\boldsymbol{f}^{\mathcal{M},int}$ and $\boldsymbol{f}^{\mathcal{M},ext}$ are the same terms for the mesh motion. $\boldsymbol{f}^{\mathcal{M},ext}$ is generally zero. Note that, although two coupled domains at the interface have matching geometries (i.e. $\Gamma^{\mathcal{S}}_{\mathcal{I}} = \Gamma^{\mathcal{F}}_{\mathcal{I}}$ ), the meshes at the fluid–structure interface usually do not node-to-node match due to the different mesh requirements for the flow and structure (i.e. $ \boldsymbol{X}^{\mathcal{F}}_{\Gamma_{\mathcal{I}}}\neq \boldsymbol{X}^{\mathcal{S}}_{\Gamma_{\mathcal{I}}}$ , $m^{\mathcal{F}}_{\Gamma_{\mathcal{I}}} \neq m^{\mathcal{S}}_{\Gamma_{\mathcal{I}}} $ ). Therefore, when applying the coupling conditions to non-matching meshes, mapping is needed before transferring the information. Here, $\boldsymbol{H}^{\mathcal{F}}(\boldsymbol{X}^{\mathcal{S}}_{\Gamma_{\mathcal{I}}},\boldsymbol{X}^{\mathcal{F}}_{\Gamma_{\mathcal{I}}})$ and $\boldsymbol{H}^{\mathcal{S}}(\boldsymbol{X}^{\mathcal{F}}_{\Gamma_{\mathcal{I}}},\boldsymbol{X}^{\mathcal{S}}_{\Gamma_{\mathcal{I}}})$ are defined for the transfer from the fluid to the structure mesh and from the structure mesh to the fluid mesh. Terms and dependencies arising due to the non-matching meshes are underlined as the convention throughout the paper. With matching fluid-structure interface meshes, the mapping matrices reduce to identity matrices and the underlined dependencies in Eq. \ref{eq:Res_Form_FSI_with_BCs} vanish. Although development and assessment of mapping algorithms are not in the scope of this paper, we investigate the difference in accuracy between them \citep{farhat1998load,DeBoer2008,Wang2016papaer}. It is worth noting that in the definition of the interface dynamic continuity (Eq. \ref{Res_Form_FSI_with_BCs_t}), the following identity from continuum mechanics \citep{belytschko2013nonlinear} is used: 
\begin{equation}
\label{eq:law_of_continuum_mechanics} 
\boldsymbol{S}^{\mathcal{S}}_{\Gamma_{\mathcal{I}}} \cdot \boldsymbol{n}^{\mathcal{S}} ~ d\Gamma^{\mathcal{S}}_{\mathcal{I}}
= \boldsymbol{\sigma}^{\mathcal{S}}_{\Gamma_{\mathcal{I}}} \cdot \boldsymbol{n}^{\mathcal{S}} ~ 
d\Gamma^{\mathcal{S}}_{\mathcal{I}} = -\boldsymbol{\sigma}^{\mathcal{F}}_{\Gamma_{\mathcal{I}}} \cdot \boldsymbol{n}^{\mathcal{F}} ~ 
d\Gamma^{\mathcal{F}}_{\mathcal{I}}.
\end{equation}

\subsection{Partitioned FSI}
\label{Partitioned_FSI}
The FSI problem stated continuously in Eq. \ref{eq:conti_FSI} and discretely in Eq. \ref{eq:Res_Form_FSI_with_BCs} constitutes a coupled set of non-linear and linear subproblems that can be solved separately and iteratively until the interface conditions, the equilibrium of tractions and kinematic continuity, are satisfied. This results in the so-called Gauss-Seidel fixed-point iterations for a strongly coupled partitioned fluid-structure interaction. Among the partitioned coupling techniques for FSI \citep{Badia2008}, we use the so-called Dirichlet–Neumann partitioned procedure which is by far the most widely used strategy, both for simplicity and because of wide range of applicability to single-disciplinary solvers. This technique treats the fluid domain as the Dirichlet partition, i.e. it takes the prescribed interface displacements as the Dirichlet boundary condition for the mesh motion problem, and the structure domain as the Neumann partition loaded with interface fluid forces. 

With the Dirichlet–Neumann partitioned procedure, we break down the stationary FSI problem into the fluid, the structure and the mesh motion subproblems which are treated by black-box solvers as 
\begin{subequations}
\label{eq:Partitioned_FSI_black_box}
\begin{align}
\quad \quad \quad \quad \quad \quad 
\boldsymbol{f}^{\mathcal{F}}_{\Gamma_{\mathcal{I}}} & = \mathscr{F}(\boldsymbol{x}^{\mathcal{F}})
\label{eq:Partitioned_FSI_black_box_f}\\
\boldsymbol{u}^{\mathcal{S}}_{\Gamma_{\mathcal{I}}} & = \mathscr{S}(\boldsymbol{f}^{\mathcal{S},ext}_{\Gamma_{\mathcal{I}}})
\label{eq:Partitioned_FSI_black_box_s}\\
\boldsymbol{x}^{\mathcal{F}} & = \mathscr{M}(\boldsymbol{u}^{\mathcal{F}}_{\Gamma_{\mathcal{I}}})\label{eq:Partitioned_FSI_black_box_m}.
\end{align}
\end{subequations}
In the equation, $\mathscr{F}$ represents the fluid solver that takes the new position of the fluid mesh $\boldsymbol{x}^{\mathcal{F}}$ as input and outputs the interface load $\boldsymbol{f}^{\mathcal{F}}_{\Gamma_{\mathcal{I}}}$ (nodal forces or tractions), $\mathscr{S}$ represents the structure solver that takes $\boldsymbol{f}^{\mathcal{S},ext}_{\Gamma_{\mathcal{I}}}$ as input and outputs $\boldsymbol{u}^{\mathcal{S}}_{\Gamma_{\mathcal{I}}}$, and the mesh motion solver $\mathscr{M}$ which outputs the deformed fluid mesh $\boldsymbol{x}^{F}$ according to $\boldsymbol{u}^{\mathcal{F}}_{\Gamma_{\mathcal{I}}}$. 
Algorithm \ref{steady_state_FSI_algorithm} details the Gauss-Seidel algorithm for a stationary FSI problem with arbitrary non-matching interface meshes. In this algorithm, $n$ denotes the current iteration, and $\boldsymbol{\hat{u}}^{\mathcal{S}}_{\Gamma_{\mathcal{I}}}$ is the relaxed interface displacements. Due to the simplicity of implementation and efficiency, the Aitken relaxation is chosen as the default relaxation scheme in this work.

\begin{algorithm}
\caption{Dirichlet-Neumann partitioned FSI work flow}
\label{steady_state_FSI_algorithm}
\begin{algorithmic}[1]

	\Statex  {\color{arsenic}//initialize the mapping matrices between $\boldsymbol{X}^{\mathcal{F}}_{\Gamma_{\mathcal{I}}} $ and $\boldsymbol{X}^{\mathcal{S}}_{\Gamma_{\mathcal{I}}} $}
    \State $\boldsymbol{H}^{\mathcal{F}}$, $\boldsymbol{H}^{\mathcal{S}}$ 
    \State $n=1$ 
    	\Statex  {\color{arsenic}//initialize interface displacements}
    \State $ _n\boldsymbol{\hat{u}}^{\mathcal{S}}_{\Gamma_{\mathcal{I}}} = \boldsymbol{0}$ 
    \Statex  {\color{arsenic}//FSI strong coupling loop}
    \While{$\left\Vert{_n\boldsymbol{\delta}}\right\Vert_{2} > \varepsilon$}
    \Statex  {\color{arsenic} // map the relaxed interface displacements}
    \State   $ _{n}\boldsymbol{\hat{u}}^{\mathcal{F}}_{\Gamma_{\mathcal{I}}} = \boldsymbol{H}^{\mathcal{S}} \cdot ~ _n \boldsymbol{\hat{u}}^{\mathcal{S}}_{\Gamma_{\mathcal{I}}}  $
    \Statex  {\color{arsenic} // solve mesh motion problem}
    \State   $ _n\boldsymbol{x}^{\mathcal{F}}= \mathscr{M}(_{n}\boldsymbol{\hat{u}}^{\mathcal{F}}_{\Gamma_{\mathcal{I}}} )$
    \Statex  {\color{arsenic} // solve fluid problem }
    \State   $_n\boldsymbol{f}^{\mathcal{F}}_{\Gamma_{\mathcal{I}}} =\mathscr{F}(_n\boldsymbol{x}^{F})$
    \Statex {\color{arsenic}// map the interface forces}
    \State $ _n\boldsymbol{f}^{\mathcal{S}}_{\Gamma_{\mathcal{I}}} = \boldsymbol{H}^{\mathcal{F}} \cdot ~ _n\boldsymbol{f}^{\mathcal{F}}_{\Gamma_{\mathcal{I}}} $   
    \Statex {\color{arsenic}// solve structure problem}
    \State  $  _n\boldsymbol{u}^{\mathcal{S}}_{\Gamma_{\mathcal{I}}} = \mathscr{S}(_n\boldsymbol{f}^{\mathcal{S}}_{\Gamma_{\mathcal{I}}})$
    \Statex {\color{arsenic}// compute interface displacement residuals}
    \State  $  _n\boldsymbol{\delta} = {}_n\boldsymbol{u}^{\mathcal{S}}_{\Gamma_{\mathcal{I}}} - {}_{n-1}\boldsymbol{u}^{\mathcal{S}}_{\Gamma_{\mathcal{I}}}$   
    \State \textbf{compute ${}_{n+1}\boldsymbol{\hat{u}}^{\mathcal{S}}_{\Gamma_{\mathcal{I}}}$ based on $\{{}_1\boldsymbol{\delta},{}_2\boldsymbol{\delta},\cdot \cdot \cdot,{}_n\boldsymbol{\delta}\}$ \ and $\{{}_1\boldsymbol{u}^{\mathcal{S}}_{\Gamma_{\mathcal{I}}},{}_2\boldsymbol{u}^{\mathcal{S}}_{\Gamma_{\mathcal{I}}},\cdot\cdot\cdot, {}_n\boldsymbol{u}^{\mathcal{S}}_{\Gamma_{\mathcal{I}}}\}$ (relaxation, etc.)} 
    \State $n=n+1$
    \EndWhile

\end{algorithmic}
\end{algorithm}
\section{Multidisciplinary adjoint-based shape sensitivity analysis }
\label{Aeroelastic_sensitivity_analysis}
Having obtained the equilibrium state of a static FSI system, we formulate the multidisciplinary shape sensitivity analysis as follows: \newline \newline
\textit{We seek to compute the gradients of a multi-objective and multi-disciplinary target function $\tilde{J}$, which is defined as a function of fluid state variables $\boldsymbol{w}^{\mathcal{F}}$ on the deformed configuration or structural state variables $\boldsymbol{u}^{\mathcal{S}}$ on the undeformed configuration or both, with respect to shape design variables, that specify the undeformed geometry of the design surface, e.g the interface.}\newline \newline
The shape optimization problem corresponding to the shape sensitivity analysis problem of interest can be defined mathematically as:

\begin{equation}
\label{e:opt_problem}
\begin{aligned}
& \underset{\boldsymbol{X}_{\mathcal{D}}}{\text{min}}
& & \tilde{J} = \alpha^{i} J^{i}, i \in \{\mathcal{F},\mathcal{S},\mathcal{I}\}, \alpha^{i} \in \mathbb{R}\\
& \text{subject to} & &  \\
& & &  \boldsymbol{r}^{\mathcal{F}}\left(\boldsymbol{w}^{\mathcal{F}},\boldsymbol{x}^{\mathcal{F}}\right)  =  \boldsymbol{0}  \\ 
& & & \boldsymbol{r}^{\mathcal{S}}\left(\boldsymbol{u}^{\mathcal{S}},\boldsymbol{X}^{\mathcal{S}},\boldsymbol{w}^{\mathcal{F}},\boldsymbol{x}^{\mathcal{F}},\underline{\boldsymbol{X}^{\mathcal{F}}_{\Gamma_{\mathcal{I}}}}\right) =  \boldsymbol{0} \\
& & & \boldsymbol{r}^{\mathcal{M}}\left(\boldsymbol{u}^{\mathcal{F}},\boldsymbol{X}^{\mathcal{F}},\boldsymbol{u}^{\mathcal{S}},\underline{\boldsymbol{X}^{\mathcal{S}}_{\Gamma_{\mathcal{I}}}}\right) =  \boldsymbol{0}  
\end{aligned}
\end{equation}
where $\tilde{J}$ is the weighted sum of the objectives $J^i$, and $\boldsymbol{X}_{\mathcal{D}} \in \{ \boldsymbol{X}^{\mathcal{F}}_{\Gamma_{\mathcal{I}}},\boldsymbol{X}^{\mathcal{F}}_{\Gamma_{\infty}},\boldsymbol{X}^{\mathcal{S}}_{\Gamma_{\mathcal{I}}},\boldsymbol{X}^{\mathcal{S}}_{\Gamma_{\mathcal{D}}},\boldsymbol{X}^{\mathcal{S}}_{\Gamma_{\mathcal{N}}}\}$ $\boldsymbol{X}_{\mathcal{D},i}$ $\in \mathbb{R}^{3}$ denotes the coordinate vector of the design surface mesh in the undeformed configuration. Note that, $J^{\mathcal{F}}(\boldsymbol{w}^{\mathcal{F}},\boldsymbol{x}^{\mathcal{F}})$ and $J^{\mathcal{S}}(\boldsymbol{u}^{\mathcal{S}},\boldsymbol{X}^{\mathcal{S}})$ respectively represent typical fluid and structure objective functions that can be found in single-disciplinary adjoint solvers. On the other hand, $J^{\mathcal{I}}(\boldsymbol{w}^{\mathcal{F}},\boldsymbol{x}^{\mathcal{F}},\boldsymbol{u}^{\mathcal{S}},\underline{\boldsymbol{X}^{\mathcal{F}}_{\Gamma_{\mathcal{I}}}},\underline{\boldsymbol{X}^{\mathcal{S}}_{\Gamma_{\mathcal{I}}}})$ is only defined on the interface and explicitly depends on all FSI state variables. A good example of such an objective function is the interface energy which is defined as:
\begin{equation}
\label{eq:interface_energy}
\begin{aligned}
J^{\mathcal{I}} = (\boldsymbol{H}^{\mathcal{S}} \cdot \boldsymbol{u}^{\mathcal{S}}_{\Gamma_\mathcal{I}})^T \cdot  \boldsymbol{f}^{\mathcal{F}}_{\Gamma_\mathcal{I}} = (\boldsymbol{u}^{\mathcal{S}}_{\Gamma_\mathcal{I}})^T\cdot(\boldsymbol{H}^{\mathcal{F}} \cdot \boldsymbol{f}^{\mathcal{F}}_{\Gamma_\mathcal{I}}).
\end{aligned}
\end{equation}     
We note that this expression enforces the conservation of the interface energy and results in the following identity:
\begin{equation}
\label{eq:interface_energy2}
\boldsymbol{H}^{\mathcal{F}} = (\boldsymbol{H}^{\mathcal{S}})^T 
\end{equation}  
where superscript $T$ denotes the transpose operator. Note that, in the case of matching meshes at the interface, the mapping matrices reduce to identity matrices and the underlined dependencies vanish.

In a manner consistent with the primal FSI problem, we define a Lagrange function that augments the objective function to incorporate the state constraints (Eq. \ref{e:opt_problem}):

\begin{equation}
\label{e:partitioned_FSI_Lagrange_Function}
\begin{aligned}
\mathcal{L}(\boldsymbol{w}^{\mathcal{F}},\boldsymbol{u}^{\mathcal{F}},\boldsymbol{u}^{\mathcal{S}},\boldsymbol{x}^{\mathcal{F}},\boldsymbol{X}^{\mathcal{F}},\boldsymbol{X}^{\mathcal{S}}) =   \tilde{J} + (\boldsymbol{\Psi}^{i})^T \cdot \boldsymbol{r}^{i}
\end{aligned}
\end{equation}
where $i \in \{\mathcal{F},\mathcal{S},\mathcal{M}\}$, $\boldsymbol{\Psi} = \left[\boldsymbol{\Psi}^{\mathcal{F}},\boldsymbol{\Psi}^{\mathcal{S}},\boldsymbol{\Psi}^{\mathcal{M}}\right]$ is the vector of the adjoint variables (Lagrange multipliers) associated with the complete residual vector ($\boldsymbol{r}=\left[\boldsymbol{r}^{\mathcal{F}},\boldsymbol{r}^{\mathcal{S}},\boldsymbol{r}^{\mathcal{M}}\right]$). Exploiting the chain rule of differentiation and the kinematic conditions in Eqs. \ref{Res_Form_FSI_with_BCs_md_bc} and \ref{Res_Form_FSI_with_BCs_md}, the total variation of $\mathcal{L}$ with respect to the undeformed shape of the design surface hence reads:

\begin{equation}
\label{e:Derivative_of_FSI_Lagrange_Function}
\begin{aligned}
& \frac{d\mathcal{L}}{d\boldsymbol{X}_{\mathcal{D}}}  =  \\ & \frac{\partial \mathcal{L}}{\partial \boldsymbol{w}^{\mathcal{F}}} \cdot  \frac{d \boldsymbol{w}^{\mathcal{F}}}{d\boldsymbol{X}_{\mathcal{D}}} + \left( \frac{\partial \mathcal{L}}{\partial \boldsymbol{u}^{\mathcal{F}}}+\frac{\partial \mathcal{L}}{\partial \boldsymbol{x}^{\mathcal{F}}}\right) \cdot  \frac{d \boldsymbol{u}^{\mathcal{F}}}{d\boldsymbol{X}_{\mathcal{D}}} + \frac{\partial \mathcal{L}}{\partial \boldsymbol{u}^{\mathcal{S}}} \cdot  \frac{d \boldsymbol{u}^{\mathcal{S}}}{d\boldsymbol{X}_{\mathcal{D}}} \\
& + \frac{\partial \mathcal{L}}{\partial \boldsymbol{x}^{\mathcal{F}}} \cdot  \frac{d \boldsymbol{X}^{\mathcal{F}}}{d\boldsymbol{X}_{\mathcal{D}}} + \frac{\partial \mathcal{L}}{\partial \boldsymbol{X}^{\mathcal{F}}} \cdot  \frac{d \boldsymbol{X}^{\mathcal{F}}}{d\boldsymbol{X}_{\mathcal{D}}} + \frac{\partial \mathcal{L}}{\partial \boldsymbol{X}^{\mathcal{S}}} \cdot  \frac{d \boldsymbol{X}^{\mathcal{S}}}{d\boldsymbol{X}_{\mathcal{D}}}.
\end{aligned}
\end{equation}
While the terms multiplying $\frac{d \boldsymbol{w}^{\mathcal{F}}}{d\boldsymbol{X}_{\mathcal{D}}}$, $\frac{d \boldsymbol{u}^{\mathcal{F}}}{d\boldsymbol{X}_{\mathcal{D}}}$ and $\frac{d \boldsymbol{u}^{\mathcal{S}}}{d\boldsymbol{X}_{\mathcal{D}}}$ are eliminated respectively by satisfying the adjoint fluid problem, the adjoint structure problem and the adjoint mesh motion problem, the expressions in the last line give rise to the coupled shape gradients. Special attention must be paid to $\frac{\partial \mathcal{L}}{\partial \boldsymbol{x}^{\mathcal{F}}}$ which is a partial derivative of the Lagrange functional w.r.t the deformed fluid mesh, including both the internal and boundary nodes. Also observe that this term contributes not only to the coupled adjoint mesh motion problem but also to the coupled shape sensitivities.

Expanding each partial shape derivative in Eq. \ref{e:Derivative_of_FSI_Lagrange_Function} results to  
\begin{equation}
\label{e:Total_FSI_Sensitivity_Equation}
\begin{aligned}
& \frac{d\mathcal{L}}{d\boldsymbol{X}_{\mathcal{D}}} = \\ &
 \left(\frac{\partial \tilde{J}}{\partial \boldsymbol{x}^{\mathcal{F}}} + (\boldsymbol{\Psi}^{\mathcal{F}})^{T} \cdot \frac{\partial \boldsymbol{r}^{\mathcal{F}}}{\partial \boldsymbol{x}^{\mathcal{F}}} + (\boldsymbol{\Psi}^{\mathcal{S}}_{\Gamma_\mathcal{I}})^{T} \cdot \boldsymbol{H}^{\mathcal{F}} \cdot \frac{\partial \boldsymbol{f}^{\mathcal{F}}_{\Gamma_\mathcal{I}}}{\partial \boldsymbol{x}^{\mathcal{F}}}\right)\cdot \frac{d \boldsymbol{X}^{\mathcal{F}}}{d\boldsymbol{X}_{\mathcal{D}}} \\
& + \left( (\boldsymbol{\Psi}^{\mathcal{M}})^{T} \cdot \frac{\partial \boldsymbol{r}^{\mathcal{M}}}{\partial \boldsymbol{X}^{\mathcal{F}}}\right) \cdot \frac{d \boldsymbol{X}^{\mathcal{F}}}{d\boldsymbol{X}_{\mathcal{D}}}  + \\ 
&  \left(\frac{\partial \tilde{J}}{\partial \boldsymbol{X}^{\mathcal{S}}} + (\boldsymbol{\Psi}^{\mathcal{S}})^{T} \cdot \frac{\partial \boldsymbol{r}^{\mathcal{S}}}{\partial \boldsymbol{X}^{\mathcal{S}}}\right)\cdot \frac{d \boldsymbol{X}^{\mathcal{S}}}{d\boldsymbol{X}_{\mathcal{D}}}.
\end{aligned}
\end{equation}
This is a valuable shape sensitivity equation for the FSI problem since it clearly states which shape sensitivities should be computed by each discipline and in which configuration. Precisely, the first parentheses in the above equation contain partial shape derivatives which can basically be computed by an adjoint fluid solver in the deformed configuration, whereas the second and third parentheses can be computed by an adjoint structural solver in the undeformed fluid and structure configurations, respectively. Remember, a structural/pseudo-structural model is used here for the fluid mesh motion problem.

In the following subsections, we discuss the coupled adjoint systems and, subsequently, the analysis of the coupled shape sensitivity equations in a fully partitioned way.

\subsection{Coupled adjoint fluid problem}
Taking into account the above-mentioned findings and the explicit dependency of the residual vectors and the general objective function $J$ on the fluid state (see Eq. \ref{e:opt_problem}), the coupled adjoint system and shape sensitivities associated with the fluid read as follows: 
\begin{subequations}
\label{e:adjoint_fluid}
\begin{align}
&\frac{\partial \mathcal{L}}{\partial \boldsymbol{w}^{\mathcal{F}}} = \nonumber \\ & \quad  \quad  \frac{\partial \tilde{J}}{\partial \boldsymbol{w}^{\mathcal{F}}} + (\boldsymbol{\Psi}^{\mathcal{F}})^T \cdot \frac{\partial \boldsymbol{r}^{\mathcal{F}}}{\partial \boldsymbol{w}^{\mathcal{F}}} + (\boldsymbol{\Psi}^{\mathcal{S}}_{\Gamma_\mathcal{I}})^{T} \cdot \boldsymbol{H}^{\mathcal{F}} \cdot \frac{\partial \boldsymbol{f}^{\mathcal{F}}_{\Gamma_\mathcal{I}}}{\partial \boldsymbol{w}^{\mathcal{F}}} = \boldsymbol{0}^T \\
& \frac{\partial \mathcal{L}}{\partial \boldsymbol{x}^{\mathcal{F}}} = \nonumber \\
& \quad  \quad \frac{\partial \tilde{J}}{\partial \boldsymbol{x}^{\mathcal{F}}} + (\boldsymbol{\Psi}^{\mathcal{F}})^{T} \cdot \frac{\partial \boldsymbol{r}^{\mathcal{F}}}{\partial \boldsymbol{x}^{\mathcal{F}}} + (\boldsymbol{\Psi}^{\mathcal{S}}_{\Gamma_\mathcal{I}})^{T} \cdot \boldsymbol{H}^{\mathcal{F}} \cdot \frac{\partial \boldsymbol{f}^{\mathcal{F}}_{\Gamma_\mathcal{I}}}{\partial \boldsymbol{x}^{\mathcal{F}}}.
\end{align}
\end{subequations}
Intuitively, one can define the following fluid shape optimization problem whose adjoint system and shape sensitivities are equivalent to Eqs. \ref{e:adjoint_fluid}:

\begin{equation}
\label{e:fluid_coupled_opt_problem}
\begin{aligned}
& \underset{\boldsymbol{x}^{\mathcal{F}}}{\text{min}}
& & \tilde{J}^{\mathcal{F}} = \alpha^{\mathcal{F}} J^{\mathcal{F}} + \alpha^{\mathcal{I}} J^{\mathcal{I}} + J^{\mathcal{F},a} \\
& \text{subject to} & &  \\
& & &  \boldsymbol{r}^{\mathcal{F}}\left(\boldsymbol{w}^{\mathcal{F}},\boldsymbol{x}^{\mathcal{F}}\right)  =  \boldsymbol{0}  \\
& \text{in which} & &  \\
& & & J^{\mathcal{F},a} = \boldsymbol{d}^T \cdot  \boldsymbol{f}^{\mathcal{F}}_{\Gamma_\mathcal{I}}; \quad \boldsymbol{d} = (\boldsymbol{H}^{\mathcal{F}})^T \cdot \boldsymbol{\Psi}^{\mathcal{S}}_{\Gamma_\mathcal{I}}
\end{aligned}
\end{equation}
where $\tilde{J}^{\mathcal{F}}$ is the weighted sum of fluid-dependent functions. $J^{\mathcal{F},a}$ is an auxiliary objective function which arises from the interaction with structure and vanishes identically for an uncoupled fluid system. Analogously to force-based objective functionals (like drag or lift), auxiliary function $J^{\mathcal{F},a}$ projects the interface force vector $\boldsymbol{f}^{\mathcal{F}}_{\Gamma_\mathcal{I}}$ onto the so-called force projection vector $\boldsymbol{d}$. We note that, in contrast to typical force-based objective functions for fluids, the force projection vector of the auxiliary objective function is spatially varying over the interface and it is computed from the interface adjoint displacements which are mapped from the structure.

This above interpretation of coupling between the adjoint fluid problem and the adjoint displacements has been partly inspired by \cite{Fazzolari2007}, where a continuous adjoint formulation for the Euler equations coupled with linear elasticity is presented. 

One can show easily that the adjoint system and shape sensitivities of Eqs. \ref{e:adjoint_fluid} are definitely equal to those of Eq.  \ref{e:fluid_coupled_opt_problem}, by defining the following Lagrange function:
\begin{equation}
\label{e:fluid_lagrange_func}
\mathcal{L}^{\mathcal{F}}(\boldsymbol{w}^{\mathcal{F}},\boldsymbol{x}^{\mathcal{F}}) = \tilde{J}^{\mathcal{F}} + (\boldsymbol{\Psi}^{\mathcal{F}})^{T} \cdot \boldsymbol{r}^{\mathcal{F}}.
\end{equation}
The first order optimality condition for the Lagrange function entails the following identities:
\begin{subequations}
\label{e:after_fluid_lagrange_func}
\begin{align}
\frac{\partial \mathcal{L}}{\partial \boldsymbol{w}^{\mathcal{F}}} & = \frac{\partial \mathcal{L}^{\mathcal{F}}}{\partial \boldsymbol{w}^{\mathcal{F}}} = \boldsymbol{0}^T\\
\frac{\partial \mathcal{L}}{\partial \boldsymbol{x}^{\mathcal{F}}} & = \frac{\partial \mathcal{L}^{\mathcal{F}}}{\partial \boldsymbol{x}^{\mathcal{F}}}.
\end{align}
\end{subequations}

The coupled fluid shape optimization problem presented in Eq. \ref{e:fluid_coupled_opt_problem} is a straightforward multi-objective adjoint optimization for fluids, however, some remarks are given here:  

\begin{remark}
The adjoint fluid solver should be capable of handling a multi-objective shape sensitivity analysis using a single adjoint solution. Otherwise, interaction between objective functions in the adjoint analysis is neglected. 
\end{remark}

\begin{remark}
The adjoint fluid solver is required to accept a non-uniform projection vector for the force-based objective functional. In the majority of derivations and implementations for fluid adjoint shape sensitivity analysis, there is no assumption or condition on the spatial uniformity of the force projection vector. Therefore, there is no need for the extra work in single-disciplinary adjoint fluid solvers to account for the adjoint coupling through the auxiliary objective function.  
\end{remark}

\begin{remark}
If a force-based objective functional is defined on the interface, one can combine the auxiliary objective function $J^{\mathcal{F},a}$ and the objective functional into a single force-based objective functional by summing up the respective force projection vectors. For example, in the case of interface energy, the fluid multi-objective functional in Eq. \ref{e:fluid_coupled_opt_problem} reads:
\begin{equation}
	\tilde{J}^{\mathcal{F}} = (\boldsymbol{d}^*)^T\cdot \boldsymbol{f}^{\mathcal{F}}_{\Gamma_\mathcal{I}}; \quad \boldsymbol{d}^* =  \alpha^{\mathcal{I}} \boldsymbol{H}^{\mathcal{S}} \cdot \boldsymbol{u}^{\mathcal{S}}_{\Gamma_\mathcal{I}} + (\boldsymbol{H}^{\mathcal{F}})^T \cdot \boldsymbol{\Psi}^{\mathcal{S}}_{\Gamma_\mathcal{I}}.
\end{equation}      
\end{remark}

\subsection{Coupled adjoint mesh motion problem}
Due to the full linearization of the FSI problem, which is referred to as the three-field-based formulation, the coupled adjoint system and shape sensitivities of the mesh motion problem appear in Eq. \ref{e:Derivative_of_FSI_Lagrange_Function}. Collecting the terms associated with the variation of the fluid displacement field and the fluid mesh in the deformed and undeformed configurations, results in the following coupled adjoint system and shape sensitivities:

\begin{subequations}
\label{e:adjoint_mesh_motion}
\begin{align}
&\frac{\partial \mathcal{L}}{\partial \boldsymbol{u}^{\mathcal{F}}}+\frac{\partial \mathcal{L}}{\partial \boldsymbol{x}^{\mathcal{F}}} = -(\boldsymbol{\Psi}^{\mathcal{M}})^T \cdot \frac{\partial \boldsymbol{f}^{\mathcal{M},int}}{\partial \boldsymbol{u}^{\mathcal{F}}} + \frac{\partial \mathcal{L}^{\mathcal{F}}}{\partial \boldsymbol{x}^{\mathcal{F}}} =\boldsymbol{0}^T \\
& \frac{\partial \mathcal{L}}{\partial \boldsymbol{X}^{\mathcal{F}}} = -(\boldsymbol{\Psi}^{\mathcal{M}})^T \cdot \frac{\partial \boldsymbol{f}^{\mathcal{M},int}}{\partial \boldsymbol{X}^{\mathcal{F}}} + \nonumber \\ &  \underline{\alpha^{\mathcal{I}} \frac{\partial  J^{\mathcal{I}}}{\partial \boldsymbol{X}^{\mathcal{F}}}} + \underline{(\boldsymbol{\Psi}^{\mathcal{M}}_{\Gamma_\mathcal{I}})^T \cdot \frac{\partial \boldsymbol{r}^{\mathcal{M}}_{\Gamma_\mathcal{I}}}{\partial \boldsymbol{X}^{\mathcal{F}}}} + \underline{(\boldsymbol{\Psi}^{\mathcal{S}}_{\Gamma_\mathcal{I}})^T \cdot \frac{\partial \boldsymbol{f}^{\mathcal{S},ext}_{\Gamma_\mathcal{I}}}{\partial \boldsymbol{X}^{\mathcal{F}}} }
\end{align}
\end{subequations}
where $\frac{\partial \boldsymbol{f}^{\mathcal{M},ext}}{\partial \boldsymbol{u}^{\mathcal{F}}} = \boldsymbol{0} $ is applied due to full Dirichlet boundary conditions for the mesh motion. Another remark is that the underlined terms arise due to the dependency of the mapping matrices/operations on non-matching interface meshes. 

A closer look reveals that the equations above are very similar to the equations resulting from the adjoint-based shape sensitivity analysis for the strain energy of structures under body forces. Therefore, we can formulate the following pseudo optimization problem to efficiently compute the un-underlined terms in Eqs. \ref{e:adjoint_mesh_motion} using a single-disciplinary adjoint structural solver:   
\begin{equation}
\label{e:mesh_motion_coupled_opt_problem}
\begin{aligned}
& \underset{\boldsymbol{X}^{\mathcal{F}}}{\text{min}}
& & \tilde{J}^{\mathcal{M}} = \tilde{J}^{\mathcal{M},a}  \\
& \text{subject to} & &  \\
& & &  \boldsymbol{r}^{\mathcal{M}}\left(\boldsymbol{u}^{\mathcal{F}},\boldsymbol{X}^{\mathcal{F}}\right)  =  \boldsymbol{0} \\ 
& \text{in which} & &  \\
& & & \tilde{J}^{\mathcal{M},a}  = (\boldsymbol{f}^{\mathcal{M},a})^T \cdot \boldsymbol{u}^{\mathcal{F}} \\
& & & \boldsymbol{f}^{\mathcal{M},a} = (\frac{\partial \mathcal{L}^{\mathcal{F}}}{\partial \boldsymbol{x}^{\mathcal{F}}})^T 
\end{aligned}
\end{equation}
where $J^{\mathcal{M},a}$ is an auxiliary objective function which arises from the interaction with the fluid mesh and it can be interpreted as a linear strain energy which is caused by the adjoint body force $\boldsymbol{f}^{\mathcal{M},a}$.
Subsequently, the Lagrange function reads as follows:
\begin{equation}
\label{e:mesh_motion_lagrange_func}
\mathcal{L}^{\mathcal{M}}(\boldsymbol{u}^{\mathcal{F}},\boldsymbol{X}^{\mathcal{F}}) = \tilde{J}^{\mathcal{M},a} + (\boldsymbol{\Psi}^{\mathcal{M}})^{T} \cdot \boldsymbol{r}^{\mathcal{M}}.
\end{equation}
Finally, it is easy to show that the following differential identities hold: 
\begin{subequations}
\label{e:mesh_motion_identities}
\begin{align}
& \frac{\partial \mathcal{L}}{\partial \boldsymbol{u}^{\mathcal{F}}}+\frac{\partial \mathcal{L}}{\partial \boldsymbol{x}^{\mathcal{F}}} = \frac{\partial \mathcal{L}^{\mathcal{M}}}{\partial \boldsymbol{u}^{\mathcal{F}}} = \boldsymbol{0}^T \label{e:mesh_motion_identities_adjoints}\\
& \frac{\partial \mathcal{L}}{\partial \boldsymbol{X}^{\mathcal{F}}} = \frac{\partial \mathcal{L}^{\mathcal{M}}}{\partial \boldsymbol{X}^{\mathcal{F}}} +  \label{e:mesh_motion_identities_sensitivities} \nonumber \\
& \underline{\alpha^{\mathcal{I}} \frac{\partial  J^{\mathcal{I}}}{\partial \boldsymbol{X}^{\mathcal{F}}}} + \underline{(\boldsymbol{\Psi}^{\mathcal{M}}_{\Gamma_\mathcal{I}})^T \cdot \frac{\partial \boldsymbol{r}^{\mathcal{M}}_{\Gamma_\mathcal{I}}}{\partial \boldsymbol{X}^{\mathcal{F}}}} + \underline{(\boldsymbol{\Psi}^{\mathcal{S}}_{\Gamma_\mathcal{I}})^T \cdot \frac{\partial \boldsymbol{f}^{\mathcal{S},ext}_{\Gamma_\mathcal{I}}}{\partial \boldsymbol{X}^{\mathcal{F}}} }
\end{align}
\end{subequations}

\subsection{Coupled adjoint structure problem}
Given adjoint fluid displacements on the fluid interface, the coupled adjoint system and shape sensitivities associated with the structure in Eq. \ref{e:Derivative_of_FSI_Lagrange_Function} read as follows:  
\begin{subequations}
\label{e:adjoint_structure}
\begin{align}
&\frac{\partial \mathcal{L}}{\partial \boldsymbol{u}^{\mathcal{S}}} = \nonumber \\ &   \frac{\partial \tilde{J}}{\partial \boldsymbol{u}^{\mathcal{S}}} + (\boldsymbol{\Psi}^{\mathcal{S}})^T \cdot \frac{\partial \boldsymbol{r}^{\mathcal{S}}}{\partial \boldsymbol{u}^{\mathcal{S}}} + (\boldsymbol{\Psi}^{\mathcal{M}}_{\Gamma_\mathcal{I}})^{T} \cdot \boldsymbol{H}^{\mathcal{S}} \cdot (\frac{\partial \boldsymbol{u}^{\mathcal{S}}}{\partial \boldsymbol{u}^{\mathcal{S}}_{\Gamma_\mathcal{I}}})^T  = \boldsymbol{0}^T \\
& \frac{\partial \mathcal{L}}{\partial \boldsymbol{X}^{\mathcal{S}}} = 
 \frac{\partial \tilde{J}}{\partial \boldsymbol{X}^{\mathcal{S}}} -(\boldsymbol{\Psi}^{\mathcal{S}})^T \cdot \frac{\partial \boldsymbol{f}^{\mathcal{S},int}}{\partial \boldsymbol{X}^{\mathcal{S}}} + \nonumber \\
 &   \underline{\alpha^{\mathcal{I}} \frac{\partial  J^{\mathcal{I}}}{\partial \boldsymbol{X}^{\mathcal{S}}}} + \underline{(\boldsymbol{\Psi}^{\mathcal{M}}_{\Gamma_\mathcal{I}})^T \cdot \frac{\partial \boldsymbol{r}^{\mathcal{M}}_{\Gamma_\mathcal{I}}}{\partial \boldsymbol{X}^{\mathcal{S}}}} + \underline{(\boldsymbol{\Psi}^{\mathcal{S}}_{\Gamma_\mathcal{I}})^T \cdot \frac{\partial \boldsymbol{f}^{\mathcal{S},ext}_{\Gamma_\mathcal{I}}}{\partial \boldsymbol{X}^{\mathcal{S}}}}.
\end{align}
\end{subequations}
As mentioned previously, the underlined terms will vanish with matching interface meshes. Following the idea of using single-disciplinary adjoint solvers for partitioned adjoint FSI analysis, one can define the following structural shape optimization problem whose adjoint system and shape sensitivities are equivalent to Eqs. \ref{e:adjoint_structure} (under the assumption of matching interface meshes):
\begin{equation}
\label{e:structure_coupled_opt_problem}
\begin{aligned}
& \underset{\boldsymbol{X}^{\mathcal{S}}}{\text{min}}
& & \tilde{J}^{\mathcal{S}} = \alpha^{\mathcal{S}} J^{\mathcal{S}} + \alpha^{\mathcal{I}} J^{\mathcal{I}} + J^{\mathcal{S},a} \\
& \text{subject to} & &  \\
& & &  \boldsymbol{r}^{\mathcal{S}}\left(\boldsymbol{u}^{\mathcal{S}},\boldsymbol{X}^{\mathcal{S}}\right)  =  \boldsymbol{0}  \\
& \text{in which} & &  \\
& & & J^{\mathcal{S},a} = (\boldsymbol{f}^{\mathcal{S},a})^T \cdot  \boldsymbol{u}^{\mathcal{S}}; \\
& & &\boldsymbol{f}^{\mathcal{S},a} = \begin{bmatrix}
\left[\boldsymbol{0}\right]_{1 \times m^{\mathcal{S}}_{\Omega}}  \\  \left[\boldsymbol{0}\right]_{1 \times m^{\mathcal{S}}_{\Gamma_{\mathcal{D}}}} \\  
\left[\boldsymbol{0}\right]_{1 \times m^{\mathcal{S}}_{\Gamma_{\mathcal{N}}}} \\  
(\boldsymbol{H}^{\mathcal{S}})^T \cdot \boldsymbol{\Psi}^{\mathcal{M}}_{\Gamma_\mathcal{I}}
\end{bmatrix}_{1 \times m^{\mathcal{S}}}
\end{aligned}
\end{equation}
where $\tilde{J}^{\mathcal{S}}$ is the weighted sum of functions depending on structural displacements. Similar to the coupled adjoint fluid and mesh motion problems, an auxiliary function $J^{\mathcal{S},a}$ is introduced to account for the adjoint coupling of structure and fluid mesh motion, using typical objective functions found in single-disciplinary adjoint solvers. Here motivated by the linear strain energy objective function, $\boldsymbol{f}^{\mathcal{S},a}$ can be interpreted as 
an adjoint force acting only on the fluid-structure interface. 

Last but not least, it can be proved that adjoint-based shape sensitivity analysis for Eq. \ref{e:structure_coupled_opt_problem} leads to the following identities:

\begin{subequations}
\label{e:structure_identities}
\begin{align}
& \frac{\partial \mathcal{L}}{\partial \boldsymbol{u}^{\mathcal{S}}} = \frac{\partial \mathcal{L}^{\mathcal{S}}}{\partial \boldsymbol{u}^{\mathcal{S}}} = \boldsymbol{0}^T \label{e:mesh_motion_identities_adjoints}\\
& \frac{\partial \mathcal{L}}{\partial \boldsymbol{X}^{\mathcal{S}}} = \frac{\partial \mathcal{L}^{\mathcal{S}}}{\partial \boldsymbol{X}^{\mathcal{S}}} +  \label{e:mesh_motion_identities_sensitivities} \nonumber \\
&  \underline{\alpha^{\mathcal{I}} \frac{\partial  J^{\mathcal{I}}}{\partial \boldsymbol{X}^{\mathcal{S}}}} + \underline{(\boldsymbol{\Psi}^{\mathcal{M}}_{\Gamma_\mathcal{I}})^T \cdot \frac{\partial \boldsymbol{r}^{\mathcal{M}}_{\Gamma_\mathcal{I}}}{\partial \boldsymbol{X}^{\mathcal{S}}}} + \underline{(\boldsymbol{\Psi}^{\mathcal{S}}_{\Gamma_\mathcal{I}})^T \cdot \frac{\partial \boldsymbol{f}^{\mathcal{S},ext}_{\Gamma_\mathcal{I}}}{\partial \boldsymbol{X}^{\mathcal{S}}}} \\
& \text{in which} \nonumber \\
& \mathcal{L}^{\mathcal{S}}(\boldsymbol{u}^{\mathcal{S}},\boldsymbol{X}^{\mathcal{S}}) = \tilde{J}^{\mathcal{S}} + (\boldsymbol{\Psi}^{\mathcal{S}})^{T} \cdot \boldsymbol{r}^{\mathcal{S}}.
\end{align}
\end{subequations}

\begin{figure*}
	\includegraphics[keepaspectratio,width=\textwidth]{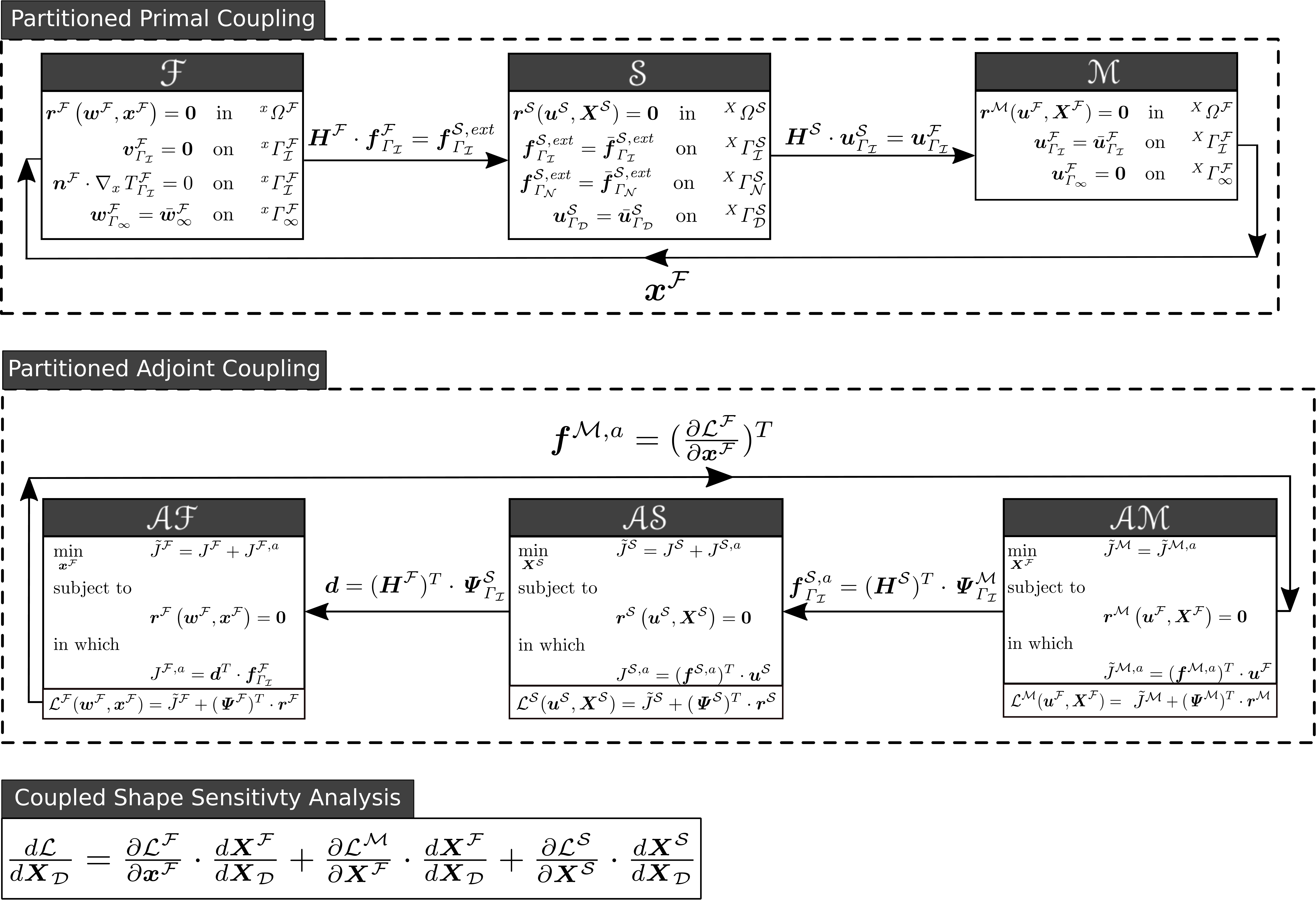}
	\caption{Partitioned multidisciplinary shape sensitivity analysis of steady-state FSI using single-disciplinary primal and adjoint solvers. Note, the simplified form of the coupled shape sensitivity equation is used.}
	\label{fig:AFSI_three_field}
\end{figure*}

\subsection{Coupled shape sensitivity equation}
Having systematically derived strongly coupled adjoint systems in a partitioned manner, one can compute the coupled shape sensitivities from the shape sensitivities delivered by the individual adjoint systems as follows: 
\begin{equation}
\label{e:partitioned_Discrete_Adj_Sens_Eq}
\begin{aligned}
& \frac{d\mathcal{L}}{d\boldsymbol{X}_{\mathcal{D}}} = 
\frac{\partial \mathcal{L}^{\mathcal{F}}}{\partial \boldsymbol{x}^{\mathcal{F}}} \cdot \frac{d \boldsymbol{X}^{\mathcal{F}}}{d\boldsymbol{X}_{\mathcal{D}}} +  \frac{\partial \mathcal{L}^{\mathcal{M}}}{\partial \boldsymbol{X}^{\mathcal{F}}} \cdot \frac{d \boldsymbol{X}^{\mathcal{F}}}{d\boldsymbol{X}_{\mathcal{D}}} + \frac{\partial \mathcal{L}^{\mathcal{S}}}{\partial \boldsymbol{X}^{\mathcal{S}}} \cdot \frac{d \boldsymbol{X}^{\mathcal{S}}}{d\boldsymbol{X}_{\mathcal{D}}} + \\
&   \resizebox{1.02\hsize}{!}{$\left( \underline{\alpha^{\mathcal{I}} \frac{\partial  J^{\mathcal{I}}}{\partial \boldsymbol{X}^{\mathcal{F}}}} + \underline{(\boldsymbol{\Psi}^{\mathcal{M}}_{\Gamma_\mathcal{I}})^T \cdot \frac{\partial \boldsymbol{r}^{\mathcal{M}}_{\Gamma_\mathcal{I}}}{\partial \boldsymbol{X}^{\mathcal{F}}}} + \underline{(\boldsymbol{\Psi}^{\mathcal{S}}_{\Gamma_\mathcal{I}})^T \cdot \frac{\partial \boldsymbol{f}^{\mathcal{S},ext}_{\Gamma_\mathcal{I}}}{\partial \boldsymbol{X}^{\mathcal{F}}} }\right) \cdot \frac{d \boldsymbol{X}^{\mathcal{F}}}{d\boldsymbol{X}_{\mathcal{D}}} + $}\\
&  \resizebox{1.02\hsize}{!}{$\left( \underline{\alpha^{\mathcal{I}} \frac{\partial  J^{\mathcal{I}}}{\partial \boldsymbol{X}^{\mathcal{S}}}} + \underline{(\boldsymbol{\Psi}^{\mathcal{M}}_{\Gamma_\mathcal{I}})^T \cdot \frac{\partial \boldsymbol{r}^{\mathcal{M}}_{\Gamma_\mathcal{I}}}{\partial \boldsymbol{X}^{\mathcal{S}}}} + \underline{(\boldsymbol{\Psi}^{\mathcal{S}}_{\Gamma_\mathcal{I}})^T \cdot \frac{\partial \boldsymbol{f}^{\mathcal{S},ext}_{\Gamma_\mathcal{I}}}{\partial \boldsymbol{X}^{\mathcal{S}}}} \right) \cdot \frac{d \boldsymbol{X}^{\mathcal{S}}}{d\boldsymbol{X}_{\mathcal{D}}}.$} 
\end{aligned}
\end{equation}
As mentioned previously, the underlined terms vanish identically if the interface meshes are matching or are not subject to shape sensitivity analysis. Therefore, the coupled shape sensitivity equation can be further simplified as
\begin{equation}
\label{e:partitioned_Discrete_Adj_Sens_Eq_simplified}
\begin{aligned}
& \frac{d\mathcal{L}}{d\boldsymbol{X}_{\mathcal{D}}} = 
\frac{\partial \mathcal{L}^{\mathcal{F}}}{\partial \boldsymbol{x}^{\mathcal{F}}} \cdot \frac{d \boldsymbol{X}^{\mathcal{F}}}{d\boldsymbol{X}_{\mathcal{D}}} +  \frac{\partial \mathcal{L}^{\mathcal{M}}}{\partial \boldsymbol{X}^{\mathcal{F}}} \cdot \frac{d \boldsymbol{X}^{\mathcal{F}}}{d\boldsymbol{X}_{\mathcal{D}}} + \frac{\partial \mathcal{L}^{\mathcal{S}}}{\partial \boldsymbol{X}^{\mathcal{S}}} \cdot \frac{d \boldsymbol{X}^{\mathcal{S}}}{d\boldsymbol{X}_{\mathcal{D}}}. 
\end{aligned}
\end{equation}   
Analogous to the primal problem, the partitioned adjoint-based FSI sensitivity analysis presented above can be realized by single-disciplinary adjoint solvers in a black-box manner. Figure \ref{fig:AFSI_three_field} illustrates the flows of information in the partitioned primal and adjoint FSI problems using a set of Dirichlet and Neumann-type coupling conditions, where
\begin{itemize}
	\item $\mathscr{AF}$ is the adjoint fluid solver that computes for the domain-based shape sensitivities of the multi-objective function $\tilde{J}^{\mathcal{F}}$ for a given force projection vector on the fluid interface mesh. The projection vector is computed by the transpose matrix-vector product of the force mapping matrix $\boldsymbol{H}^{\mathcal{F}}$ and the interface structural adjoint displacements.
	\item $\mathscr{AM}$ is the adjoint mesh motion solver which can be viewed as an adjoint structural solver computing shape sensitivities for the volumetric strain of the pseudo structure (i.e. the fluid mesh) under the adjoint body force $\boldsymbol{f}^{\mathcal{M},a}$.
	\item  $\mathscr{AS}$ is the adjoint structural solver which can compute for the shape sensitivities of a multi-objective function containing the pseudo interface strain energy $J^{\mathcal{S},a}$ induced by the structural adjoint force $\boldsymbol{f}^{S,a}$.   
\end{itemize}
\begin{table*}
	\caption{Partitioning of exemplary multi-disciplinary target functions to coupled fluid, structure and pseudo-structural (mesh motion) objective functions to be used in single-disciplinary adjoint solvers {$\mathscr{AF}, \mathscr{AS}, \mathscr{AM}$}. $\boldsymbol{D}$ and $\boldsymbol{A}^\mathcal{F}$ are respectively drag and area normal vector fields.} % title of Table
	\centering  % used for centering table
	\resizebox{\textwidth}{!}{\begin{tabular}{c c c c c} % centered columns (4 columns)
			\hline 
			\hline \\
			\multicolumn{1}{l}{\scalebox{2.5}{function}} & \multicolumn{1}{c}{\scalebox{2.5}{$\tilde{J}$}} & \multicolumn{1}{c}{\scalebox{2.5}{$\tilde{J}^{\mathcal{F}}$}} & \multicolumn{1}{c}{\scalebox{2.5}{$\tilde{J}^{\mathcal{S}}$}} & \multicolumn{1}{c}{\scalebox{2.5}{$\tilde{J}^{\mathcal{M}}$}}  \\
			\hline \\
			\multicolumn{1}{l}{\scalebox{2}{Interface drag.}} & 
			\multicolumn{1}{c|}{\scalebox{2.0}{$\boldsymbol{D}^T\cdot \boldsymbol{f}^{\mathcal{F}}_{\Gamma_\mathcal{I}}$}}  &  \multicolumn{1}{|c|}{\scalebox{2.0}{$\boldsymbol{D}^T \cdot \boldsymbol{f}^{\mathcal{F}}_{\Gamma_\mathcal{I}}+\boldsymbol{d}^T \cdot \boldsymbol{f}^{\mathcal{F}}_{\Gamma_\mathcal{I}}$}} &  \multicolumn{1}{|c|}{\scalebox{2.0}{$(\boldsymbol{f}^{\mathcal{S},a})^T \cdot  \boldsymbol{u}^{\mathcal{S}}$}} &  \multicolumn{1}{|c}{\scalebox{2.0}{$(\boldsymbol{f}^{\mathcal{M},a})^T \cdot \boldsymbol{u}^{\mathcal{F}}$}} \\\\\\[0.25cm]
			\multicolumn{1}{l}{\scalebox{2}{Total power loss.}} & \multicolumn{1}{c|}{\scalebox{2.0}{$\sum\limits_{i=1}^{ m^{\mathcal{F}}_{\Gamma_\infty}}\left[(p^{\mathcal{F}}+\frac{1}{2}\norm{\boldsymbol{v}^{\mathcal{F}}}^2)(\boldsymbol{A}^{\mathcal{F}}\cdot\boldsymbol{v}^{\mathcal{F}})\right]_{\Gamma_{\infty,i}}$}}  &  \multicolumn{1}{|c|}{\specialcell{\scalebox{2.0}{$\sum\limits_{i=1}^{ m^{\mathcal{F}}_{\Gamma_\infty}}\left[(p^{\mathcal{F}}+\frac{1}{2}\norm{\boldsymbol{v}^{\mathcal{F}}}^2)(\boldsymbol{A}^{\mathcal{F}}\cdot\boldsymbol{v}^{\mathcal{F}})\right]_{\Gamma_{\infty,i}}$} \\[0.4cm]\scalebox{2.0}{$ +\boldsymbol{d}^T\cdot \boldsymbol{f}^{\mathcal{F}}_{\Gamma_\mathcal{I}}$}}} &  \multicolumn{1}{|c|}{\scalebox{2.0}{$(\boldsymbol{f}^{\mathcal{S},a})^T \cdot  \boldsymbol{u}^{\mathcal{S}}$}} &  \multicolumn{1}{c}{\scalebox{2.0}{$(\boldsymbol{f}^{\mathcal{M},a})^T \cdot \boldsymbol{u}^{\mathcal{F}}$}} \\\\\\[0.25cm]
			\multicolumn{1}{l}{\scalebox{2}{Fluid interface energy.}} & 
			\multicolumn{1}{c|}{\scalebox{2.0}{$(\boldsymbol{H}^{\mathcal{S}} \cdot \boldsymbol{u}^{\mathcal{S}}_{\Gamma_\mathcal{I}})^T\cdot \boldsymbol{f}^{\mathcal{F}}_{\Gamma_\mathcal{I}}$}}  &  \multicolumn{1}{|c|}{\scalebox{2.0}{$(\boldsymbol{H}^{\mathcal{S}} \cdot \boldsymbol{u}^{\mathcal{S}}_{\Gamma_\mathcal{I}})^T\cdot \boldsymbol{f}^{\mathcal{F}}_{\Gamma_\mathcal{I}}+\boldsymbol{d}^T\cdot \boldsymbol{f}^{\mathcal{F}}_{\Gamma_\mathcal{I}}$}} &  \multicolumn{1}{|c|}{\scalebox{2.0}{$((\boldsymbol{H}^{\mathcal{S}})^T \cdot \boldsymbol{f}^{\mathcal{F}}_{\Gamma_\mathcal{I}})^T \cdot  \boldsymbol{u}^{\mathcal{S}}_{\Gamma_\mathcal{I}}+(\boldsymbol{f}^{\mathcal{S},a})^T \cdot  \boldsymbol{u}^{\mathcal{S}}$}} &  \multicolumn{1}{|c}{\scalebox{2.0}{$(\boldsymbol{f}^{\mathcal{M},a})^T \cdot \boldsymbol{u}^{\mathcal{F}}$}} \\\\\\[0.25cm]
			\multicolumn{1}{l}{\scalebox{2}{Structural interface energy.}} & 
			\multicolumn{1}{c|}{\scalebox{2.0}{$( \boldsymbol{u}^{\mathcal{S}}_{\Gamma_\mathcal{I}})^T\cdot (\boldsymbol{H}^{\mathcal{F}}\cdot\boldsymbol{f}^{\mathcal{F}}_{\Gamma_\mathcal{I}})$}}  &  \multicolumn{1}{|c|}{\scalebox{2.0}{$((\boldsymbol{H}^{\mathcal{F}})^T \cdot \boldsymbol{u}^{\mathcal{S}}_{\Gamma_\mathcal{I}})^T\cdot \boldsymbol{f}^{\mathcal{F}}_{\Gamma_\mathcal{I}}+\boldsymbol{d}^T\cdot \boldsymbol{f}^{\mathcal{F}}_{\Gamma_\mathcal{I}}$}} &  \multicolumn{1}{|c|}{\scalebox{2.0}{$(\boldsymbol{H}^{\mathcal{F}} \cdot \boldsymbol{f}^{\mathcal{F}}_{\Gamma_\mathcal{I}})^T \cdot  \boldsymbol{u}^{\mathcal{S}}_{\Gamma_\mathcal{I}}+(\boldsymbol{f}^{\mathcal{S},a})^T \cdot  \boldsymbol{u}^{\mathcal{S}}$}} &  \multicolumn{1}{c}{\scalebox{2.0}{$(\boldsymbol{f}^{\mathcal{M},a})^T \cdot \boldsymbol{u}^{\mathcal{F}}$}} 
	\end{tabular}}
	\label{table:summary_cost_functions} % is used to refer this table in the text
\end{table*}
Based on the partitioned primal and adjoint FSI workflows in Figure \ref{fig:AFSI_three_field} and the single-disciplinary solvers therein, the following remarks can be added:	
\begin{remark}
Both the data flow and the matrix operations in the adjoint problem are reversed compared to the primal problem. This observation is in correct agreement with the general adjoint-based sensitivity analysis.  
\end{remark}
\begin{remark}
The single-disciplinary adjoint solvers $\{\mathscr{AF}, \mathscr{AM}, \mathscr{AS}\}$ are coupled 
with each other by augmenting the associated objective functions with the auxiliary objective functions $\{J^{\mathcal{F},a}, J^{\mathcal{M},a}, J^{\mathcal{S},a}\}$, respectively. The auxiliary functions are either force-based or displacement-based functionals which are readily available in well-established single-disciplinary adjoint solvers. Although, the presented coupling scheme is independent of the derivation and implementation of the underlying adjoint solvers, it imposes on them the requirement of handling shape sensitivity analysis for a weighted sum of objectives, $\{\tilde{J}^{\mathcal{F}}, \tilde{J}^{\mathcal{M}}, \tilde{J}^{\mathcal{S}}\}$, using a single adjoint solution. For the sake of clarity, Table. \ref{table:summary_cost_functions} lists  exemplary target functions which are partitioned by means of the presented scheme. Note that the fluid and structural interface energies are equal in case of the energy conservative spatial mapping (see Eqs. \ref{eq:interface_energy},\ref{eq:interface_energy2}).
   
\end{remark}

\begin{remark}
Considering the coupled adjoint mesh motion problem \ref{e:mesh_motion_coupled_opt_problem}, it is observed that the presented partitioning requires the partial derivatives of the fluid Lagrange function w.r.t the internal and boundary nodes of the fluid mesh, i.e.,
\begin{equation}
\frac{\partial \mathcal{L}^{\mathcal{F}}}{\partial \boldsymbol{x}^{\mathcal{F}}} = 
\begin{bmatrix}
\frac{\partial \mathcal{L}^{\mathcal{F}}}{\partial \boldsymbol{x}^{\mathcal{F}}_{\Omega}}  & & 
\frac{\partial \mathcal{L}^{\mathcal{F}}}{\partial \boldsymbol{x}^{\mathcal{F}}_{\Gamma_{\infty}}}  & & \frac{\partial \mathcal{L}^{\mathcal{F}}}{\partial \boldsymbol{x}^{\mathcal{F}}_{\Gamma_{\mathcal{I}}}} 
\end{bmatrix}_{1 \times m^{\mathcal{F}}}.
\end{equation}
The domain term $\frac{\partial \mathcal{L}^{\mathcal{F}}}{\partial \boldsymbol{x}^{\mathcal{F}}_{\Omega}}$ might not be computed and available by every fluid adjoint solver, e.g., due to the so-called reduced gradient or boundary-based formulations \citep{Kavvadias2015,Lozano2017}. However under the condition that the fluid solution is invariant w.r.t the fluid interior mesh, the domain geometric derivatives can be assumed to be zero. As a result, the coupled adjoint mesh motion problem reads
\begin{equation}
\begin{aligned}
& \begin{bmatrix} 
(\frac{\partial \boldsymbol{f}^{\mathcal{M},int}_{\Omega}}{\partial \boldsymbol{u}^{\mathcal{F}}_{\Omega}})^T & \left[\boldsymbol{0}\right]_{ m^{\mathcal{F}}_{\Omega} \times m^{\mathcal{F}}_{\Gamma_{\infty}}} & \left[\boldsymbol{0}\right]_{ m^{\mathcal{F}}_{\Omega} \times m^{\mathcal{F}}_{\Gamma_{\mathcal{I}}}} \\ 
(\frac{\partial \boldsymbol{f}^{\mathcal{M},int}_{\Omega}}{\partial \boldsymbol{u}^{\mathcal{F}}_{\Gamma_{\infty}}})^T & \left[\boldsymbol{I}\right]_{ m^{\mathcal{F}}_{\Gamma_{\infty}} \times m^{\mathcal{F}}_{\Gamma_{\infty}}} & \left[\boldsymbol{0}\right]_{ m^{\mathcal{F}}_{\Gamma_{\infty}} \times m^{\mathcal{F}}_{\Gamma_{\mathcal{I}}}} \\ 
(\frac{\partial \boldsymbol{f}^{\mathcal{M},int}_{\Omega}}{\partial \boldsymbol{u}^{\mathcal{F}}_{\Gamma_{\mathcal{I}}}})^T & \left[\boldsymbol{0}\right]_{ m^{\mathcal{F}}_{\Gamma_{\mathcal{I}}} \times m^{\mathcal{F}}_{\Gamma_{\infty}}} & \left[\boldsymbol{I}\right]_{ m^{\mathcal{F}}_{\Gamma_{\mathcal{I}}} \times m^{\mathcal{F}}_{\Gamma_{\mathcal{I}}}}
\end{bmatrix} \begin{bmatrix}
\boldsymbol{\Psi}^{\mathcal{M}}_{\Omega} \\ \\
\boldsymbol{\Psi}^{\mathcal{M}}_{\Gamma_{\infty}} \\ \\
\boldsymbol{\Psi}^{\mathcal{M}}_{\Gamma_{\mathcal{I}}}
\end{bmatrix} = \\ 
&  \begin{bmatrix}
\left[\boldsymbol{0}\right]_{m^{\mathcal{F}}_{\Omega} \times 1}  \\ \\ 
(\frac{\partial \mathcal{L}^{\mathcal{F}}}{\partial \boldsymbol{x}^{\mathcal{F}}_{\Gamma_{\infty}}})^T \\ \\ (\frac{\partial \mathcal{L}^{\mathcal{F}}}{\partial \boldsymbol{x}^{\mathcal{F}}_{\Gamma_{\mathcal{I}}}})^T 
\end{bmatrix}
\end{aligned} 
\end{equation}
which can be solved analytically, giving
\begin{equation}
\label{eq:two_field_formulation_adjoint_displcament}
\boldsymbol{\Psi}^{\mathcal{M}} =
\begin{bmatrix}
\left[\boldsymbol{0}\right]_{1 \times m^{\mathcal{F}}_{\Omega}}  & & 
\frac{\partial \mathcal{L}^{\mathcal{F}}}{\partial \boldsymbol{x}^{\mathcal{F}}_{\Gamma_{\infty}}}  & & \frac{\partial \mathcal{L}^{\mathcal{F}}}{\partial \boldsymbol{x}^{\mathcal{F}}_{\Gamma_{\mathcal{I}}}} 
\end{bmatrix}^T_{1 \times m^{\mathcal{F}}}.
\end{equation}
Clearly elimination of the domain term $\frac{\partial \mathcal{L}^{\mathcal{F}}}{\partial \boldsymbol{x}^{\mathcal{F}}_{\Omega}}$ reduces the number of adjoint problems that should be solved numerically. In other words, the so-called reduced boundary gradients of the fluid result in a reduced formulation of the presented partitioned scheme for the adjoint FSI problem.     

\end{remark}	
	
\begin{figure}
  \centering
  \begin{tabular}{@{}c@{}}
    \includegraphics[width=0.50\textwidth,keepaspectratio]{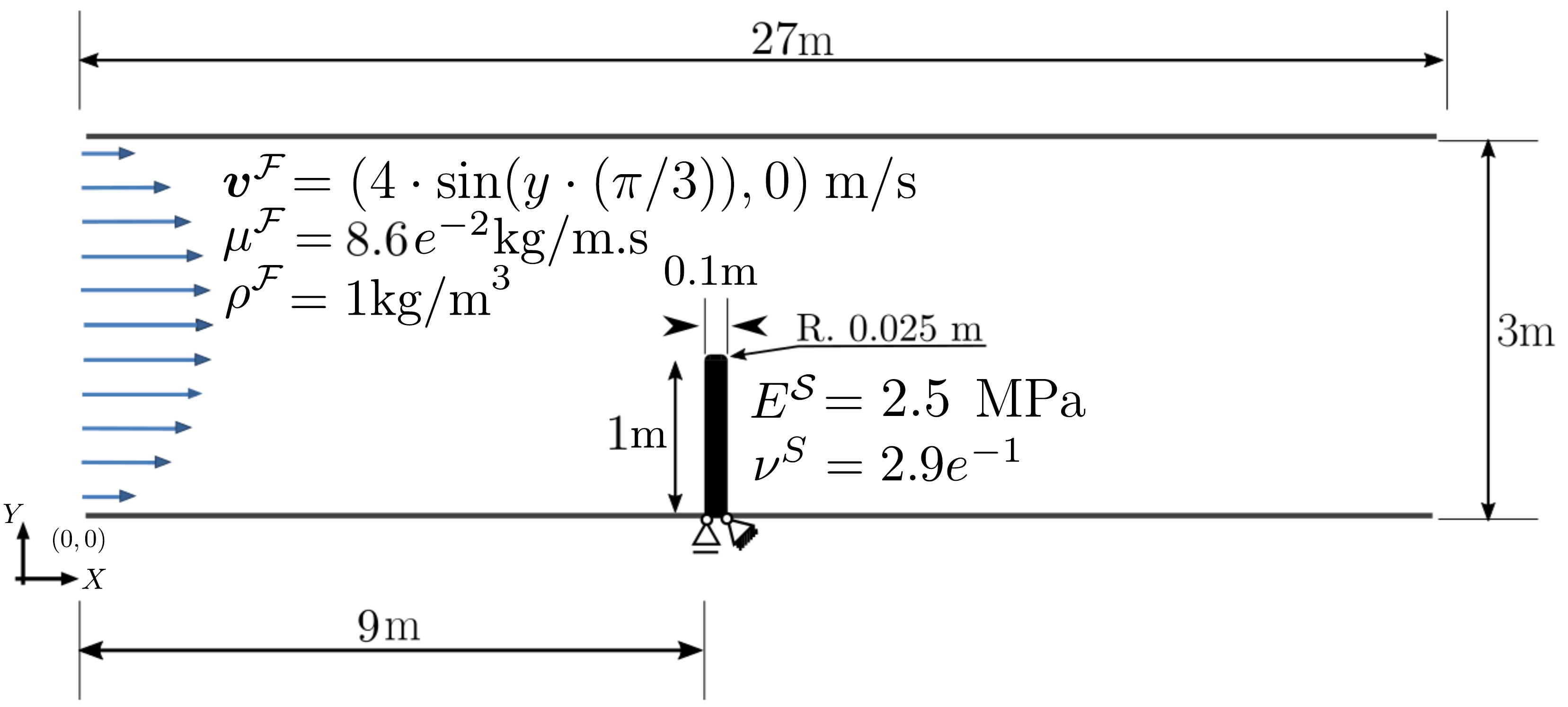} \\[\abovecaptionskip]
    \small (a) Problem setup.
  \end{tabular}

  \vspace{\floatsep}

  \begin{tabular}{@{}c@{}}
    \includegraphics[width=0.50\textwidth,keepaspectratio]{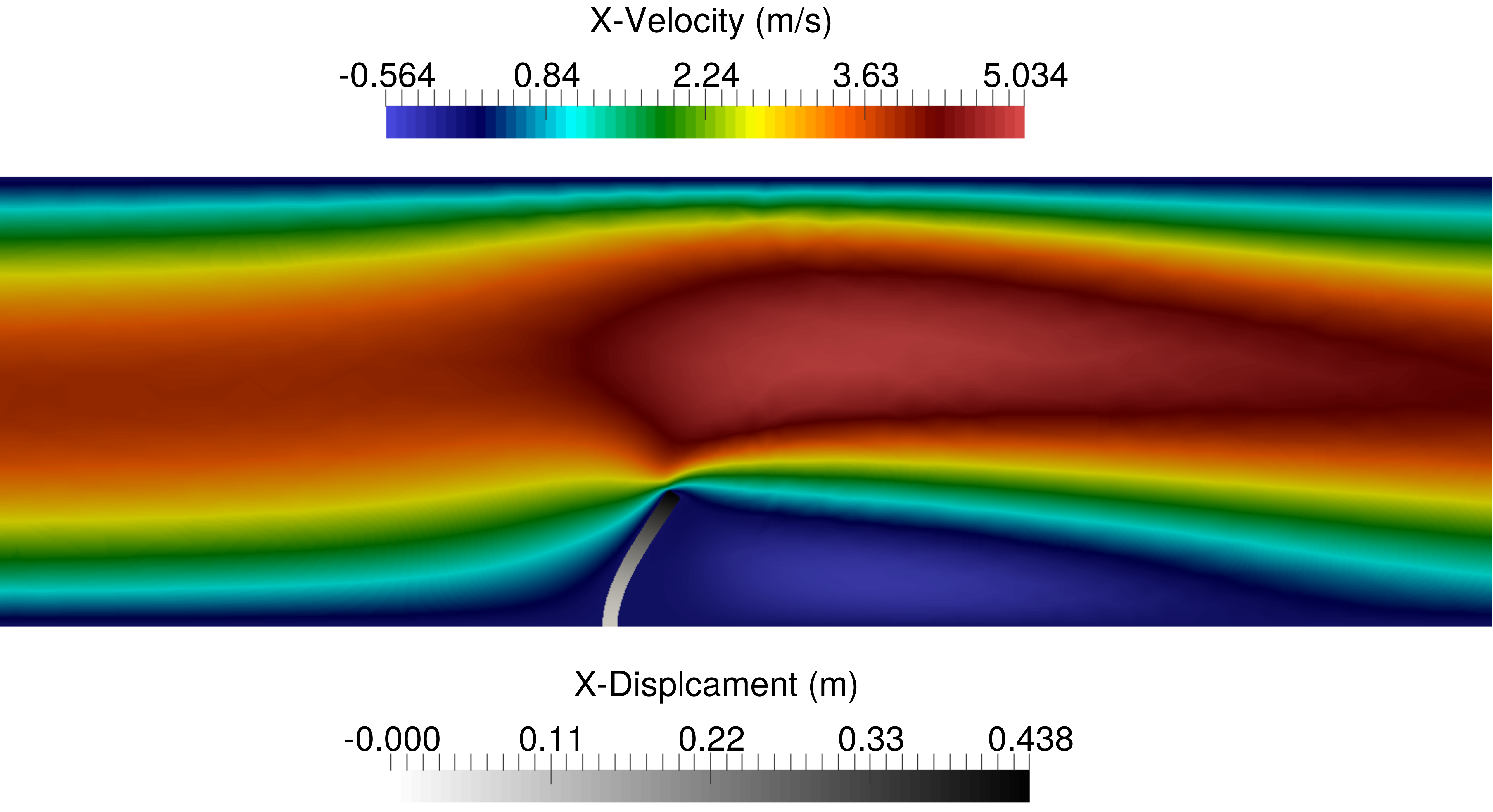} \\[\abovecaptionskip]
    \small (b) Steady-state FSI solution.
  \end{tabular}

  \caption{Flexible beam in a channel.}\label{fig:flexible_beam_in_channel}
\end{figure}

\section{Numerical studies}
\label{numerical_studies}
This section demonstrates the accuracy and general applicability of the proposed partitioned scheme. For this purpose, two multiphysics frameworks are considered: one fully FEM-based and another one hybrid FEM-FVM-based. Herein, the FEM-based analyses, including primal and adjoint-based shape sensitivity analyses, are performed using the open-source software KRATOS Multiphysics (for a detailed description refer to \citep{Dadvand2010,KR2018}). Whereas, the FVM-based computations are done through the open source SU2 suite \citep{Economon2016,SU22018}. Here, the spatial couplings of non-matching meshes are realized by an extended version of the open source coupling tool EMPIRE \citep{Wang2016,empire2018}. This tool offers field mapping technologies which are commonly found in literature to transfer information between non-matching meshes in FSI computations. The following paragraphs contain brief descriptions of the solution strategies devised in the above mentioned software packages for shape sensitivity analysis.   

In KRATOS, adjoint-based shape sensitivity analysis for fluids and structures are performed by discrete analytic method and semi-analytic discrete method, respectively. This means that although the adjoint fluid and structural solvers are derived discretely using the exact analytic Jacobian of the underlying nonlinear equations w.r.t the state variables, the partial shape derivatives (local form) of the fluid and structural optimization problems are computed analytically and approximately by a finite difference scheme, respectively.

SU2 is equipped with continuous and AD-based discrete adjoint fluid solvers, each of which has notable properties. The continuous adjoint solver has shown to be robust and efficient in terms of applicability to large-scale problems with complex geometries \citep{palacios2015large}, however the quality of the computed shape gradients is somewhat low and depends strongly on the mesh quality. This can be explained by the reduced boundary formulation \citep{Economon2014}. On the other hand, the discrete adjoint solver of SU2 provides the numerically exact shape gradients by reformulating the adjoint problem as a fixed-point problem in order to exploit the fixed-point structure of the flow solvers \citep{albring2015development,AlSaGa2016}.

Last but not least, it is should be mentioned that the presented partitioned scheme directly inherits and retains the accuracy, scalability and computational efficiency of the underlying single-disciplinary adjoint solvers.

\subsection{FEM-based shape sensitivity analysis for FSI}
\label{FEM-based_shape_sensitivity_analysis_for_FSI}
The first test case is a flow in a channel obstructed by a flexible beam as illustrated in Figure \ref{fig:flexible_beam_in_channel}. Different setups of this problem have also been served for test purposes in literature \citep{Sanchez2017,richter2012goal,hetu1992fast,carvalho1987predictions}. Here, the case is computed in 2D, and consists of a cantilever beam immersed in a flow with Reynolds number of 10 and driven by a sinusoidal inflow profile with average velocity 0.45 m/s. Furthermore, in order to avoid the serious influence of geometric singularities on the shape gradient accuracy, the sharp corners of the interface are cured with a fillet of 25 mm radius. The reader is referred to \citep{anderson1999aerodynamic,castro2007systematic,Lozano2017,lozano2019watch} for a thorough discussions on the influence of geometric and solution singularities on the behavior of shape sensitivities.   

The fluid is modelled by the incompressible Navier-Stokes equations and the beam is modelled by a hyperelastic continuum under plain-strain conditions. A stabilized finite element method based on SUPG/PSPG \citep{TEZDUYAR1992221} stabilization with first order triangular elements is used for the spatial discretization of the fluid, while the structure domain is discretized with standard triangular elements. For the sake of error reduction, the fluid and structure domains are discretized with a conforming interface mesh.  For this test case, KRATOS finite element framework is used for the shape sensitivity analysis of the cost functions defined in Table \ref{table:summary_cost_functions}.

The steady-state solution of the coupled FSI problem is depicted in Figure \ref{fig:flexible_beam_in_channel} (b). As can be seen, the structure undergoes large displacements due to the laminar flow field. It was observed that the drag force on the structure drops by $23\%$ through the transition from the undeformed state to the equilibrium state. Also, the linear strain energy (interface energy) of the structure varies from practically $0 \,\text{N.m}$ in the undeformed state to $1.58\,\text{N.m}$ in the equilibrium state.  

\begin{figure}
  \centering
  \begin{tabular}{@{}c@{}}
    \includegraphics[width=\linewidth,keepaspectratio]{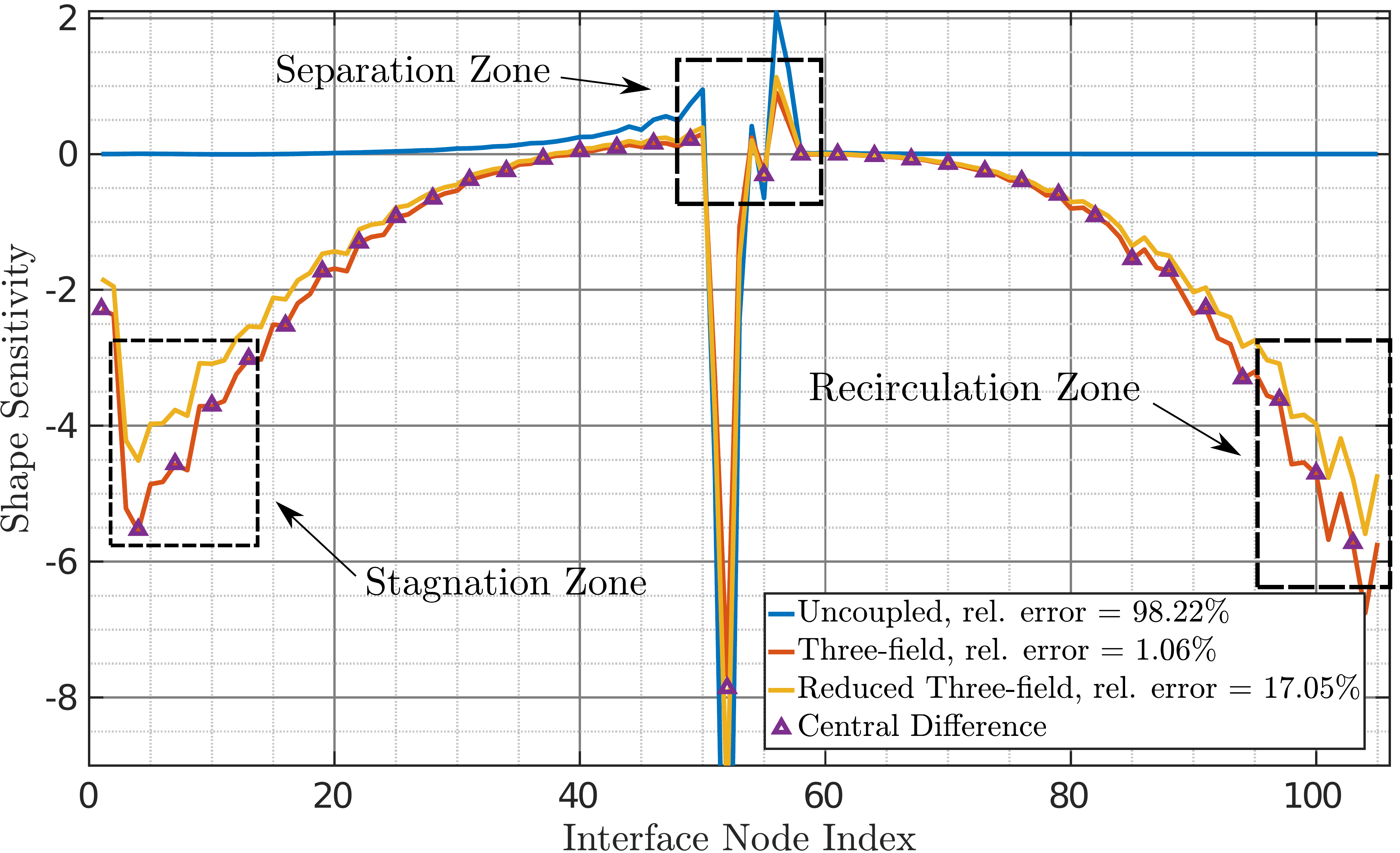} \\[\abovecaptionskip]
    \small (a) Interface drag.
  \end{tabular}

  \vspace{\floatsep}

  \begin{tabular}{@{}c@{}}
    \includegraphics[width=\linewidth,keepaspectratio]{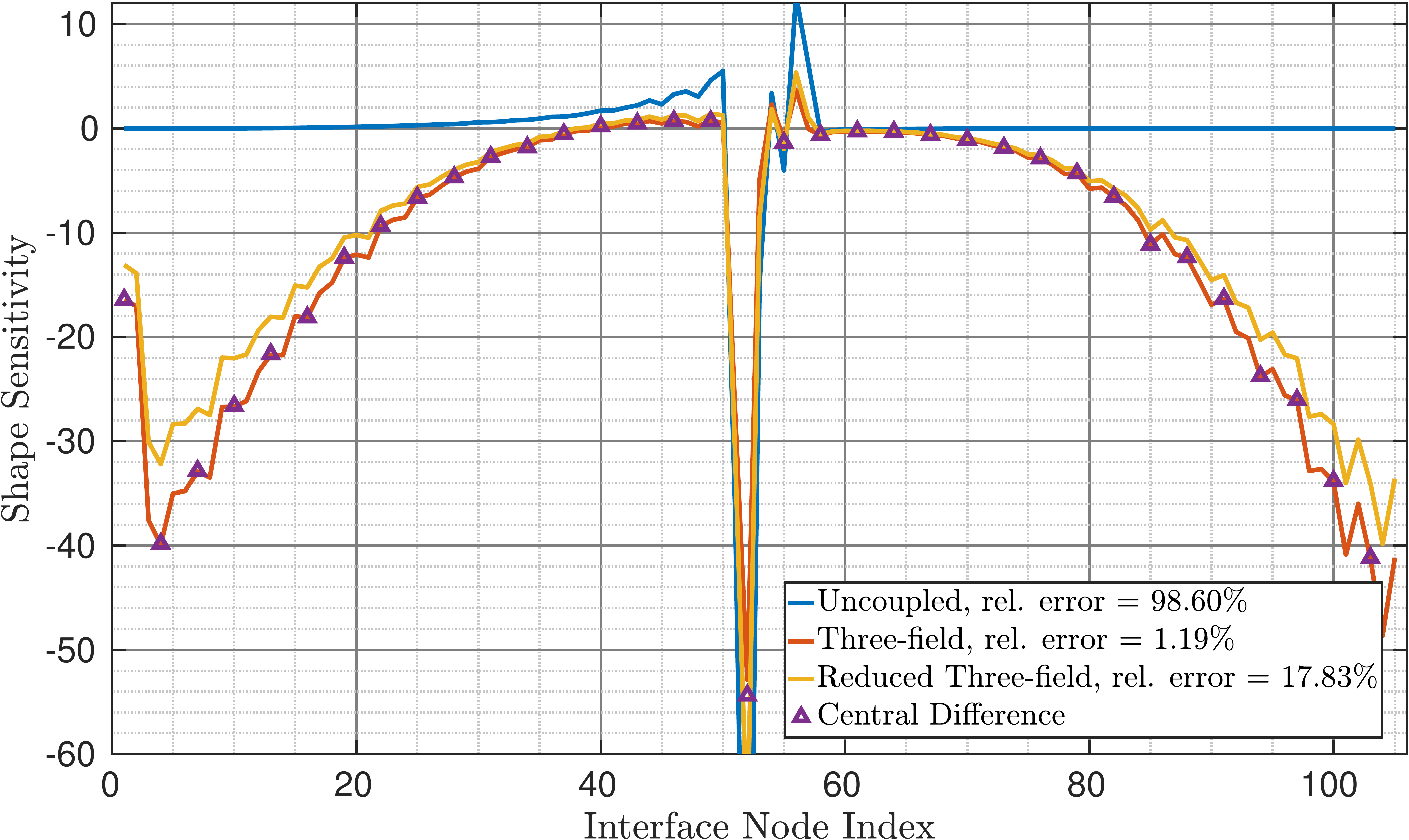} \\[\abovecaptionskip]
    \small (b) Power dissipated by channel.
  \end{tabular}

  \vspace{\floatsep}

  \begin{tabular}{@{}c@{}}
    \includegraphics[width=\linewidth,keepaspectratio]{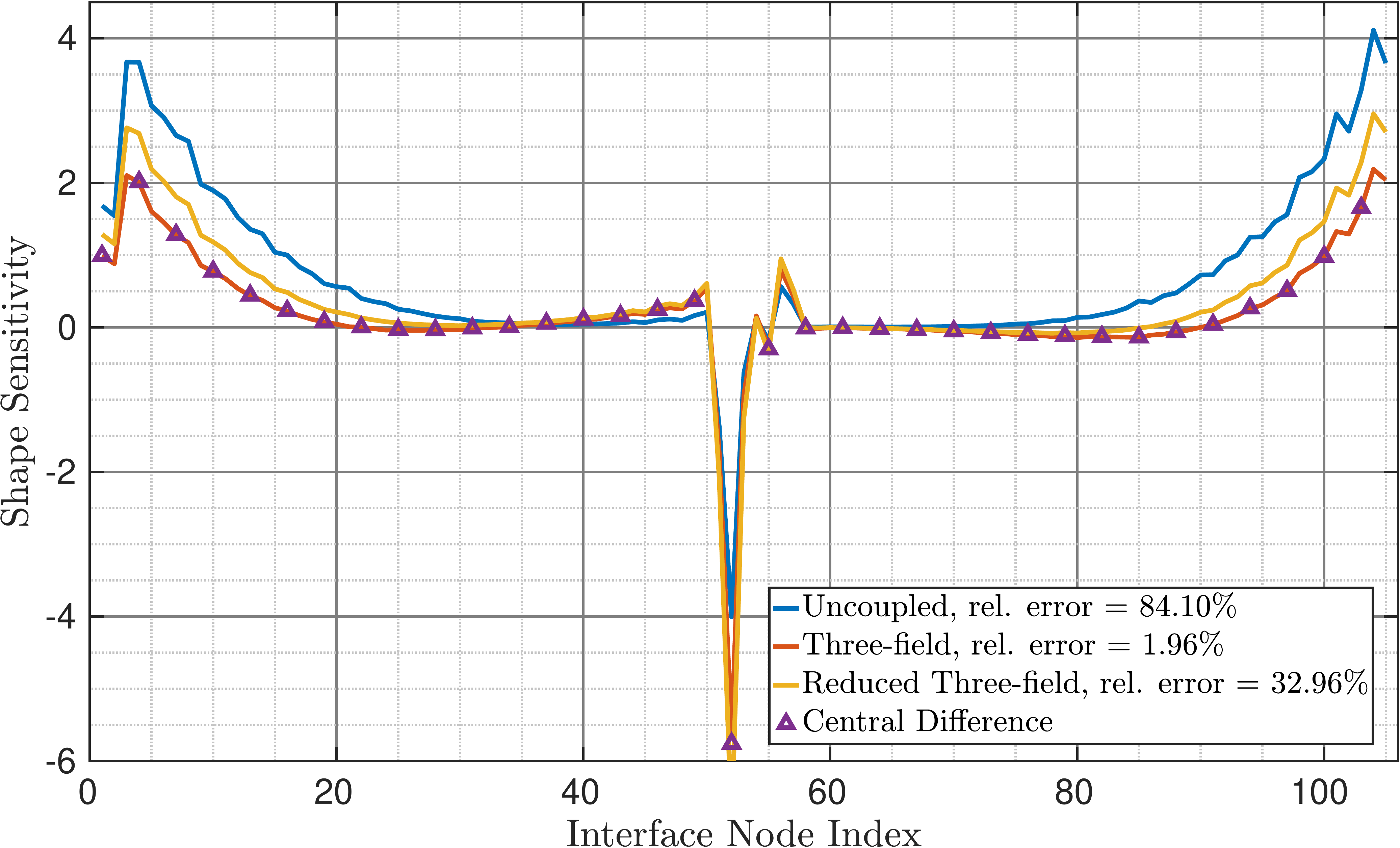} \\[\abovecaptionskip]
    \small (c) Interface energy.
  \end{tabular}

  \caption{Verification and comparison of shape gradients for different objective functions. Shape sensitivities are computed w.r.t the undeformed interface shape whereas the objectives are evaluated at the deformed equilibrium configuration. Perturbations and sensitivities are projected onto the interface normal.}\label{fig:case_1_FD_comparisons}
\end{figure}

\begin{figure*}
	\includegraphics[keepaspectratio,width=\textwidth]{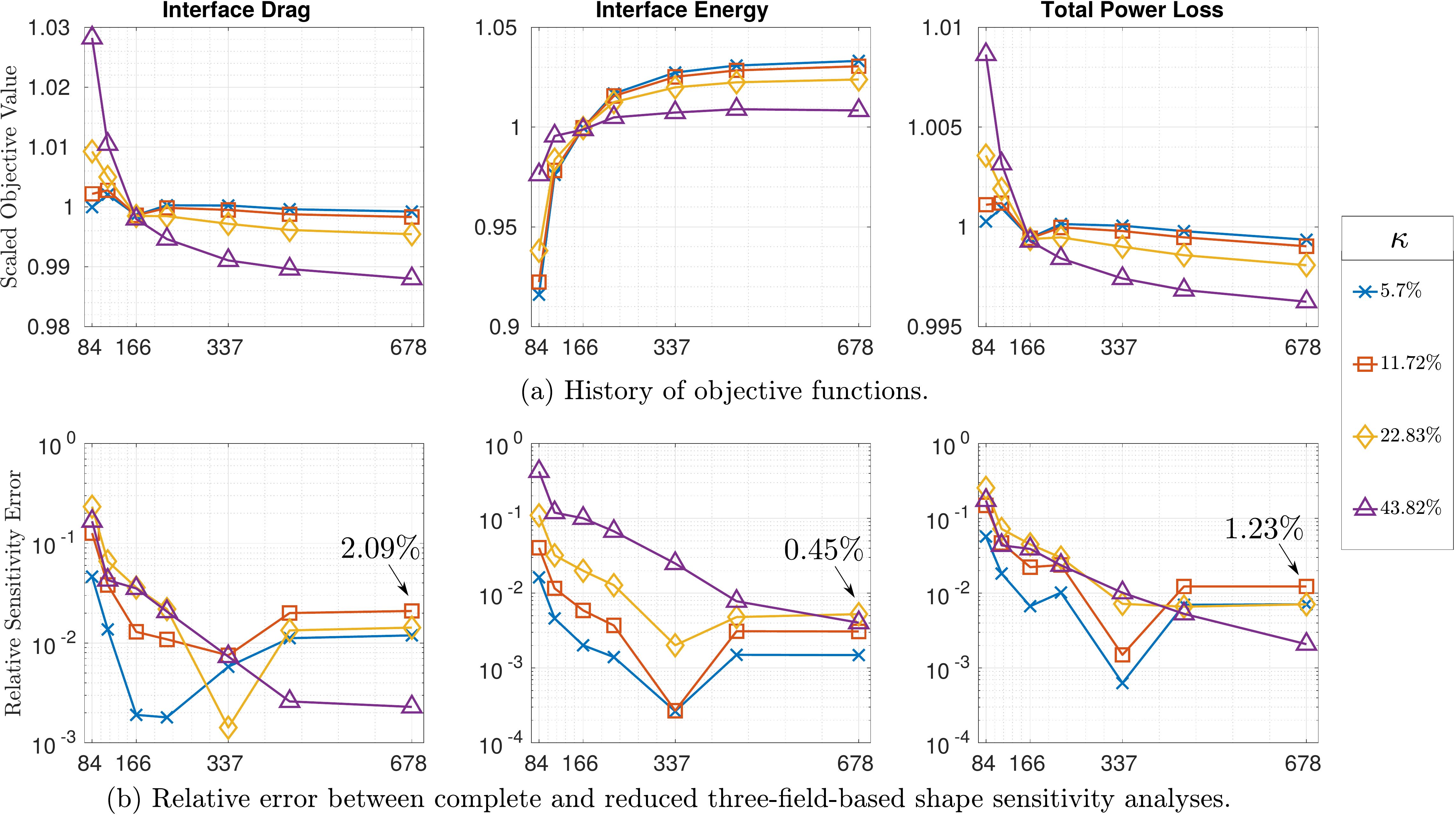}
	\caption{Mesh refinement studies for several levels of the interface flexibility. Horizontal axes represent number of the interface elements.}
	\label{fig:mesh_convergence_studies}
\end{figure*}

\subsubsection{Numerical verification and comparison studies}
\label{sec:Numerical verification and comparison studies}
Verification of the adjoint-based shape sensitivity analysis in Section \ref{Aeroelastic_sensitivity_analysis} is numerically performed against the central difference method (CD). Node coordinates on the complete interface are chosen as design variables and only boundary-normal perturbations are considered. A fixed step size of $\epsilon=10^{-5}$ has been used in all cases. Computations are performed on a coarse mesh with 84 interface elements, permitting the usage of CD for a number of interface nodes.  

Particular focus is placed on the necessity of the strongly coupled adjoint FSI analysis as well as the difference in accuracy between the complete and reduced three-field gradient formulations (refer to Remark 6 for details). Figure \ref{fig:case_1_FD_comparisons} compares the accuracy between different schemes to the reference (CD results) for various objective functions. As a measure of accuracy, we use a relative error based on the L2-norm of absolute error and the reference. As expected, it can clearly be seen that the uncoupled adjoint-based shape sensitivities have the wrong sign and pattern. It can also be observed from the plots that while the complete three-field-based gradients match, qualitatively and quantitatively, with the reference, the reduced three-field-based gradients resemble qualitatively the correct gradients. Significant discrepancies appear around the stagnation, separation, and recirculation zones. This behavior has been observed for fluids by \cite{Lozano2017} and it is generally concluded that reduced/boundary formulations are inaccurate and strongly mesh dependent in such regions, unless the mesh is refined.

\subsubsection{Mesh studies}
Difference between the complete and reduced gradient formulations stems from the lack of the interior fluid sensitivities ($\frac{\partial \mathcal{L}^{\mathcal{F}}}{\partial \boldsymbol{x}^{\mathcal{F}}_{\Omega}}$). Based on the comprehensive mesh sensitivity analyses performed by \cite{Lozano2017} and discussions in \citep{castro2007systematic,anderson1999aerodynamic}, and also the fact that ideally the discrete FSI solution must be invariant w.r.t the fluid mesh, a series of mesh studies are carried out to investigate the inconsistencies observed between these formulations in \ref{sec:Numerical verification and comparison studies}. Mesh refinement is performed on the undeformed fluid and structure geometries. Furthermore, since the fluid solution in FSI is computed on the deformed mesh (see Eq. \ref{eq:Res_Form_FSI_with_BCs}), the interface flexibility is also chosen as a parameter in the investigations. For this purpose, we use level of the interface flexibility and it is defined as $\kappa = max(\boldsymbol{u}^{\mathcal{S}}_{\Gamma_\mathcal{I}})/l$, where $l=1 m$ is the beam length. 

Figure \ref{fig:mesh_convergence_studies}(a) presents convergence histories of the considered objective functions with respect to the number of interface elements, for different levels of the interface flexibility. The plots show that, the objective functions converge to the finest mesh results. Subsequently, Figure \ref{fig:mesh_convergence_studies}(b) shows that overall the integrated relative error between the two formulations drops to acceptable levels as the mesh is refined. In some cases, the error increases and stagnates from a certain level of refinement on. This behavior may be explained by the strong dependency of semi-analytic sensitivities, computed by the single-disciplinary adjoint structural solver of KRATOS, on the finite-difference step size (especially on fine meshes). The results also show trends of smaller error and faster decay in the small-strain structure case (i.e. $\kappa=5.7\%$). This observation may confirm the validity of the reduced three-field-based and two-field-based \citep{Heners2017,Stavropoulou2016,Fazzolari2007} shape sensitivity analyses for FSI with small strains. Nevertheless, in fluids community, the reduced gradient formulations are found to be a good compromise between computational costs and accuracy. Remember that the reduced formulations eliminate the computational cost of the domain geometric variations of the fluid ($\frac{\partial \mathcal{L}^{\mathcal{F}}}{\partial \boldsymbol{x}^{\mathcal{F}}_{\Omega}}$) and the adjoint mesh motion problem subsequently.

\begin{figure*}
\includegraphics[keepaspectratio,width=\textwidth]{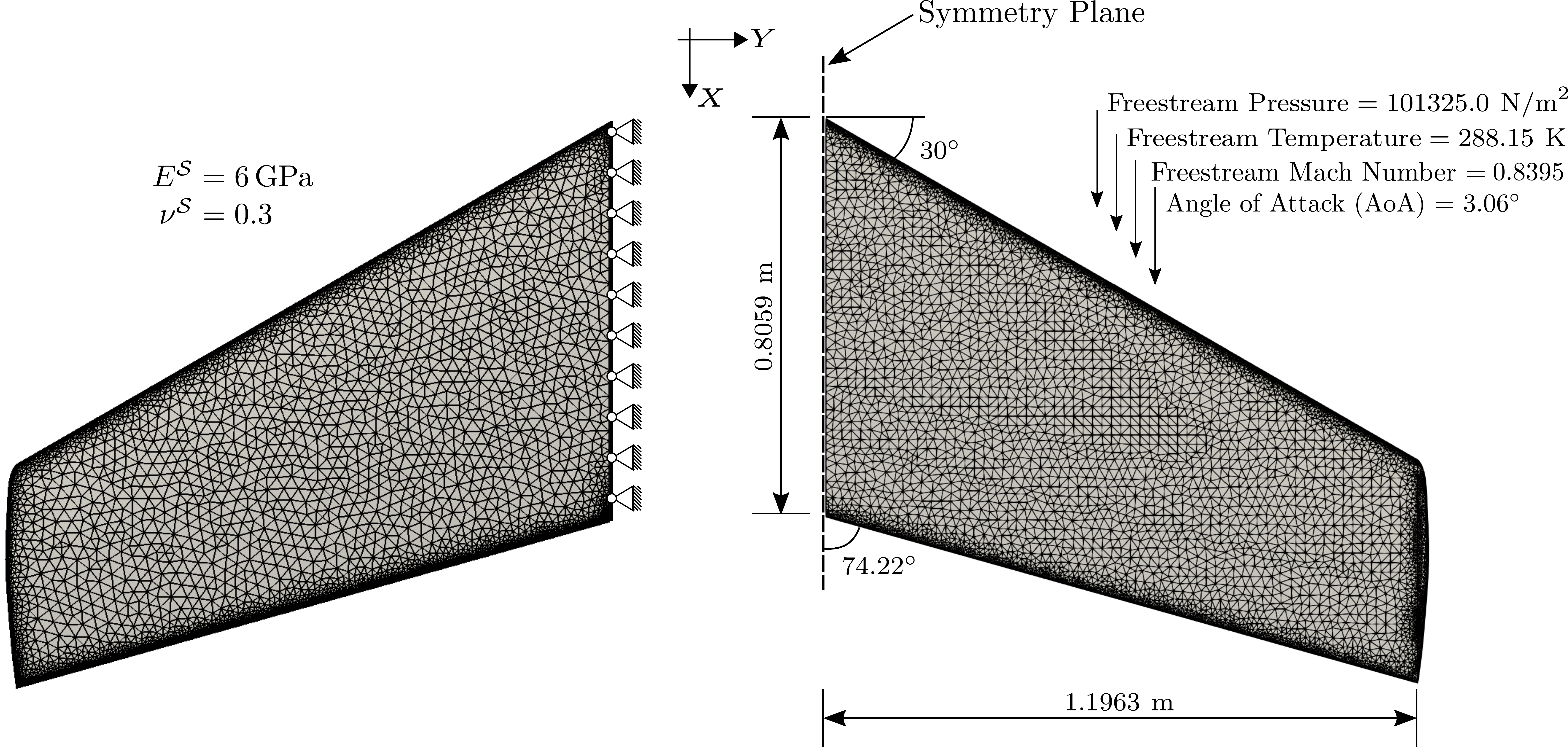}
\caption{Description and surface discretization of ONERA M6 for FSI. Left: structural model, right: fluid model. The fluid and structure interface meshes consist of 18,285 and 18,039 nodes, respectively.}
\label{fig:oneraM6}
\end{figure*} 

\subsection{Hybrid FEM-FVM-based shape sensitivity analysis of ONERA M6}
\label{Open_source_framework_for_shape_sensitivity_analysis}
Using the framework of EMPIRE-KRATOS-SU2, we performed the adjoint-based shape sensitivity analysis of the ONERA M6 wing immersed in a compressible inviscid fluid flow. In contrast to the usual analyses in the literature, we do not consider the wing to be rigid, but model it as a flexible solid structure clamped at the wing root, Figure \ref{fig:oneraM6}. Doing so, we introduce an artificial fluid-structure interaction in the model so that the corresponding shape sensitivity analysis becomes a coupled problem. The rather simple wing structure is chosen since we are focusing here on the performance of the approaches derived in Section \ref{Aeroelastic_sensitivity_analysis}. For both the fluid analysis (CFD) and the structural analysis (CSD) we assume steady cruise conditions. The details of the fluid and structural models are provided in the following paragraphs. 

\subsubsection{Fluid model}

The steady-state transonic flow over the ONERA M6 wing at Mach 0.8395 and angle of attack of 3.06 degrees is computed using non-linear Euler equations. A tetrahedral grid composed of 582,752 total elements and 108,396 nodes is used for the inviscid simulation. Figure \ref{fig:oneraM6} demonstrates a close-up view of the unstructured CFD surface mesh of the wing. The boundary conditions for the computational domain are the following: Euler slip condition on the wing surface, a symmetry plane to reflect the flow about the wing root to mimic the effect of the full wing planform, and a characteristic-based condition at a cubical far-field boundary. 

\subsubsection{CFD validation studies}
\label{sec:CFD validation studies}
\begin{figure}
\centering
\subfloat
{
  \includegraphics[width=0.49\textwidth,keepaspectratio]{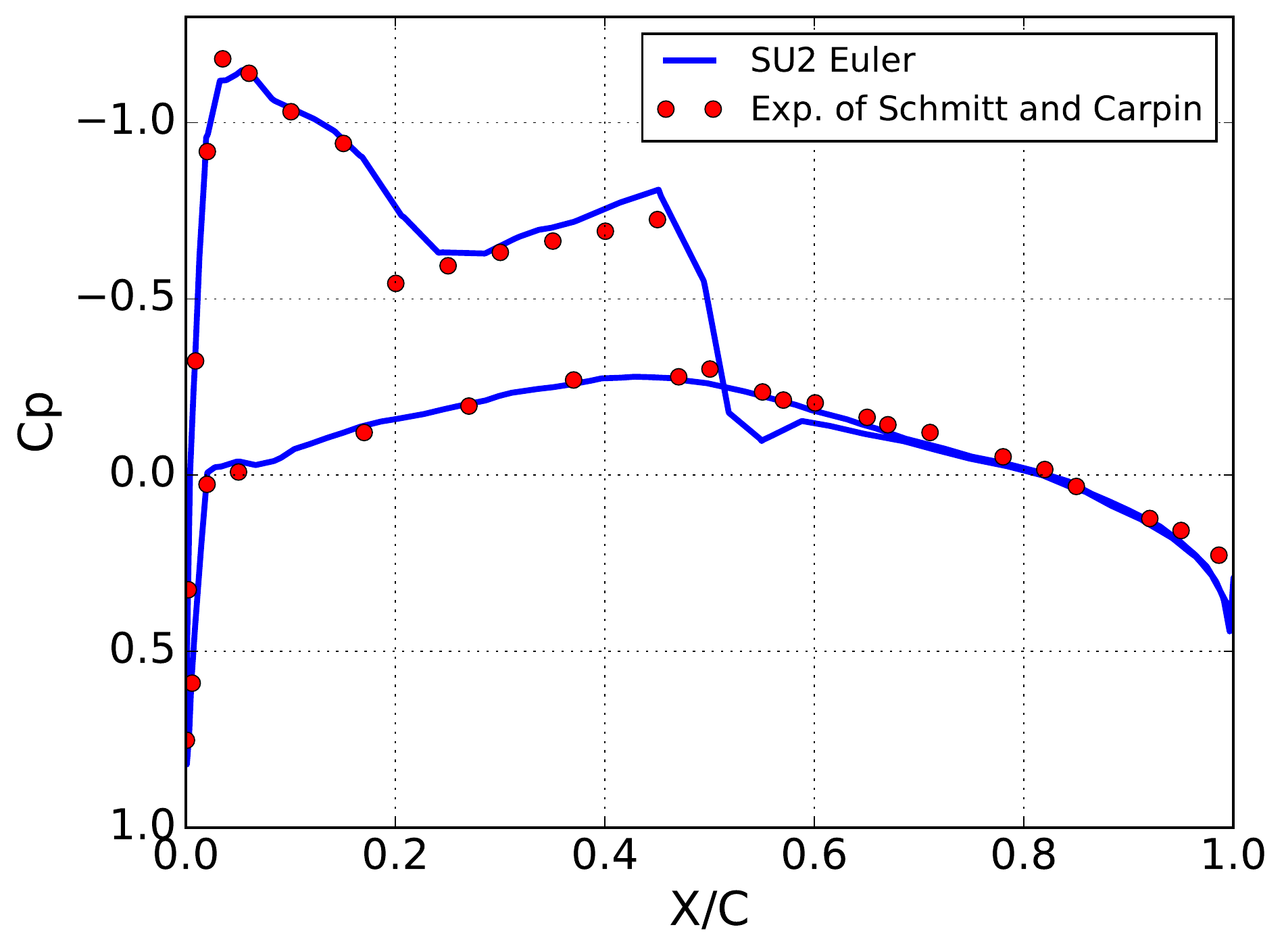}
}\\
\subfloat
{
  \includegraphics[width=0.50\textwidth,keepaspectratio]{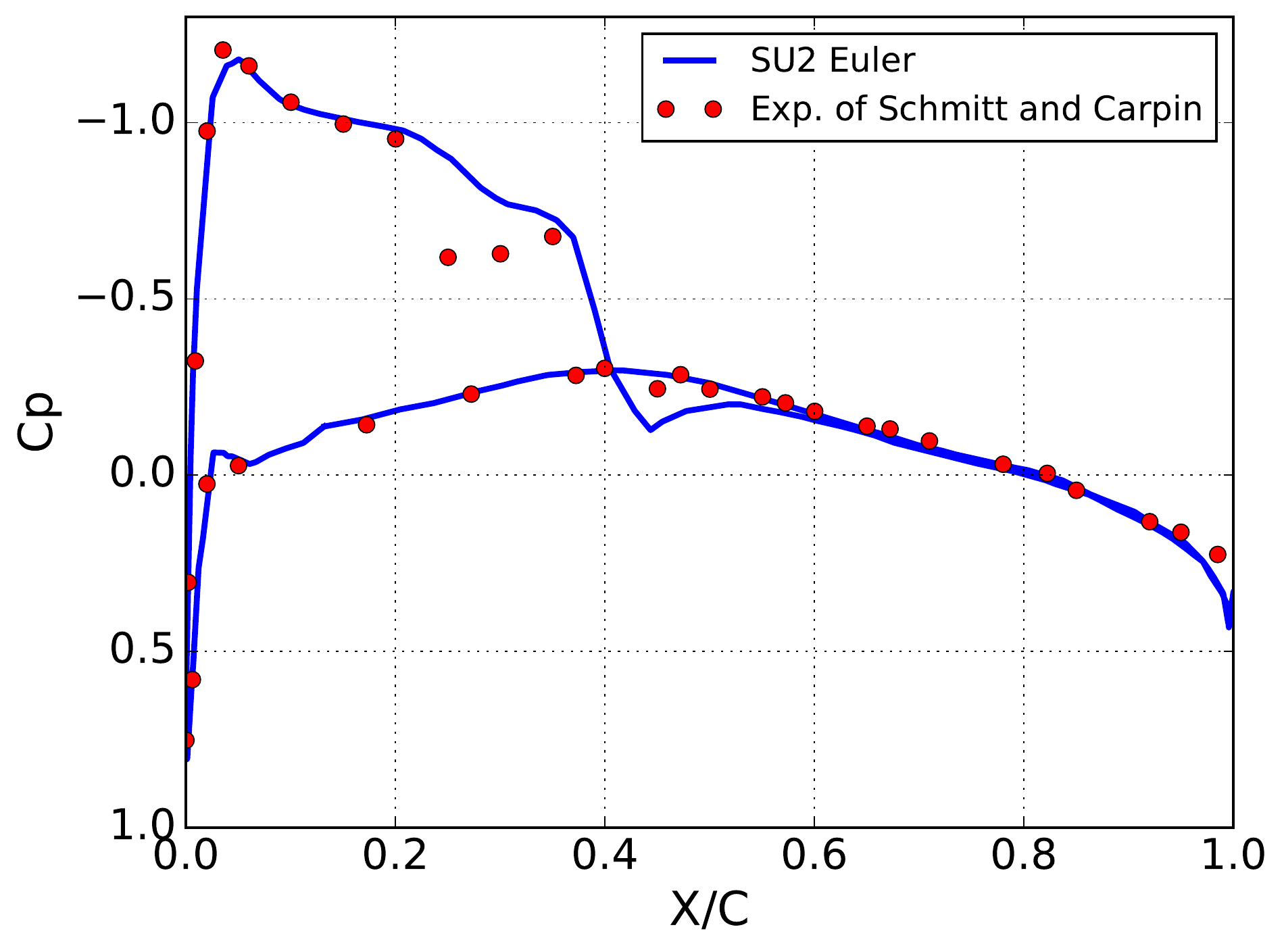}
}
\caption{Comparison of Cp profiles from the experimental results of Schmitt and Carpin (blue circles) against SU2
computational results at different sections along the span of the wing, b. Top: Y/b = 0.65, Bottom: Y/b = 0.8.}
\label{fig:ONERA_CFD_verfication}
\end{figure}

Although SU2 is comprehensively verified and validated in \cite{Palacios2014}, for the sake of completeness, direct and adjoint Euler solvers from SU2 are verified and validated against the experimental data and the central difference approximation. Assuming a rigid wing structure, Figure \ref{fig:ONERA_CFD_verfication} shows surface pressure coefficient distributions at two different span-wise stations of the wing. Overall, the computed results are in good agreement with the experimental data from \cite{Schmitt1979}, particularly along the lower surface and leading edge. Note that the flow develops strong shock in the outboard region, so the discrepancies may be attributable to the inviscid analysis.

\begin{figure}
  \centering
  \begin{tabular}{@{}c@{}}
    \includegraphics[width=0.55\linewidth,keepaspectratio]{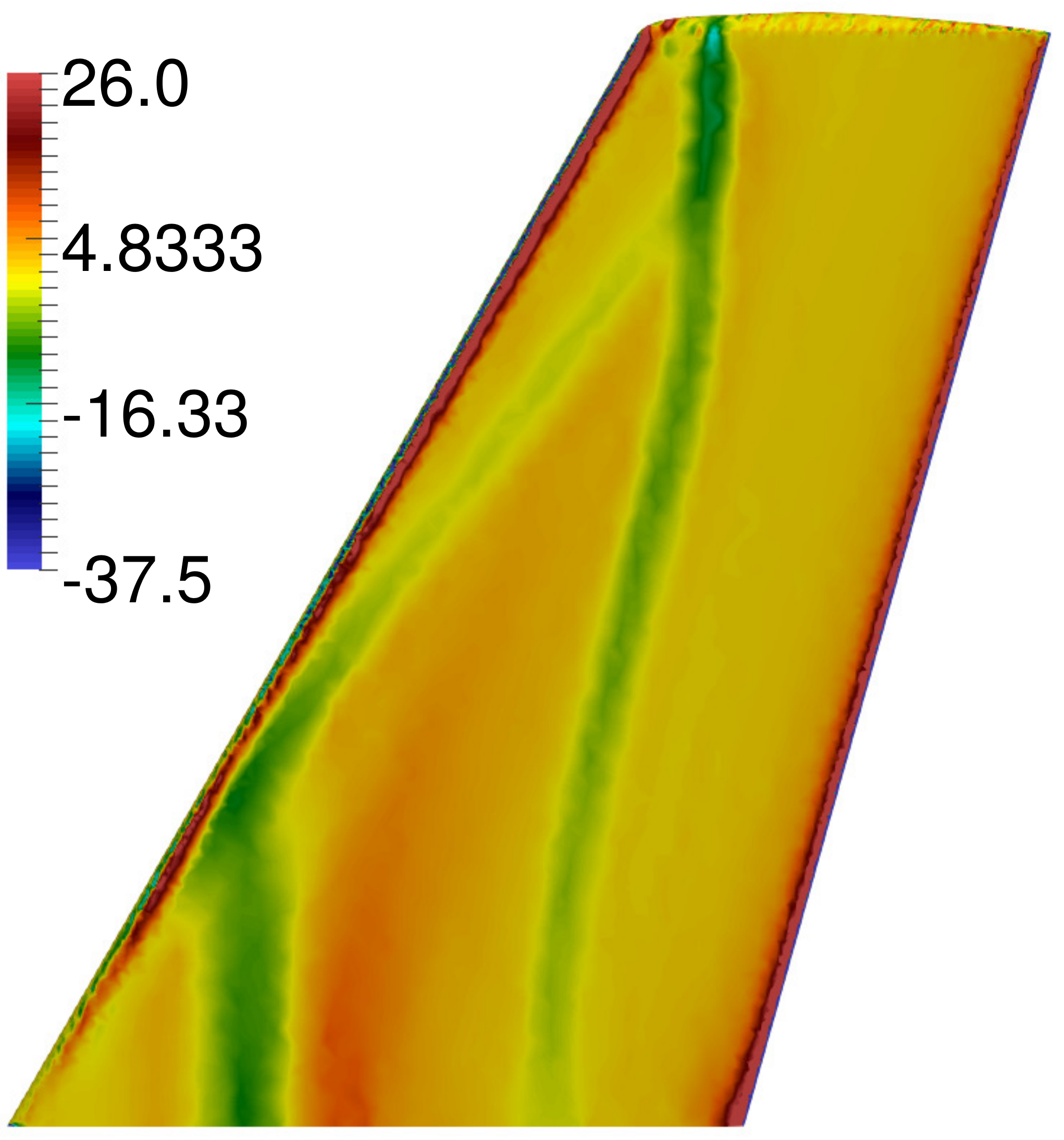} \\[\abovecaptionskip]
    \small (a) Discrete adjoint-based surface sensitivity contour \\ for a drag objective function (upper surface).
  \end{tabular}

  \vspace{\floatsep}

  \begin{tabular}{@{}c@{}}
    \includegraphics[width=\linewidth,keepaspectratio]{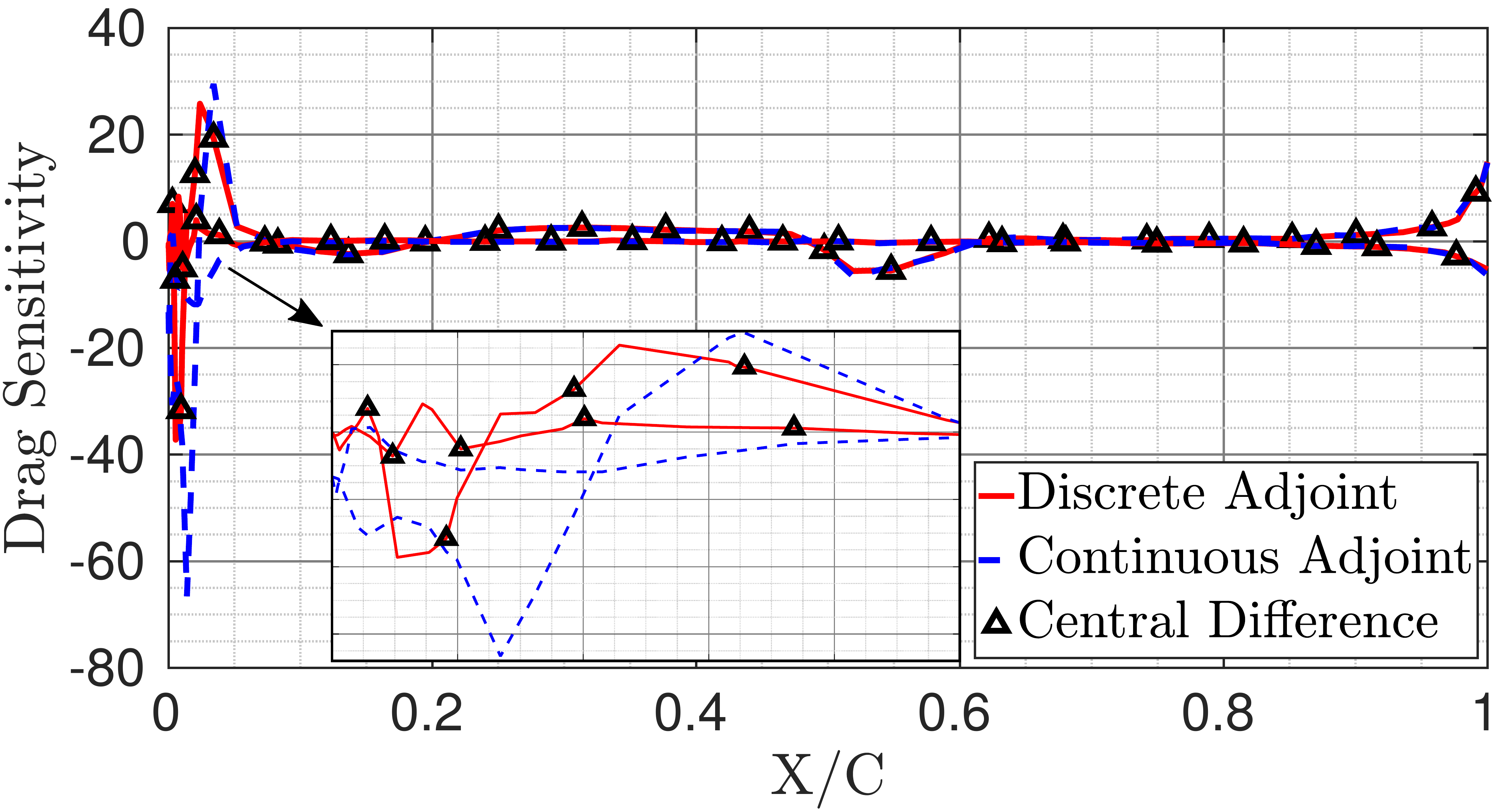} \\[\abovecaptionskip]
    \small (b) Comparison of the gradients at the section at \\ $0.65$ of the span.
  \end{tabular}

  \caption{Surface sensitivity and validation studies for a rigid ONERA M6 wing (upper surface). Perturbations and sensitivities are projected onto the surface normal.}\label{fig:CFD_Sens_verification}
\end{figure}

Having in mind that the adjoint sensitivities of the fluid appear in both the coupled adjoint mesh motion problem (Eq. \ref{e:mesh_motion_coupled_opt_problem}) and the coupled shape sensitivities (Eqs. \ref{e:partitioned_Discrete_Adj_Sens_Eq},\ref{e:partitioned_Discrete_Adj_Sens_Eq_simplified}), assessment of the  accuracy of the sensitivity information obtained by the adjoint fluid solver becomes an increasingly important part of the validation of the proposed scheme. In Figure \ref{fig:CFD_Sens_verification}, the continuous and discrete adjoint-based gradients of the wing drag with respect to a set of the grid points lying at the wing span station Y/b = 0.65 are benchmarked against finite-difference approximations and excellent agreement is observed for the AD-based discrete adjoint method. Since the continuous adjoint solver of SU2 is based on the reduced/boundary formulation (ignoring the interior mesh dependencies), the computed shape gradients are in qualitative agreement with the reference. Based on these observations, the AD-based discrete adjoint solver is used for the upcoming numerical comparisons. We also refer the readers to \cite{Economon2014,AlSaGa2016,Palacios2013,Economon2016} for details about the derivation and implementation of the adjoint solvers in SU2.

\subsubsection{Structural model}
The wing structure is modelled as a solid using 4-node tetrahedral non-linear solid elements which allow the wing to undergo large deformations. For the purposes of the following studies, two finite-element meshes of the wing are used: First, an unstructured grid which consists of 28,627 nodes and  113,096 elements with a nonmatching interface discretization, Figure \ref{fig:oneraM6}. Second, an unstructured grid which has a matching interface with the fluid mesh and it consists of 30,569 nodes and 123,245 elements. The later mesh 
serves to verify the coupled aero-structural sensitivities since it removes the mapping error at the interface, while the former 
is used for the assessment of the mapping algorithms and criteria for non-matching meshes in the FSI and adjoint FSI. It is assumed that the wing undergoes large deformations and it is made of hyperelastic material characterized by a Young’s modulus $E^{\mathcal{S}}=6\,\text{GPa}$ and Poisson ratio $\nu^{\mathcal{S}}=0.3$.

\subsubsection{Steady-state aeroelastic analysis}
\label{Steady_state_aeroelastic_analysis_of_ONERA_M6_wing}
Having set up the fluid and structural models in the baseline configuration, the steady-state aeroelastic solution was achieved by applying the primal coupling conditions to the individual domains as boundary conditions, the so-called Dirichlet-Neumann partitioning (Algorithm \ref{steady_state_FSI_algorithm}). In the following, we compare aeroelastic performance metrics of the flexible ONERA M6 wing involving matching and non-matching interface meshes. Among the various spatial mapping algorithms for surface meshes, the nearest element interpolation (NE) and the mortar method, which have similar formulation and popularity in practice, are employed in this work..

\begin{figure}
  \centering
  \begin{tabular}{@{}c@{}}
    \includegraphics[width=0.53\linewidth,keepaspectratio]{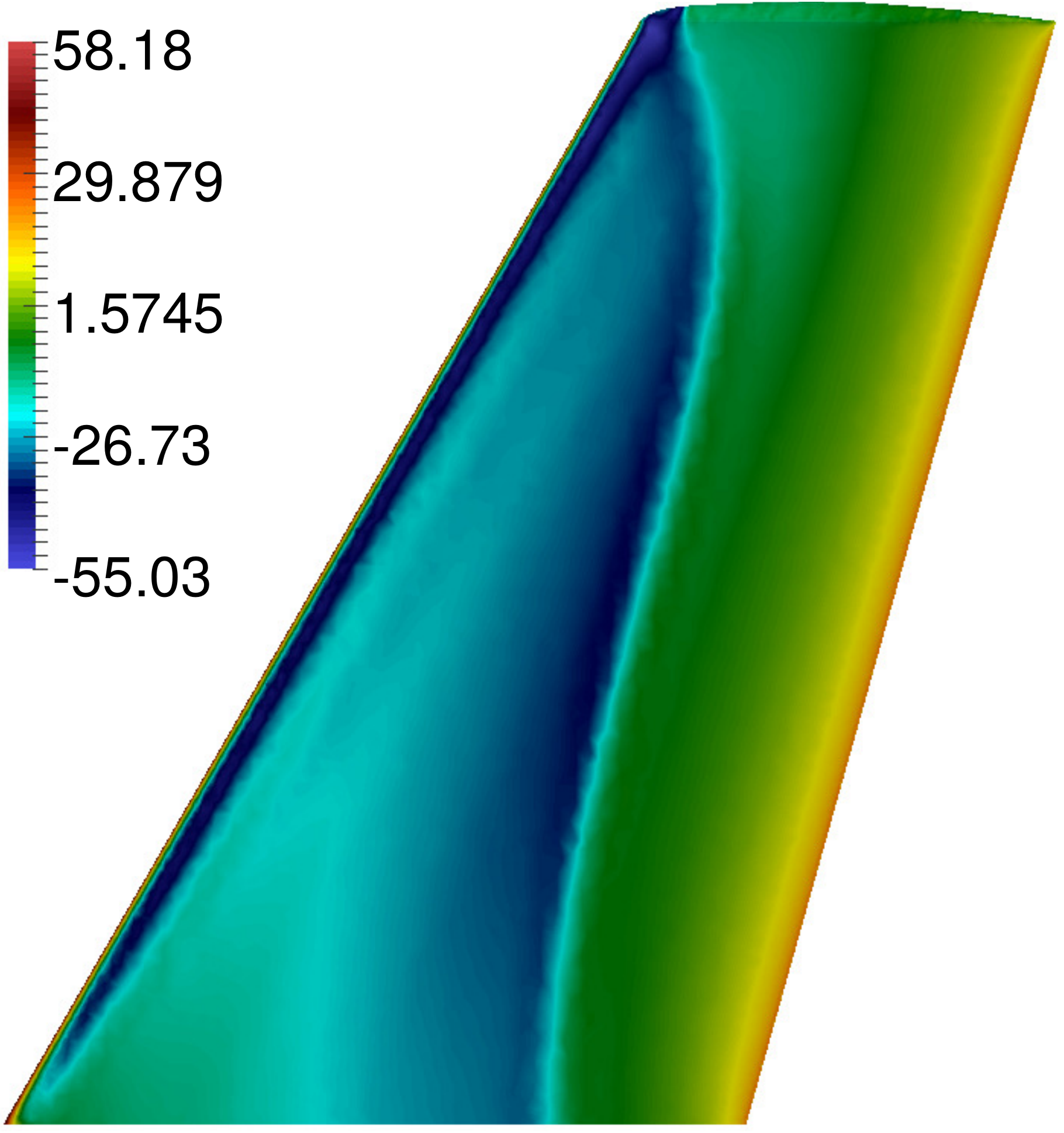} \\[\abovecaptionskip]
    \small (a) Matching interfaces.
  \end{tabular}

  \vspace{\floatsep}

  \begin{tabular}{@{}c@{}}
    \includegraphics[width=0.53\linewidth,keepaspectratio]{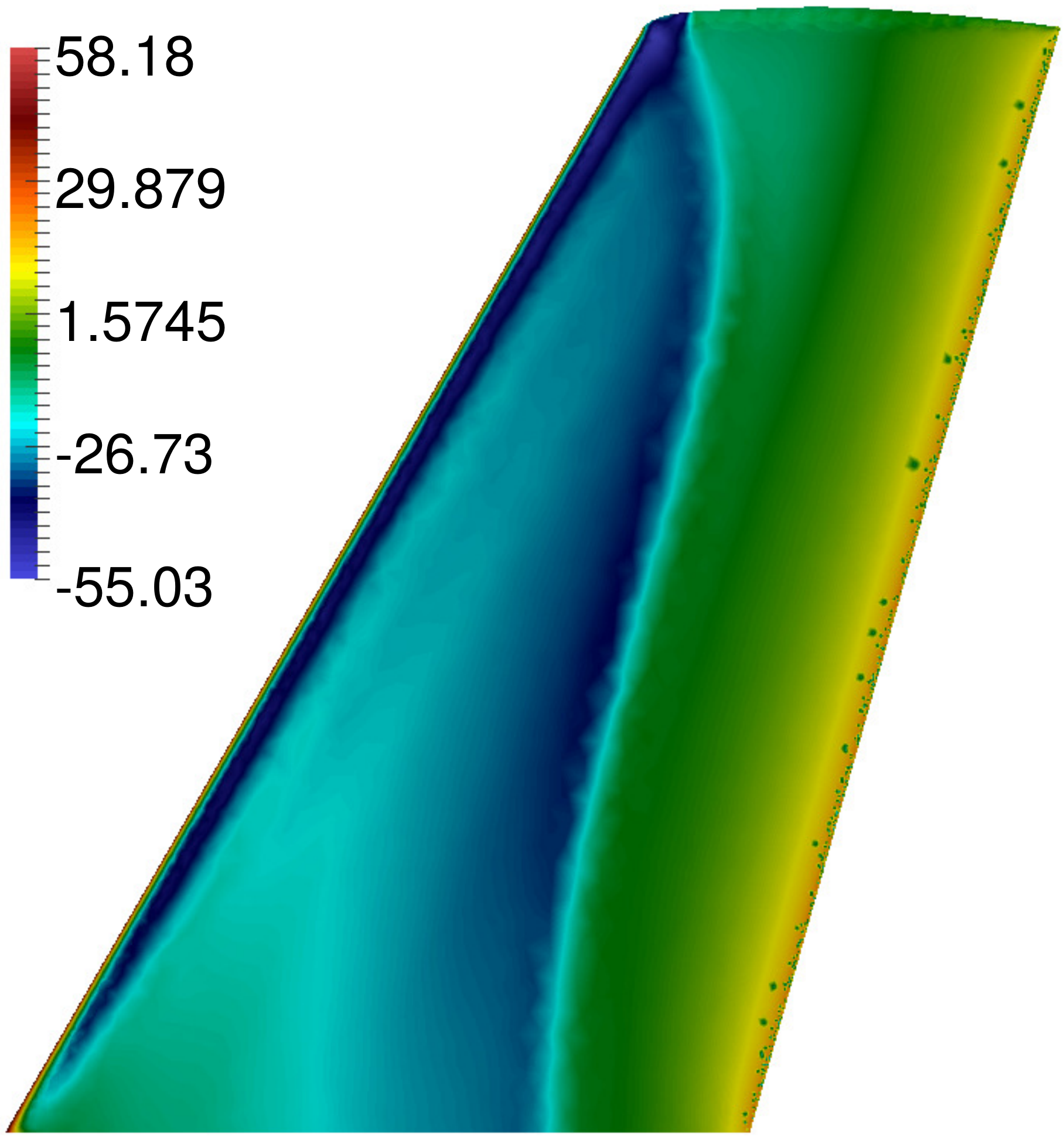} \\[\abovecaptionskip]
    \small (b) Direct mapping with nearest element interpolation.   
  \end{tabular}
  
   \vspace{\floatsep}

  \begin{tabular}{@{}c@{}}
    \includegraphics[width=0.53\linewidth,keepaspectratio]{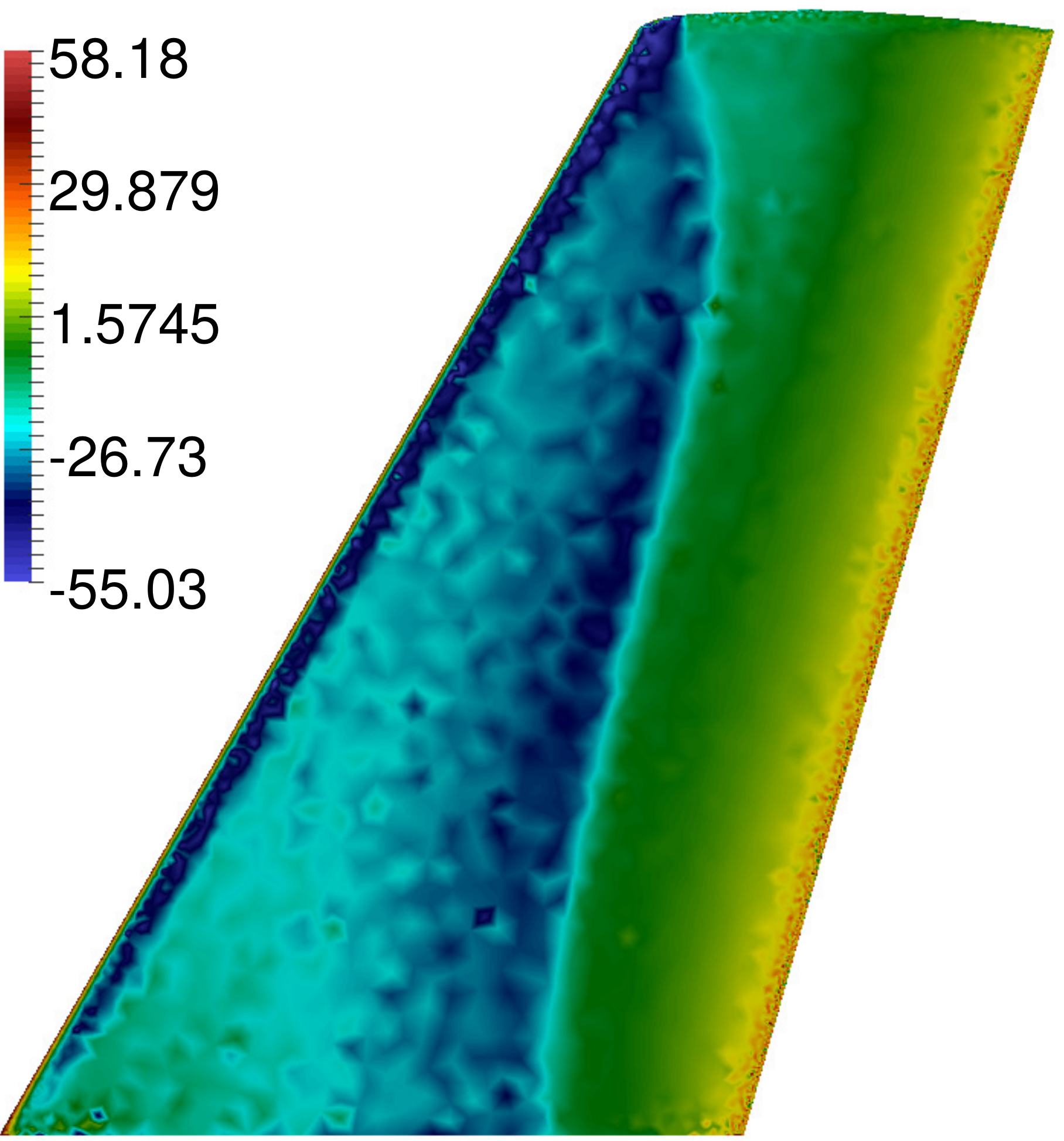} \\[\abovecaptionskip]
    \small (c) Conservative mapping with nearest element interpolation.
  \end{tabular}
  
   \vspace{\floatsep}

  \begin{tabular}{@{}c@{}}
    \includegraphics[width=0.53\linewidth,keepaspectratio]{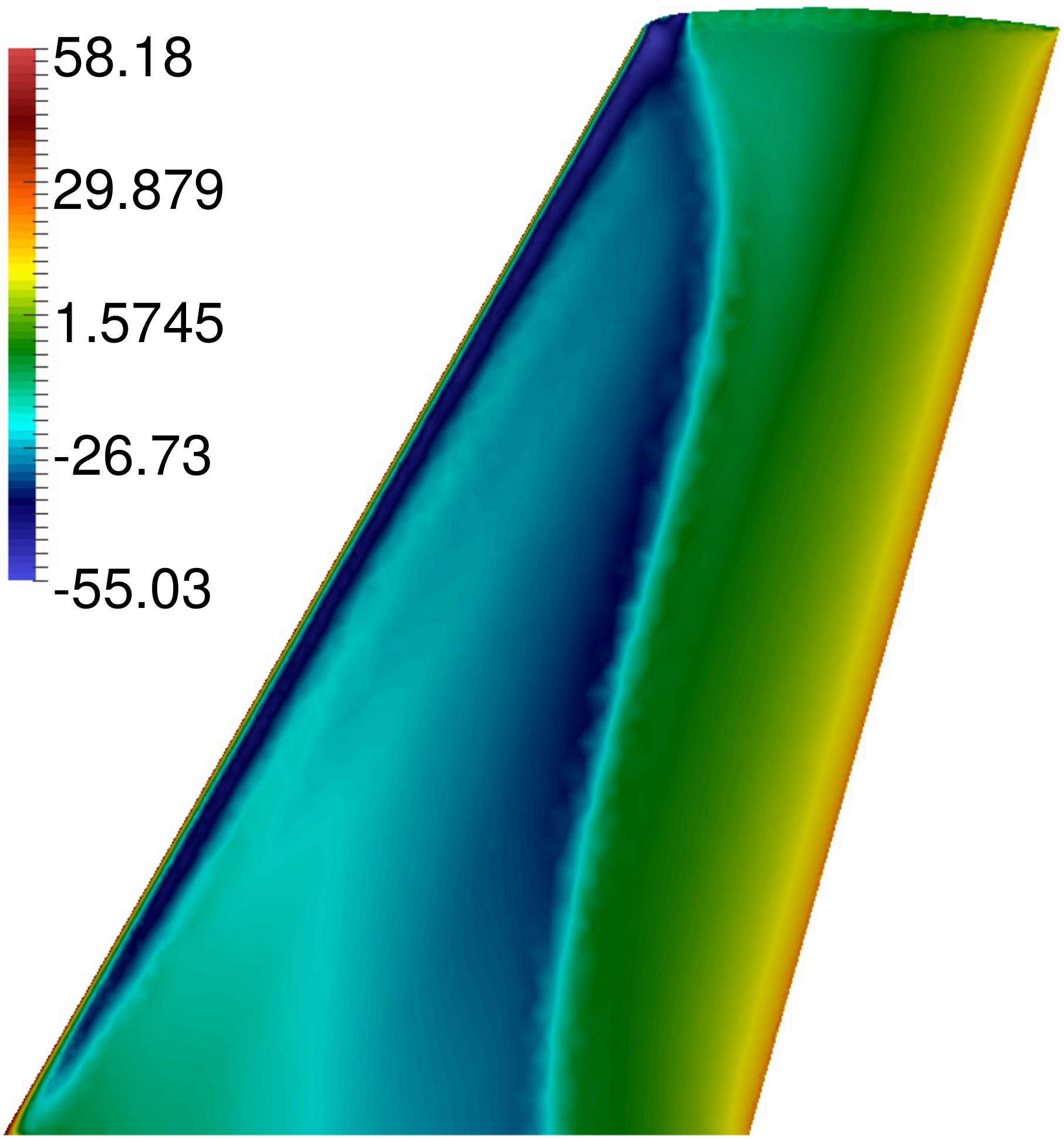} \\[\abovecaptionskip]
    \small (d) Conservative mapping with enhanced mortar method.
  \end{tabular}

  \caption{Interface pressure field (kPa) on the structure mesh for the flexible ONERA M6 wing. The results are shown for matching and non-matching interfaces using different mapping techniques.}\label{fig:ONERA_M6_FSI_traction_compariosn}
\end{figure}

\begin{figure}
  \centering
  \begin{tabular}{@{}c@{}}
    \includegraphics[width=\linewidth,keepaspectratio]{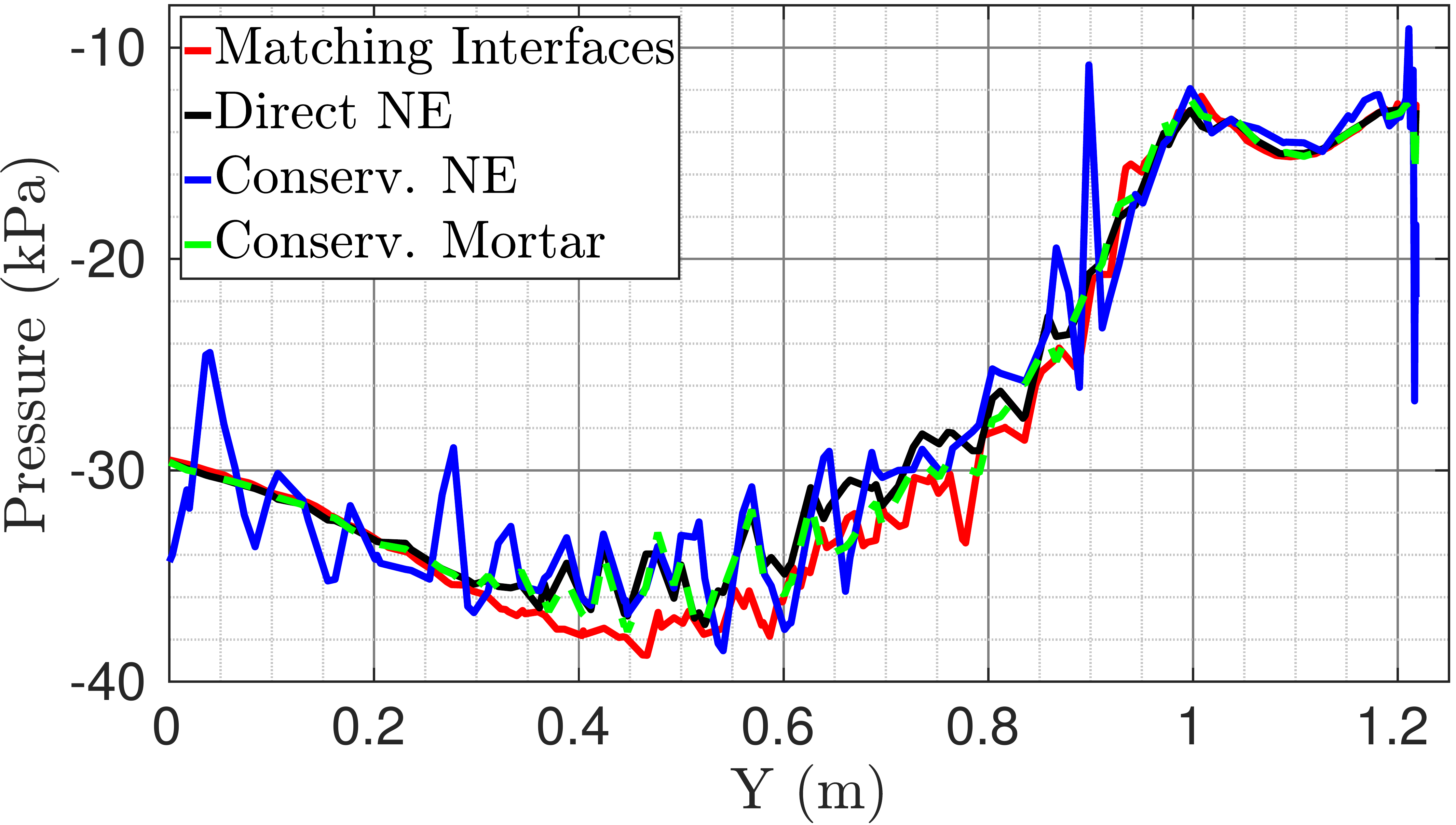} \\[\abovecaptionskip]
  \end{tabular}

  \vspace{\floatsep}

  \begin{tabular}{@{}c@{}}
    \includegraphics[width=\linewidth,keepaspectratio]{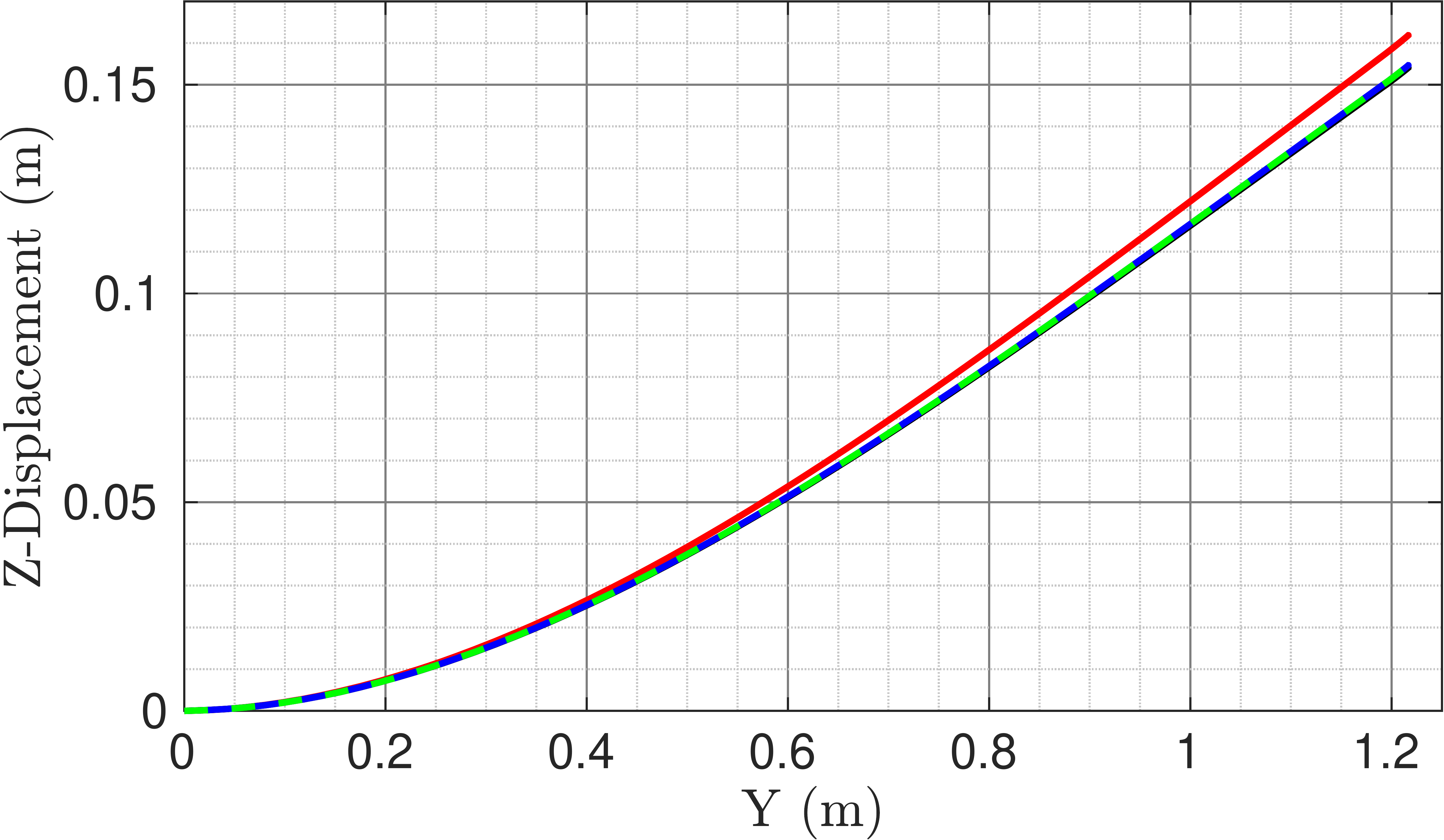} \\[\abovecaptionskip]
  \end{tabular}

  \caption{Spanwise pressure and Z-displacement fields of the upper surface in Y-Z plane at X = 0.5. NE denotes the nearest element interpolation.}\label{fig:ONERA_M6_FSI_displacement_compariosn}
\end{figure}

First, the static aeroelastic analysis is performed with the matching discrete interfaces and its solution is taken as reference for the upcoming numerical comparisons. From the undeformed state the Block Gauss-Seidel (BGS) method for the strong coupling took 24 iteration steps until the equilibrium state is reached. It was also observed that the drag and lift coefficients vary from 0.011739 and 0.286269 to 0.00502 and 0.1813, respectively, through the transition from the undeformed state to the equilibrium state. The wing-tip displacement computed with this analysis is 0.1693 meters, which is $14.15 \% $ of the span.

In the context of the spatial coupling for non-matching meshes in FSI, the direct use of the mapping algorithms is referred to as consistent mapping while using a mapping operator derived from the energy conservation is called conservative mapping. 
Consistency is an essential and basic property of the mapping algorithms which ensures that a constant field is mapped exactly.
While displacements are usually mapped using a direct/consistent mapping, surface forces/tractions are mapped either directly or conservatively (see \cite{Wang2016} and references therein). For example, the conservative fluid force transfer with NE produces nonphysical oscillations, whereas the direct mode provides accurate results. On the other hand, conservative displacement-force transfer with the mortar method delivers oscillation-free traction field on the structure interface. This behavior is linked to the weak enforcement of the coupling conditions in \cite{DeBoer2008}.    

Figure \ref{fig:ONERA_M6_FSI_traction_compariosn} shows the interface traction field on the structure mesh for the conforming and non-conforming interfaces. For the sake of quantitative comparison and completeness, Figure \ref{fig:ONERA_M6_FSI_displacement_compariosn} illustrates the interface pressure and Z-displacement distributions at a selected section in spanwise direction. Note that the plotted fields are evaluated at the static aero-elastic equilibrium of the wing.
As seen in the figures and also reported in earlier works \citep{Wang2016papaer,DeBoer2008}, the conservative mapping of forces with the nearest element interpolation gives nonphysical oscillations in the traction field mapped on the structure mesh, in contrast to the direct traction mapping. Nevertheless, the displacement fields computed for non-matching interface meshes with all three techniques are overlaying and overall in good agreement with the reference (i.e. the matching interface). This means that the structure is insensitive to local changes in the interface traction field, maybe due to the modelling of structure as solid.   

Finally, Table \ref{table:summary_FSI_ONERAM6} collects the quantitative results corresponding to the cases presented in Fig. \ref{fig:ONERA_M6_FSI_displacement_compariosn}. An important observation common to all non-matching simulations is that the interface energy (aeroelastic response) over the non-matching interfaces is in a very good agreement with the reference value (maximum error $\approx0.4\%$). This means that energy is not induced or lost due to the spatial coupling across the non-matching interfaces. However, aerodynamic and structural responses show maximum $4.21\%$ and $4.9\%$ deviations from the reference, respectively. Observed spatial coupling (mapping) errors can be generally reduced by mesh refinement \citep{DeBoer2008}, however it might be prohibitive in practical applications with moving boundaries.

\begin{table*}
\caption{Summary of aeroelastic metrics computed for flexible ONERA M6 with matching and non-matching interfaces.} % title of Table
\centering  % used for centering table
\resizebox{\textwidth}{!}{\begin{tabular}{c c c c c c c} % centered columns (4 columns)
\hline 
\hline \\
 \multicolumn{1}{l}{Spatial coupling type} & \multicolumn{1}{c}{Drag coefficient} & \multicolumn{1}{c}{Lift coefficient} & \multicolumn{1}{c}{\begin{minipage}{0.6in}Tip deflection, m\end{minipage}} & \multicolumn{1}{c}{\begin{minipage}{0.8in}Fluid interface energy, kN.m\end{minipage}} & \multicolumn{1}{c}{\begin{minipage}{0.75in}Structure interface energy, kN.m\end{minipage}} & \multicolumn{1}{c}{\begin{minipage}{0.7in}Number of Gauss-Seidel iterations, $n$\end{minipage}} \\
                        \hline \\
\multicolumn{1}{l}{\begin{minipage}{1.5in}(a) Matching interfaces.\end{minipage}} & \multicolumn{1}{c}{\(\displaystyle5.052e^{-3}\)}  &  \multicolumn{1}{c}{\(\displaystyle0.1813\)} &  \multicolumn{1}{c}{\(\displaystyle0.1693\)} &  \multicolumn{1}{c}{\(\displaystyle0.2880\)} &  \multicolumn{1}{c}{\(\displaystyle0.2880\)} &  \multicolumn{1}{c}{\(\displaystyle24\)} \\[0.25cm]
\multicolumn{1}{l}{\begin{minipage}{1.75in}(b) Direct mapping with nearest element interpolation.\end{minipage}} & \multicolumn{1}{c}{\(\displaystyle5.208e^{-3}\)}  &  \multicolumn{1}{c}{\(\displaystyle0.1856\)} &  \multicolumn{1}{c}{\(\displaystyle0.1609\)} &  \multicolumn{1}{c}{\(\displaystyle0.2848\)} &  \multicolumn{1}{c}{\(\displaystyle0.2867\)} &  \multicolumn{1}{c}{\(\displaystyle18\)} \\[0.3cm]
\multicolumn{1}{l}{\begin{minipage}{1.5in}(c) Conservative mapping with nearest element interpolation.\end{minipage}} & \multicolumn{1}{c}{\(\displaystyle5.265e^{-3}\)}  &  \multicolumn{1}{c}{\(\displaystyle0.1870\)} &  \multicolumn{1}{c}{\(\displaystyle0.1613\)} &  \multicolumn{1}{c}{\(\displaystyle0.2891\)} &  \multicolumn{1}{c}{\(\displaystyle0.2891\)} &  \multicolumn{1}{c}{\(\displaystyle21\)}\\[0.3cm]
\multicolumn{1}{l}{\begin{minipage}{1.5in}(d) Conservative mapping with enhanced mortar method.\end{minipage}} & \multicolumn{1}{c}{\(\displaystyle5.2566e^{-3}\)}  &  \multicolumn{1}{c}{\(\displaystyle0.1869\)} &  \multicolumn{1}{c}{\(\displaystyle0.1613\)} &  \multicolumn{1}{c}{\(\displaystyle0.2891\)} &  \multicolumn{1}{c}{\(\displaystyle0.2891\)} &  \multicolumn{1}{c}{\(\displaystyle22\)}
\end{tabular}}
\label{table:summary_FSI_ONERAM6} % is used to refer this table in the text
\end{table*}

\begin{figure}
  \centering
  \begin{tabular}{@{}c@{}}
    \includegraphics[width=0.53\linewidth,keepaspectratio]{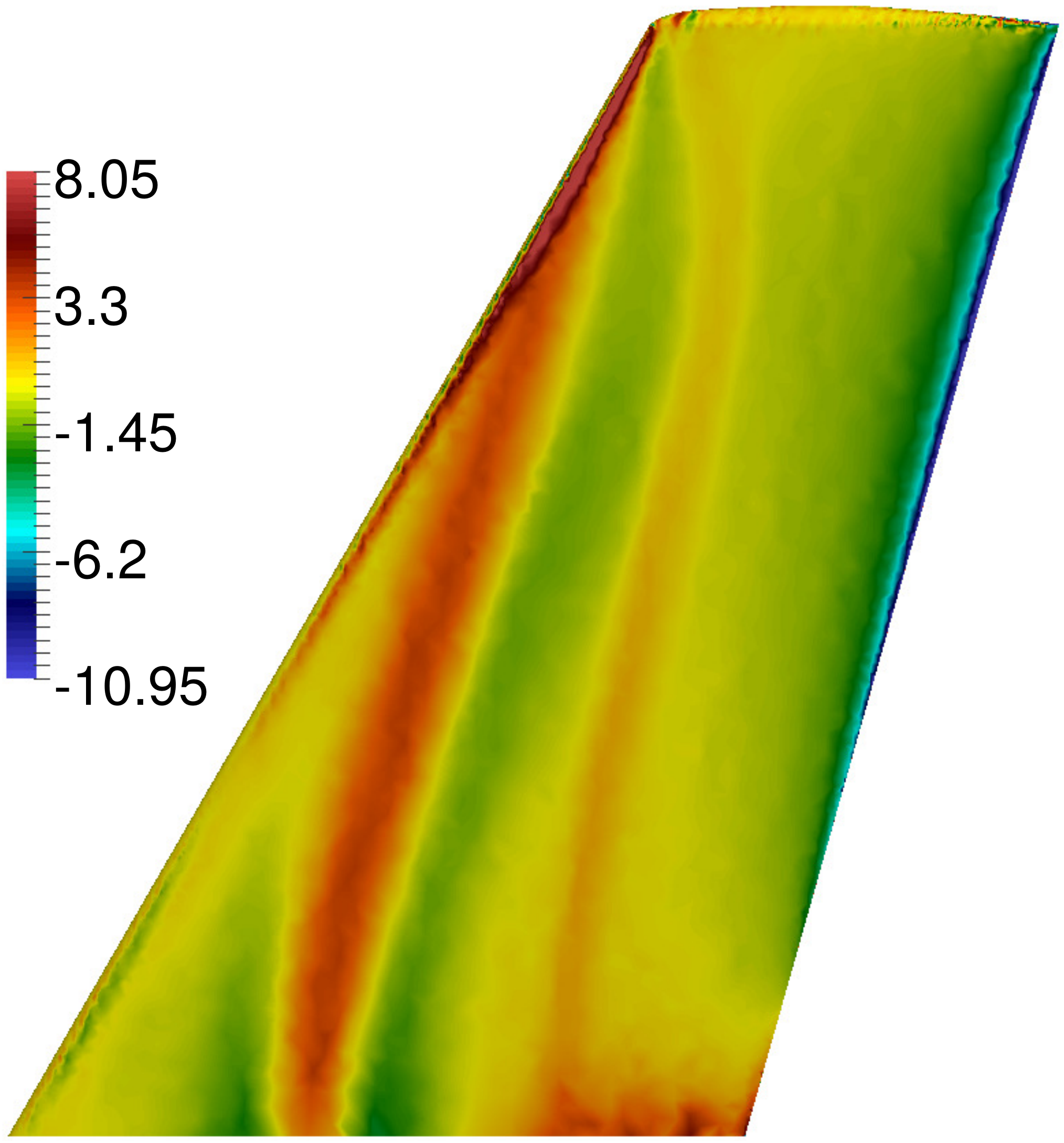} \\[\abovecaptionskip]
    \small (a) Three-field-based surface sensitivity contour of \\ the upper surface.
  \end{tabular}

  \vspace{\floatsep}

  \begin{tabular}{@{}c@{}}
    \includegraphics[width=0.92\linewidth,keepaspectratio]{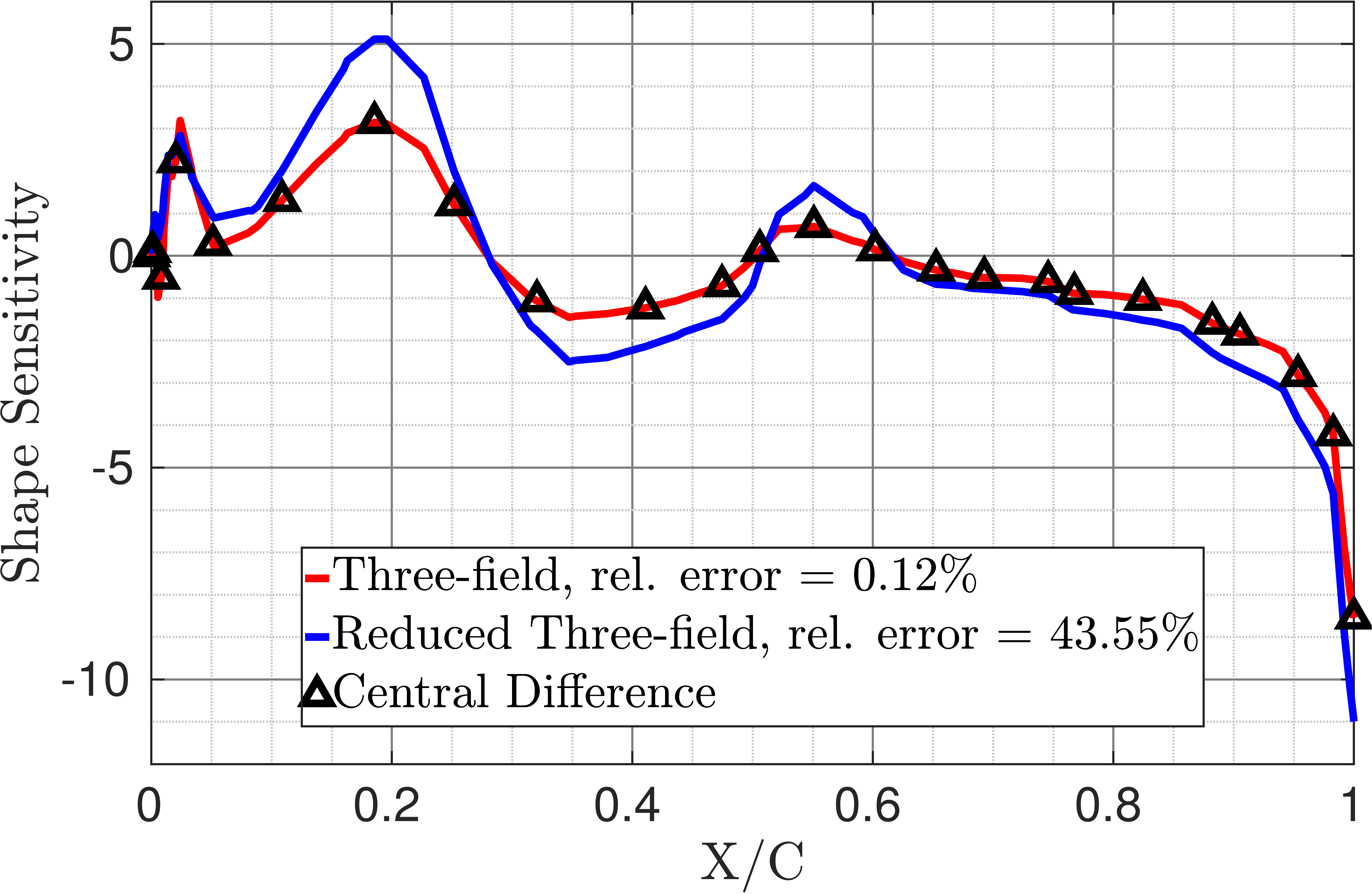} \\[\abovecaptionskip]
    \small (b) Comparison of the gradients at the section at \\ Y/b=$0.65$ (upper surface).
  \end{tabular}
  
   \vspace{\floatsep}

  \begin{tabular}{@{}c@{}}
    \includegraphics[width=0.53\linewidth,keepaspectratio]{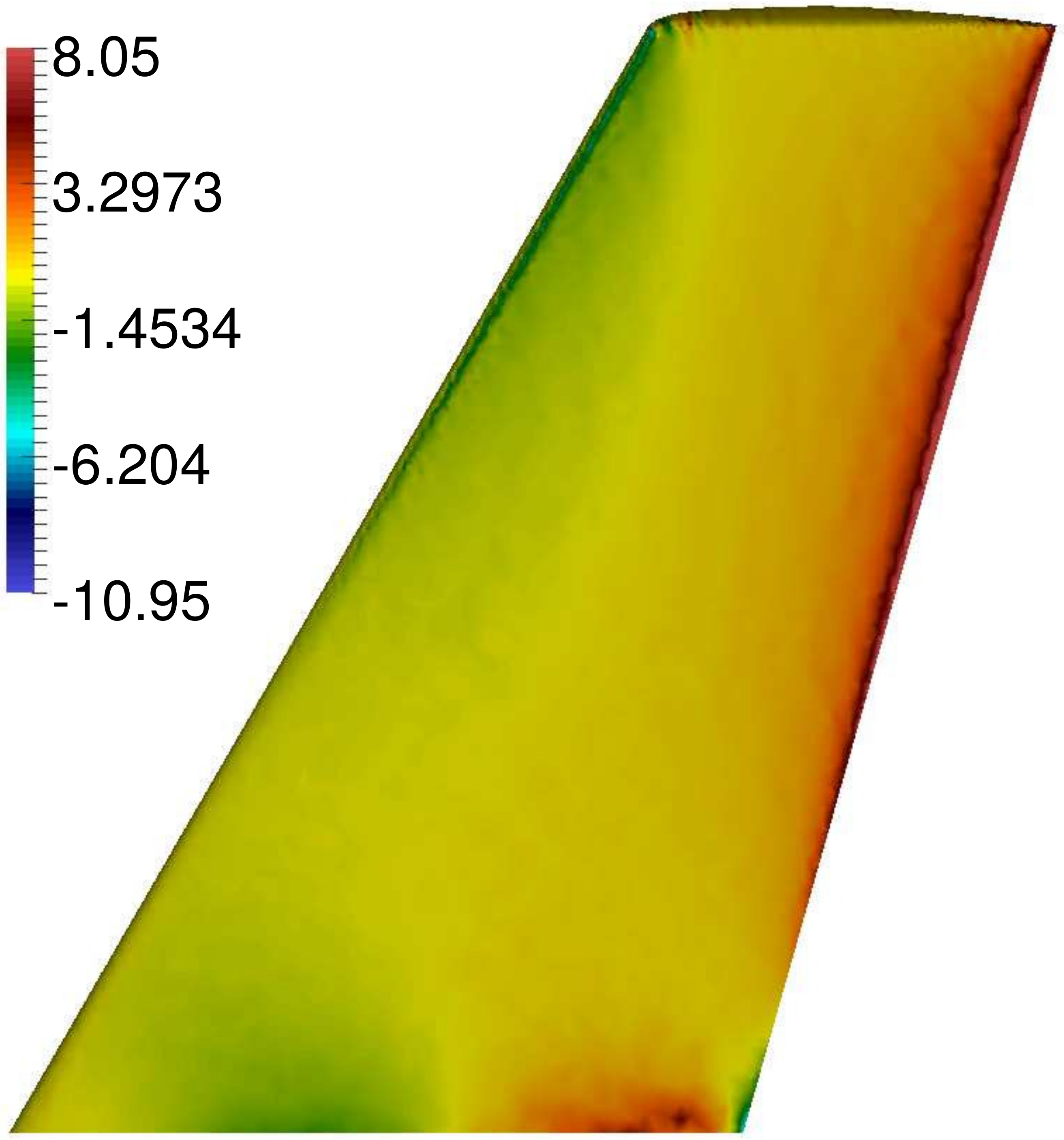} \\[\abovecaptionskip]
    \small (c) Three-field-based surface sensitivity contour of \\ the lower surface.
  \end{tabular}
  
   \vspace{\floatsep}

  \begin{tabular}{@{}c@{}}
    \includegraphics[width=0.92\linewidth,keepaspectratio]{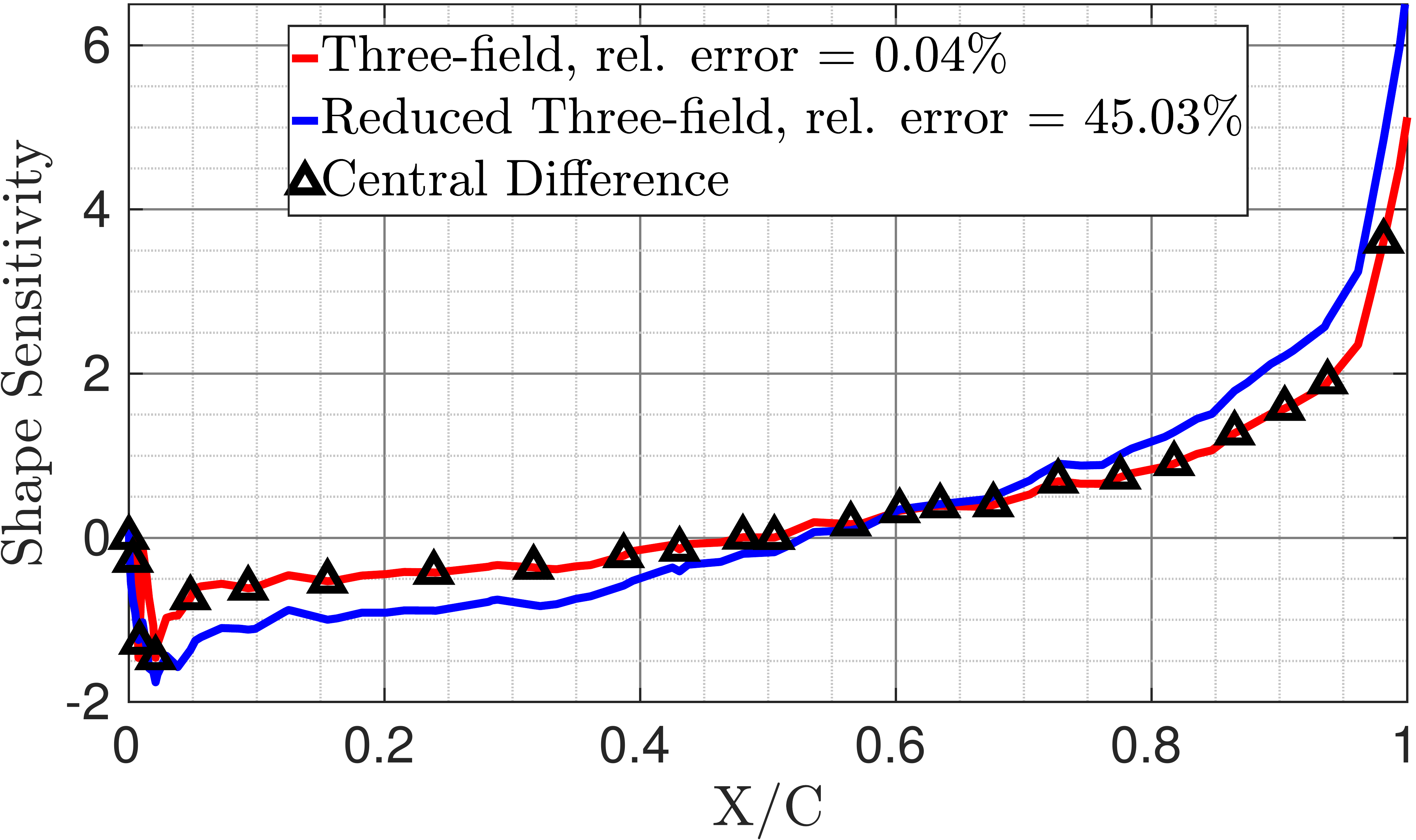} \\[\abovecaptionskip]
    \small (d) Comparison of the gradients at the section at \\ Y/b=$0.65$ (lower surface).
  \end{tabular}

  \caption{Coupled shape sensitivity analysis for the interface energy.}\label{fig:FSI_Sens_verification_ONERA_M6_strain}
\end{figure}

\begin{figure}
  \centering
  \begin{tabular}{@{}c@{}}
    \includegraphics[width=0.54\linewidth,keepaspectratio]{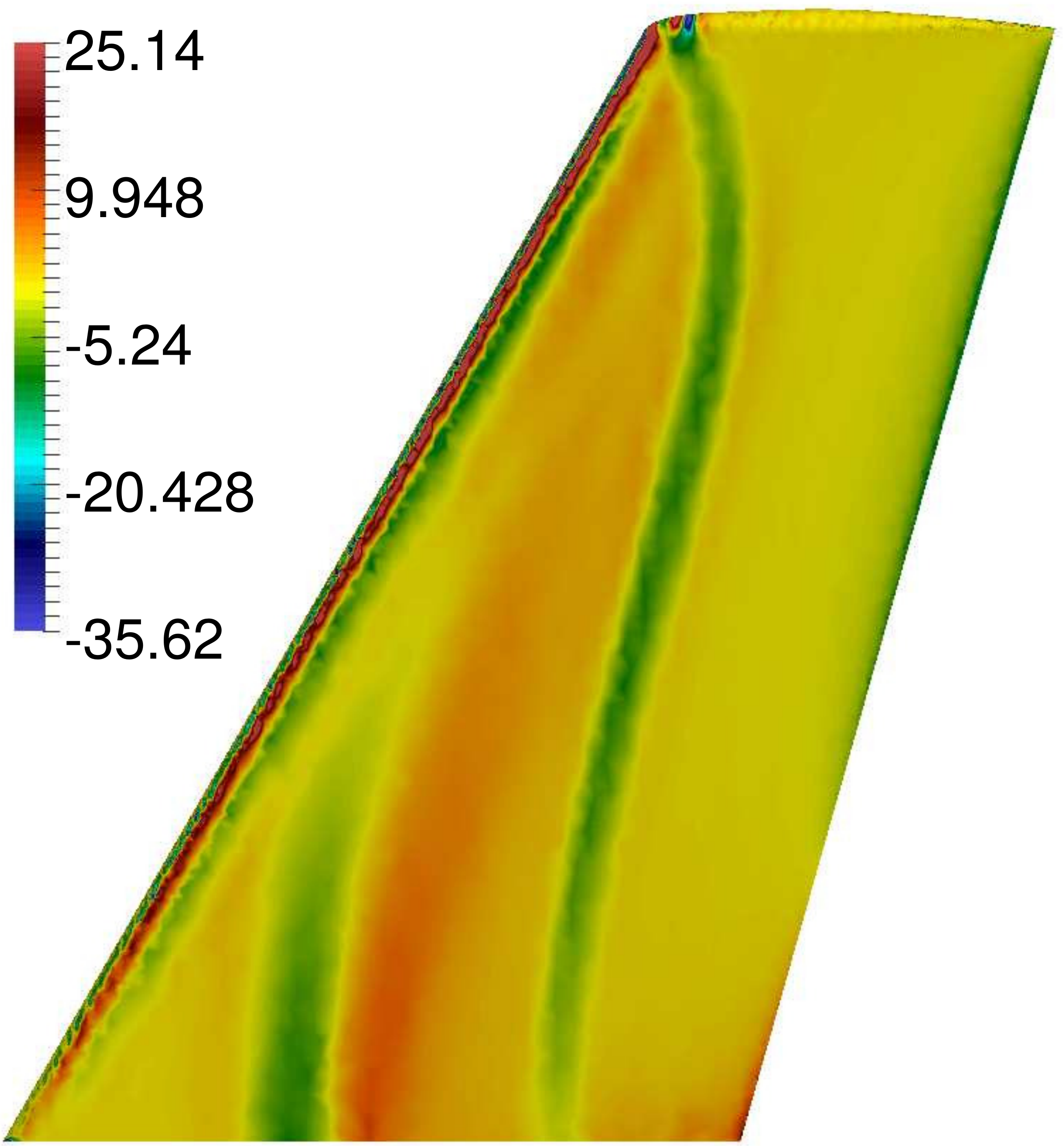} \\[\abovecaptionskip]
    \small (a) Three-field-based surface sensitivity contour of \\ the upper surface.
  \end{tabular}

  \vspace{\floatsep}

  \begin{tabular}{@{}c@{}}
    \includegraphics[width=0.9\linewidth,keepaspectratio]{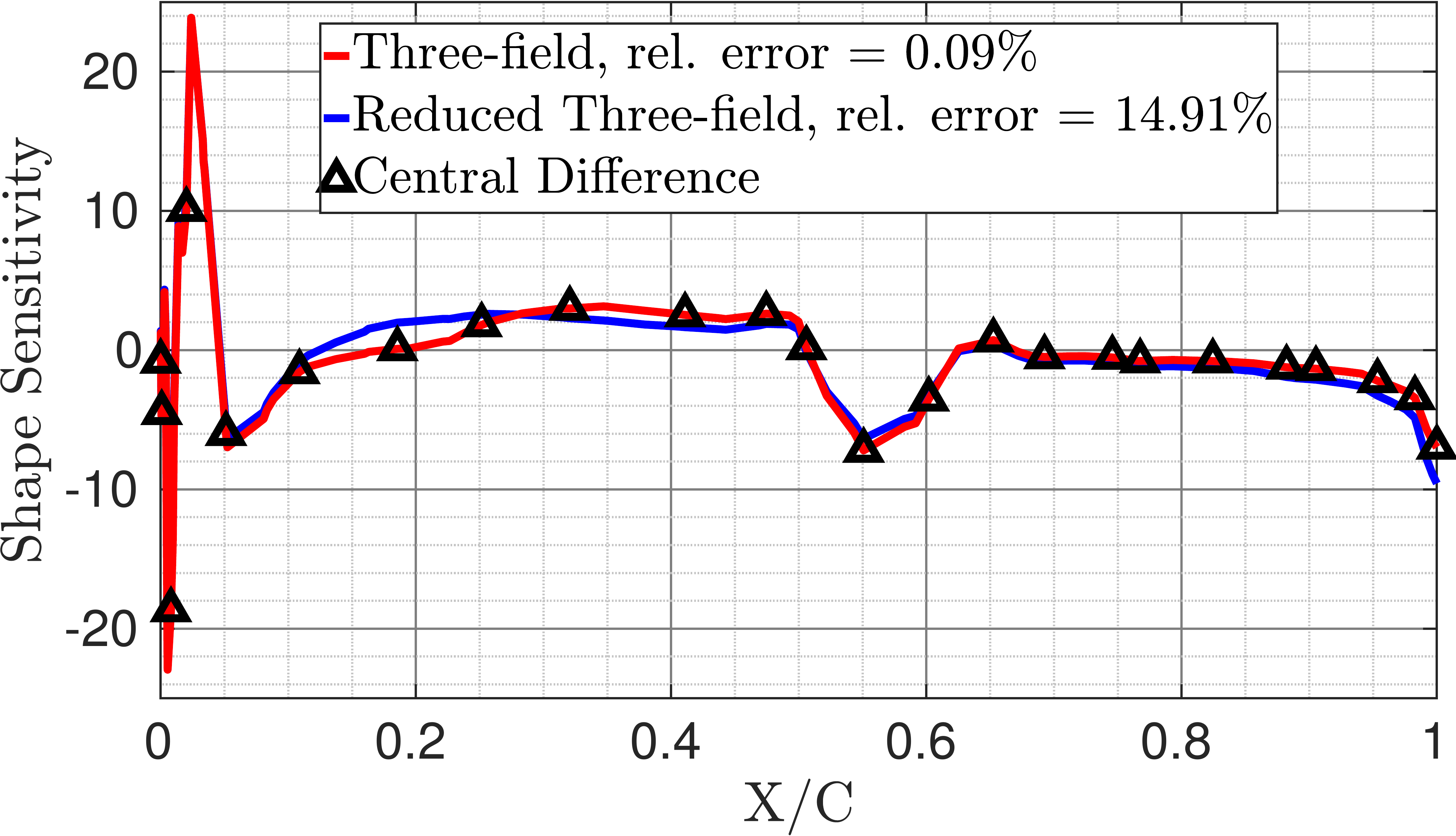} \\[\abovecaptionskip]
    \small (b) Comparison of the gradients at the section at \\ Y/b=$0.65$ (upper surface).
  \end{tabular}
  
   \vspace{\floatsep}

  \begin{tabular}{@{}c@{}}
    \includegraphics[width=0.54\linewidth,keepaspectratio]{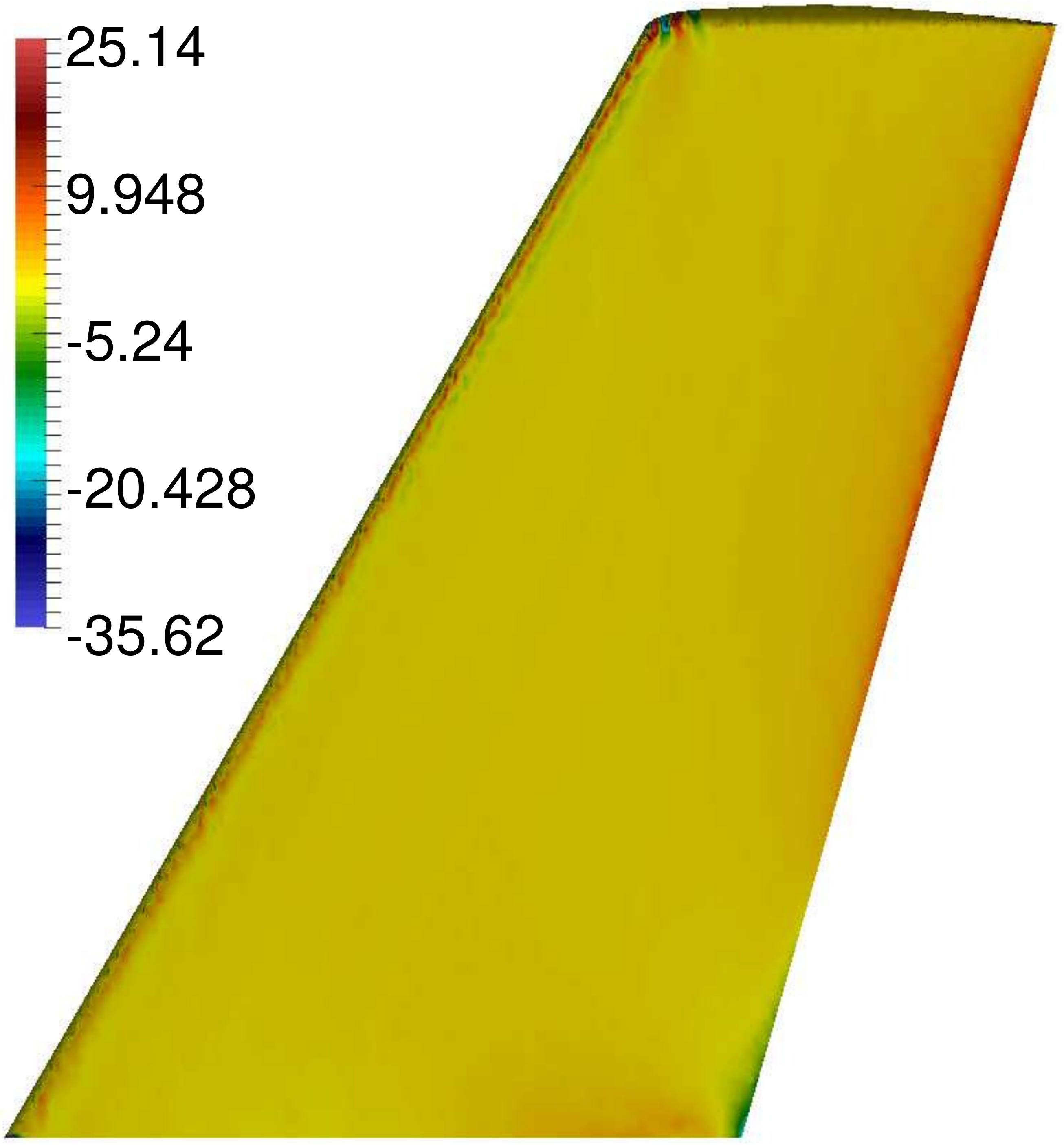} \\[\abovecaptionskip]
    \small (c) Three-field-based surface sensitivity contour of \\ the lower surface.
  \end{tabular}
  
   \vspace{\floatsep}

  \begin{tabular}{@{}c@{}}
    \includegraphics[width=0.9\linewidth,keepaspectratio]{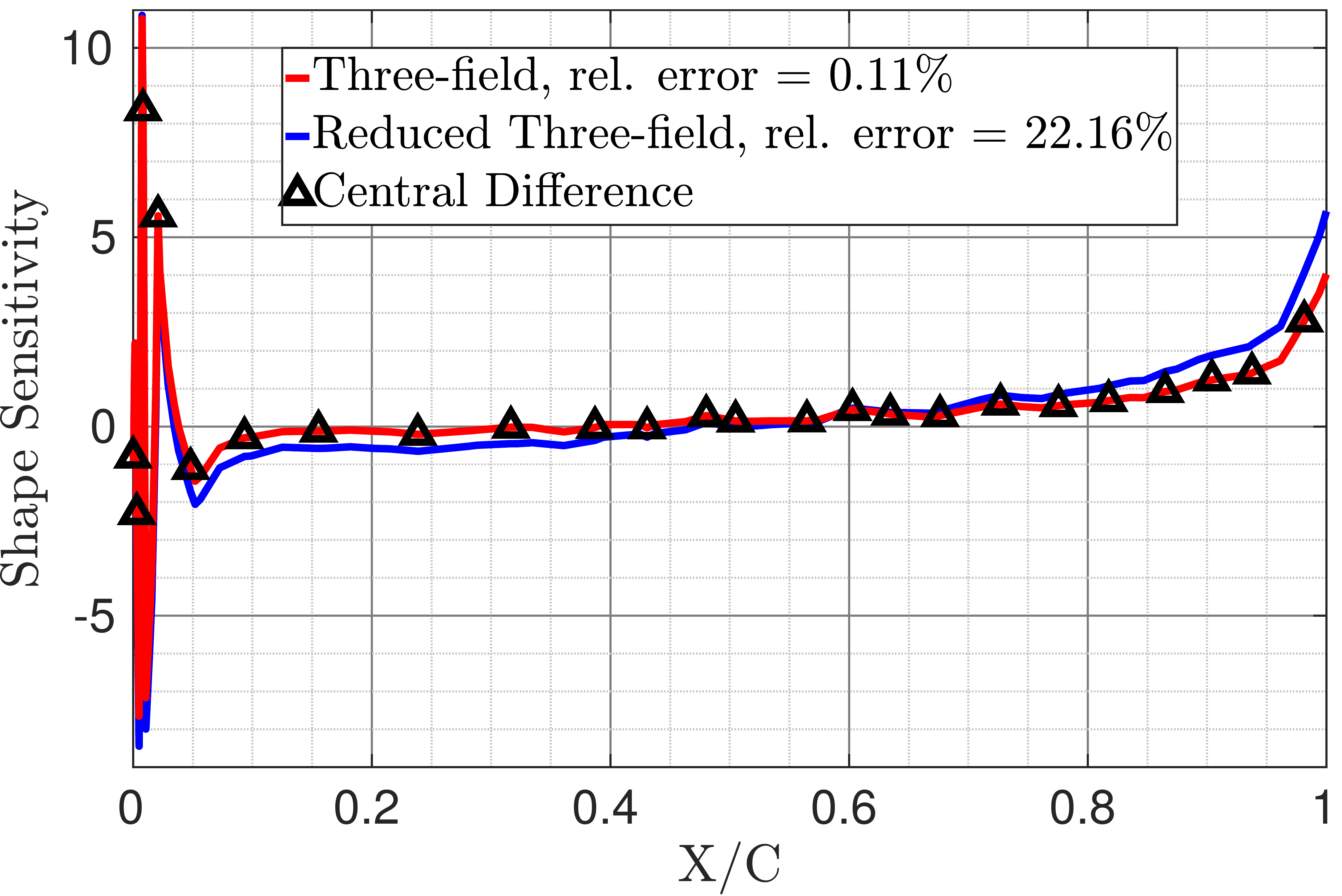} \\[\abovecaptionskip]
    \small (d) Comparison of the gradients at the section at\\  Y/b=$0.65$ (lower surface).
  \end{tabular}

  \caption{Coupled shape sensitivity analysis for the interface drag objective.}\label{fig:FSI_Sens_verification_ONERA_M6_drag}
\end{figure}
 
%\begin{figure*}
%\includegraphics[keepaspectratio,width=\textwidth]{images/AFSI_two_three_fields_formulations_comparisons/final.eps}
%\caption{Surface sensitivity contours on the upper surface of flexible ONERA M6 wing, including a comparison between three-field and two field-formulations for shape sensitivity analysis at y/b = 0.6. Left: strain energy objective. Right: drag objective}
%\label{fig:FSI_Sens_verification_ONERA_M6}
%\end{figure*}

\begin{figure}
  \includegraphics[keepaspectratio,width=0.49\textwidth]{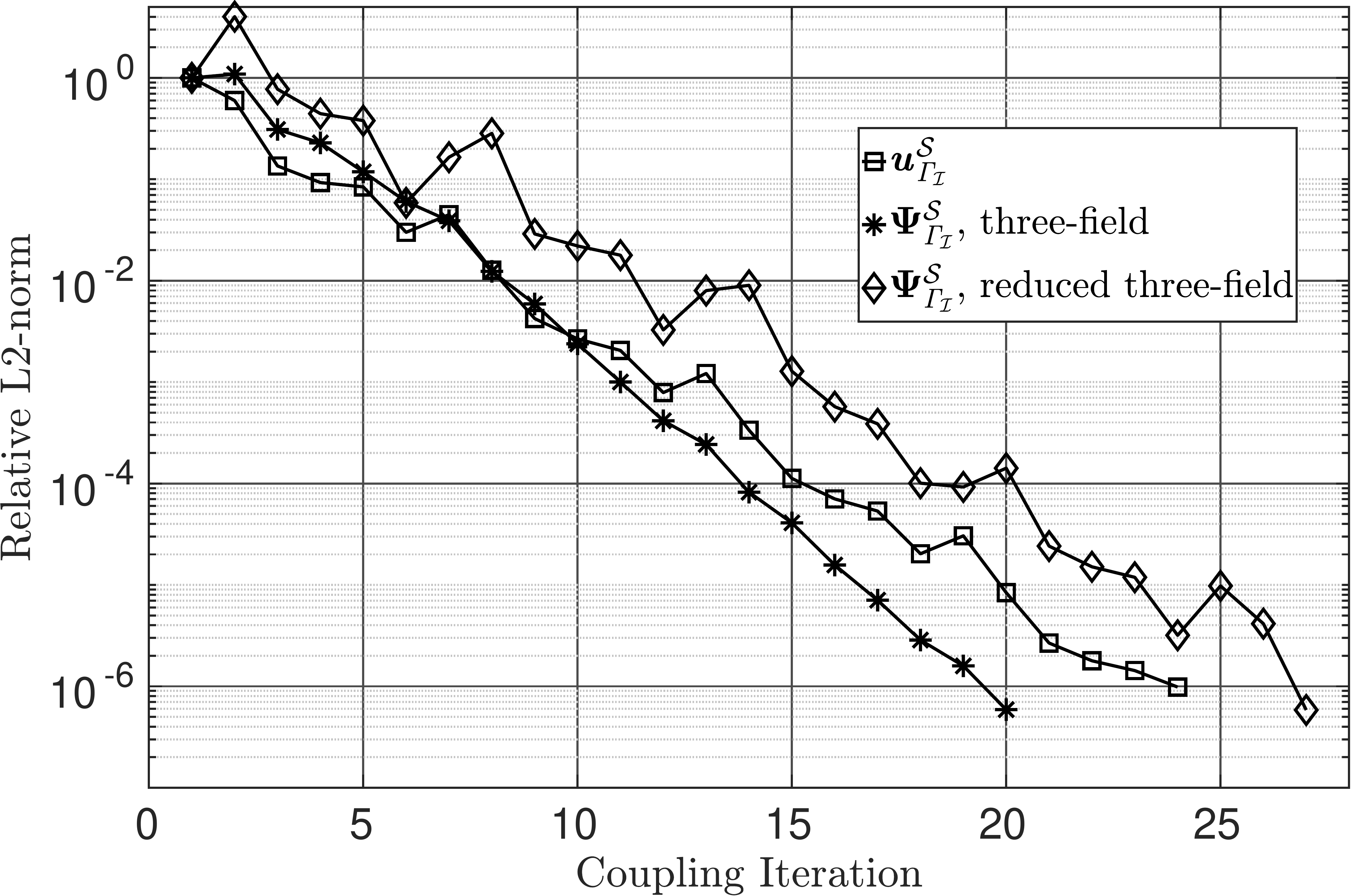}
\caption{Convergence histories of interface residuals for the direct and adjoint FSI problems.}
\label{fig:convergence_history}      
\end{figure}   

\subsubsection{Coupled shape sensitivity analysis}
This section demonstrates the applicability of the presented partitioned scheme (Fig. \ref{fig:AFSI_three_field}) to the aeroelastic shape sensitivity analysis of the flexible ONERA M6 wing, using adjoints and shape sensitivities distributed throughout different codes with specific formulations. Namely, the AD-based adjoint solver of SU2 is strongly coupled to the discrete adjoint solver of KRATOS via the coupling tool EMIPRE. The AD-based adjoint solver is chosen instead of the continuous one due to the accuracy of the computed gradients (refer to Section \ref{sec:CFD validation studies}) and the explicit availability of the domain-based shape gradients. 

As a first step, we evaluate the accuracy of the complete and reduced three-field-based shape sensitivities  against the central difference results. Considering matching interface meshes, Figures \ref{fig:FSI_Sens_verification_ONERA_M6_strain} and \ref{fig:FSI_Sens_verification_ONERA_M6_drag}, respectively, display the interface energy and drag sensitivity maps. They also compare both formulations against the reference for cross-section Y/b = 0.65. As can be seen and expected, there is a satisfactory match between the three-field-based and reference shape gradients. Moreover, similar to observations in Section \ref{FEM-based_shape_sensitivity_analysis_for_FSI}, discrepancies between the two formulations seem to be quantitative rather than qualitative. Deviations of the reduced gradient formulation are more pronounced around the shock and sharp trailing edge regions, where the validity of the reduced/boundary gradient formulation has been challenged intensively by \cite{lozano2019watch,lozano2018singular,Lozano2017}. Lozano has concluded that at sharp trailing
edges, inviscid adjoint solutions and sensitivities of force-based objectives are strongly mesh dependent and do not converge as the mesh is refined. However, this is not the case for viscous adjoint sensitivity analysis. Regarding the convergence properties, Figure \ref{fig:convergence_history} compares in a semi-logarithmic diagram the convergence histories of the interface displacement and the interface adjoint displacements for both formulations. Obviously, the complete three-field-based adjoint FSI analysis shows faster and smoother convergence behavior than the FSI and reduced three-field-based adjoint FSI analyses. This is due to the facts that FSI is a nonlinear problem whereas adjoint FSI is a linear problem, and also, unlike the reduced three-field formulation, the complete three-field approach linearizes the primal problem exactly (i.e. without any assumption). 

\begin{figure}
  \centering
  \begin{tabular}{@{}c@{}}
    \includegraphics[width=0.53\linewidth,keepaspectratio]{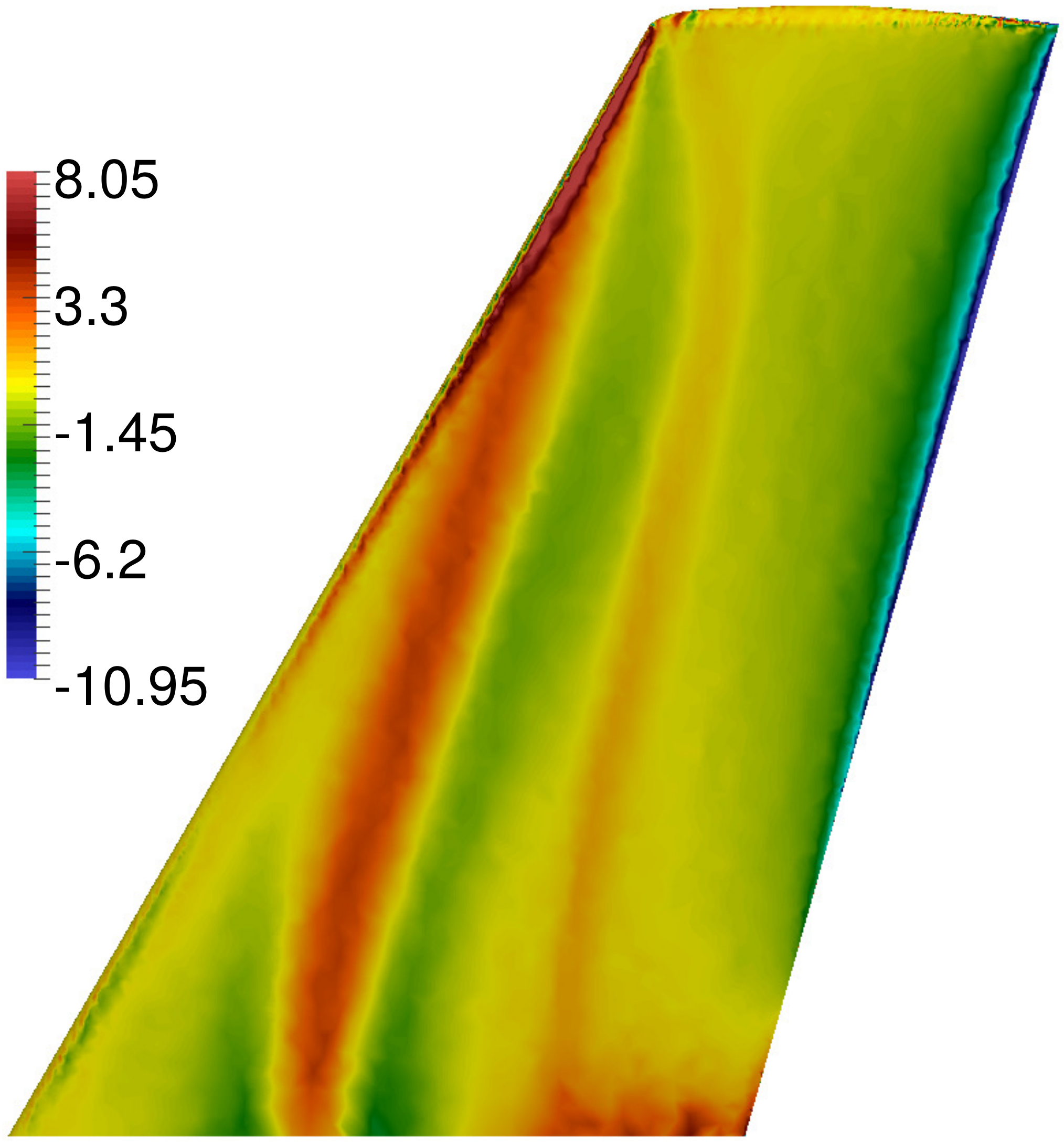} \\[\abovecaptionskip]
    \small (a) Matching interfaces.
  \end{tabular}

  \vspace{\floatsep}

  \begin{tabular}{@{}c@{}}
    \includegraphics[width=0.53\linewidth,keepaspectratio]{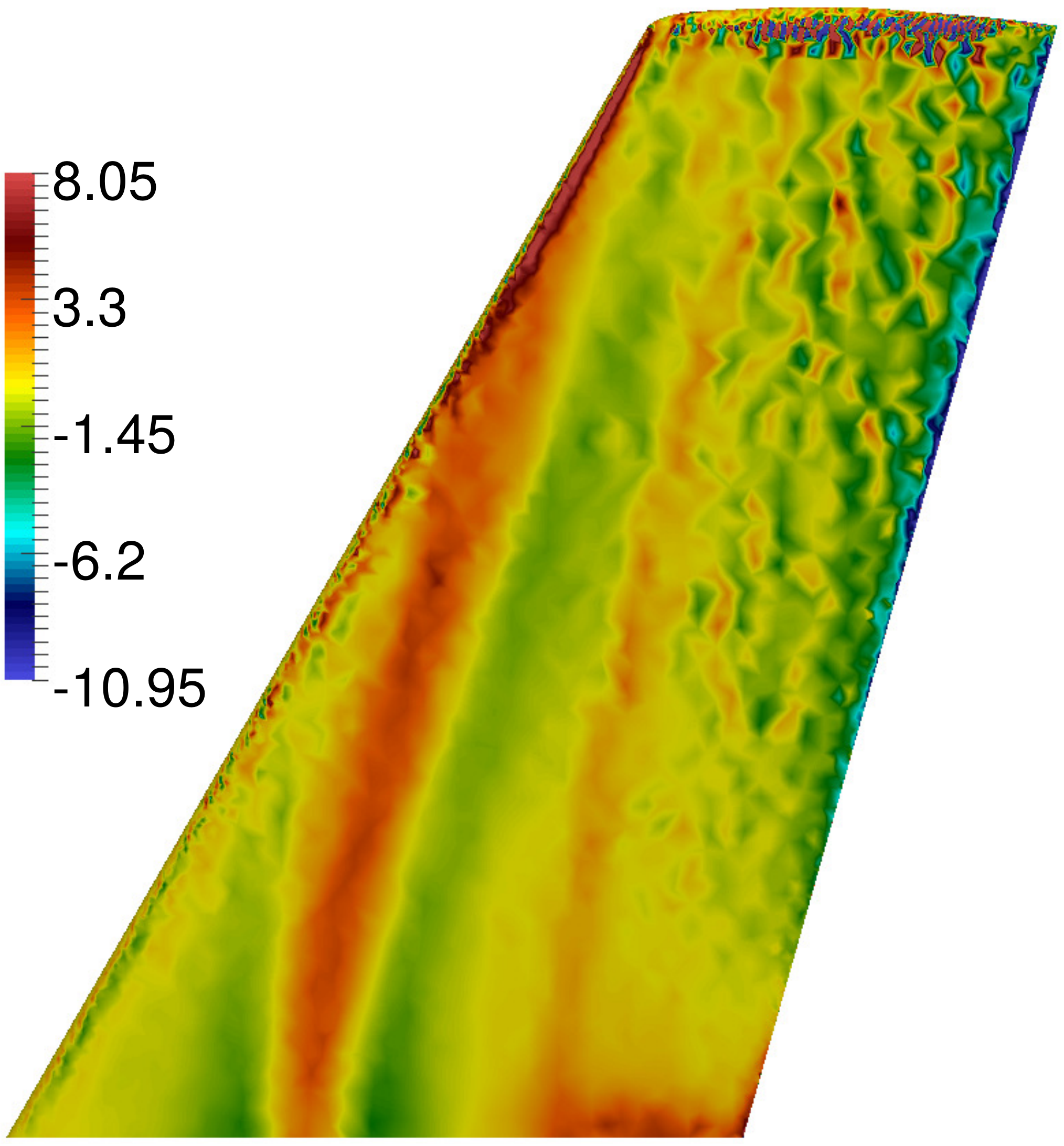} \\[\abovecaptionskip]
    \small (b) Direct mapping with nearest element interpolation.
  \end{tabular}
  
   \vspace{\floatsep}

  \begin{tabular}{@{}c@{}}
    \includegraphics[width=0.53\linewidth,keepaspectratio]{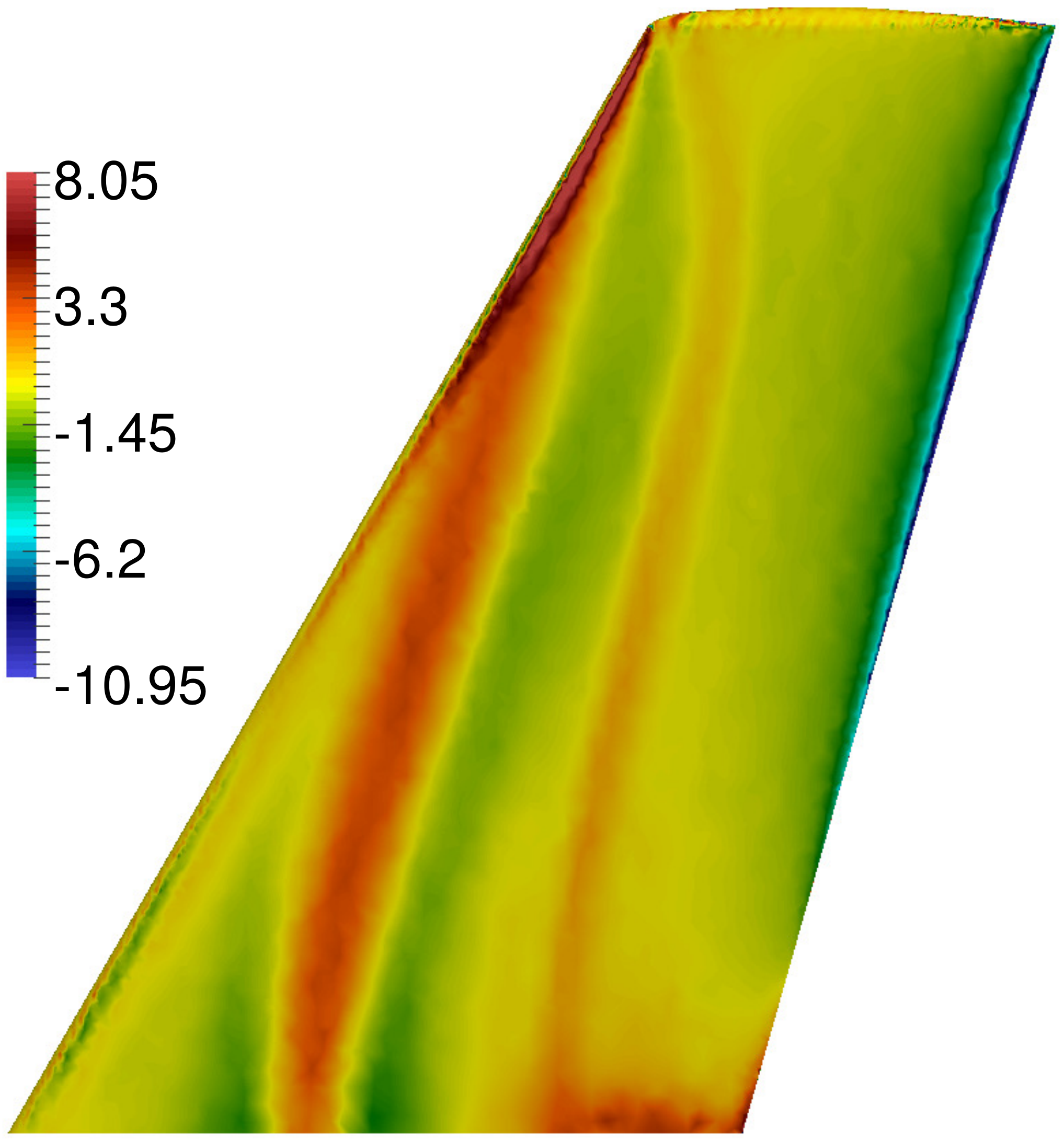} \\[\abovecaptionskip]
    \small (c) Conservative mapping with nearest element interpolation.
  \end{tabular}
  
   \vspace{\floatsep}

  \begin{tabular}{@{}c@{}}
    \includegraphics[width=0.53\linewidth,keepaspectratio]{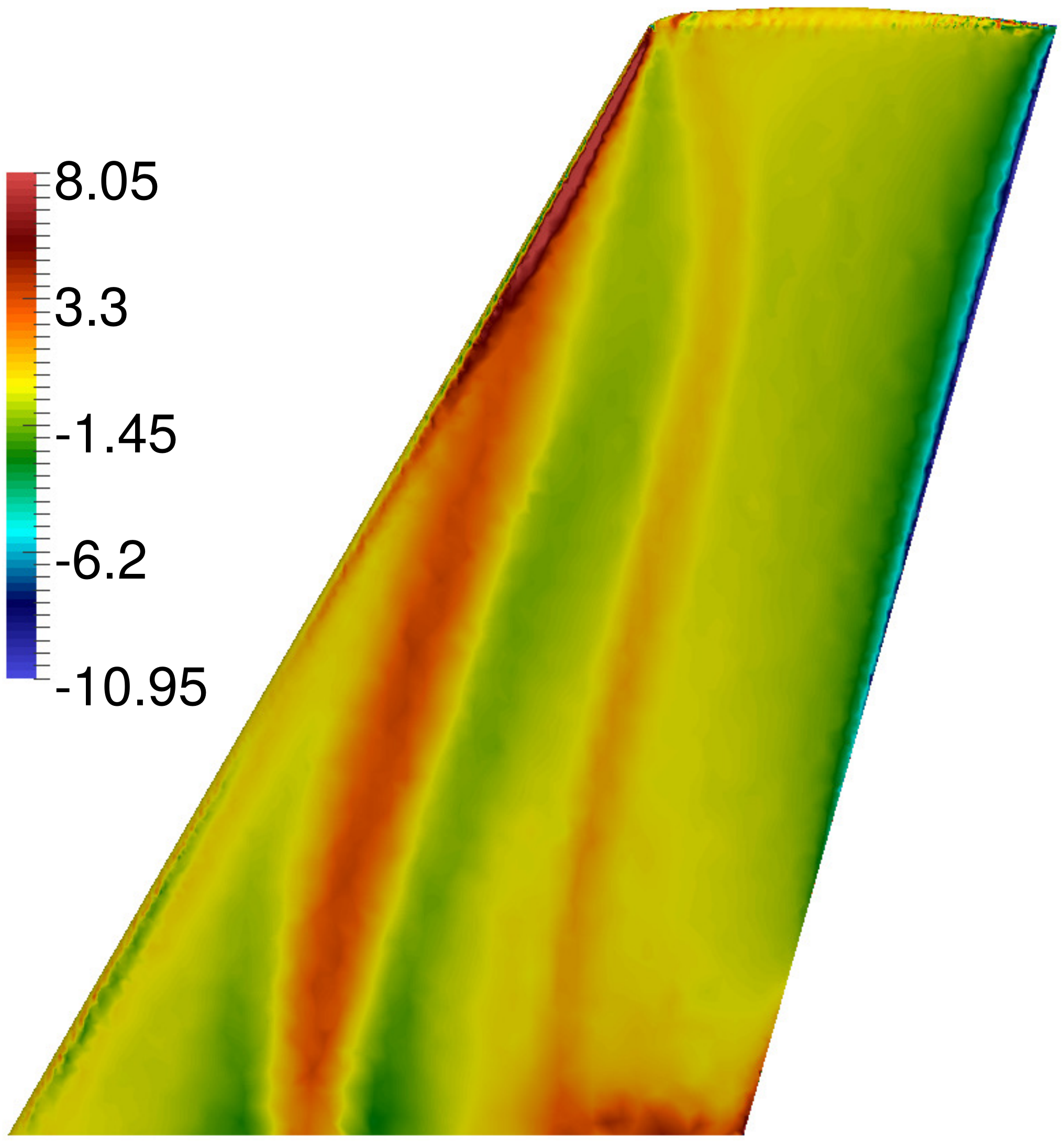} \\[\abovecaptionskip]
    \small (d) Conservative mapping with enhanced mortar method.
  \end{tabular}

  \caption{Interface energy sensitivity contours on the fluid mesh at the upper surface of the flexible ONERA M6. The results are shown for matching and non-matching interfaces using different mapping techniques.}\label{fig:oneraM6_AFSI_mapping_comparison}
\end{figure}

\begin{figure}
  \centering
  \begin{tabular}{@{}c@{}}
    \includegraphics[width=1\linewidth,keepaspectratio]{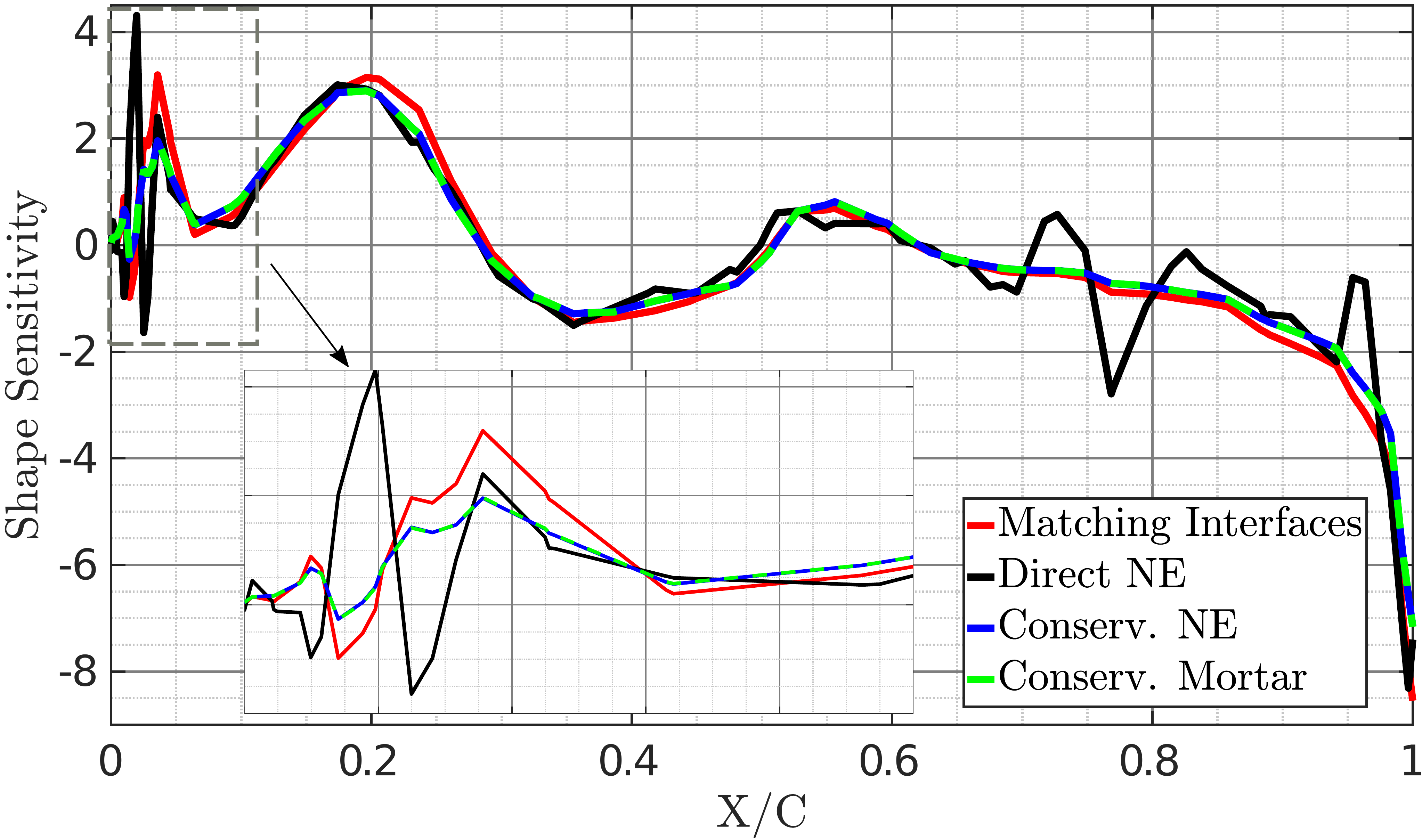} \\[\abovecaptionskip]
    \small Upper surface.
  \end{tabular}

  \vspace{\floatsep}

  \begin{tabular}{@{}c@{}}
    \includegraphics[width=1\linewidth,keepaspectratio]{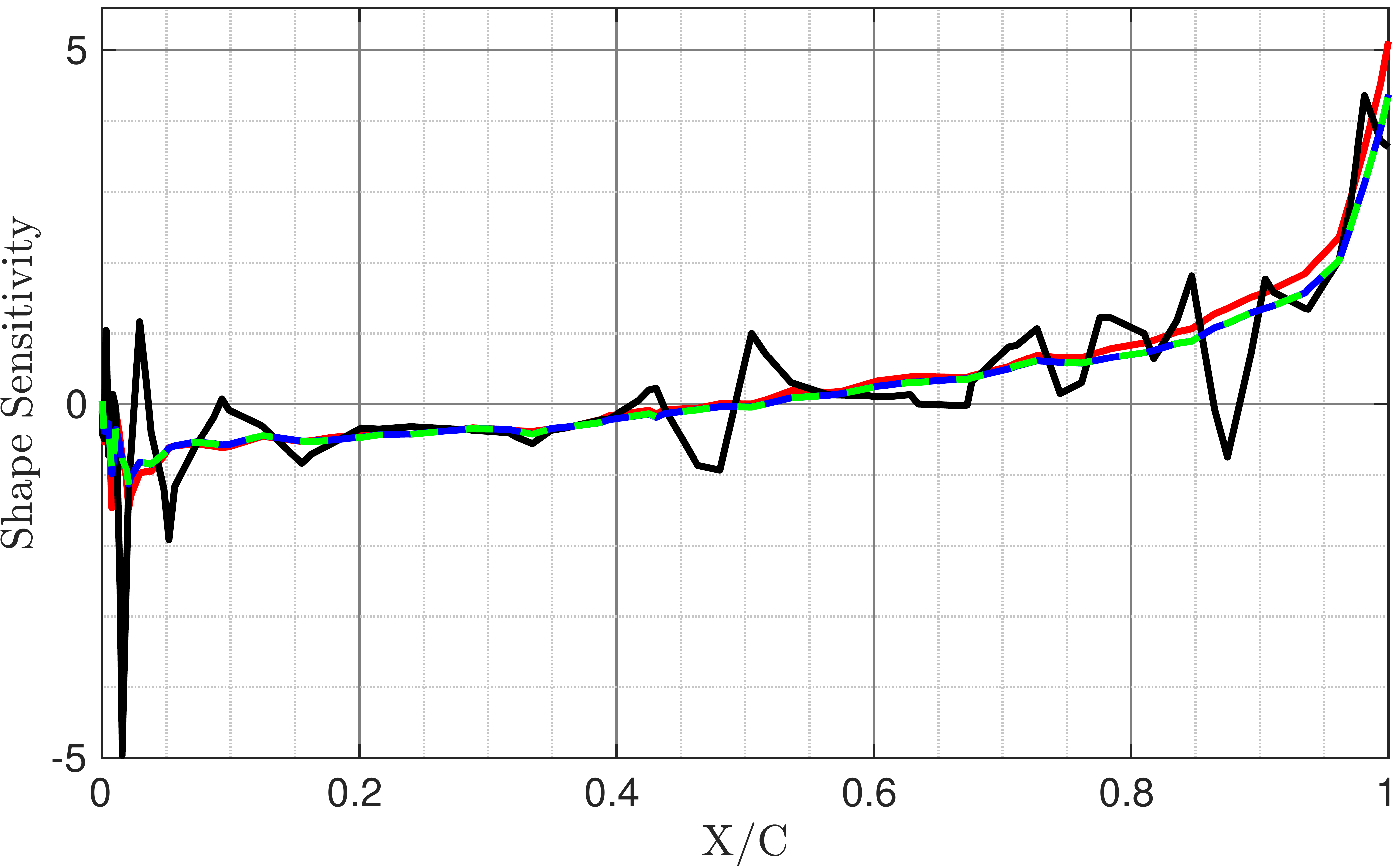} \\[\abovecaptionskip]
    \small Lower surface.
  \end{tabular}
  
  \caption{Profiles of interface energy shape sensitivity in \\ Fig. \ref{fig:oneraM6_AFSI_mapping_comparison} at Y/b=$0.65$.}\label{fig:oneraM6_AFSI_mapping_comparison_cross_sections}
\end{figure}

Lastly, we assess the accuracy of the interface sensitivity information obtained with the complete three-field formulation for the cases with non-matching interface meshes. For this purpose, nodes on the fluid interface mesh in the undeformed configuration are taken as design variables, i.e. $\boldsymbol{X}_{\mathcal{D}} = \boldsymbol{X}^{\mathcal{F}}_{\Gamma_{\mathcal{I}}} $. Since interface shape derivatives of the spatial mapping matrices, i.e. $\frac{\partial \boldsymbol{H}^{i}}{\partial \boldsymbol{X}^{j}_{\Gamma_{\mathcal{I}}}}, \{i,j\}\in \{\mathcal{F},\mathcal{S}\}$, are not normally available in coupling tools like EMPIRE, the simplified coupled sensitivity equation (Eq. \ref{e:partitioned_Discrete_Adj_Sens_Eq_simplified}) is used, which reads as follows:    
\begin{equation}
\label{e:ONERA_SENS_EQU}
\frac{d\mathcal{L}}{d\boldsymbol{X}^{\mathcal{F}}_{\Gamma_{\mathcal{I}}}} \approx 
\frac{\partial \mathcal{L}^{\mathcal{F}}}{\partial \boldsymbol{x}^{\mathcal{F}}_{\Gamma_{\mathcal{I}}}} +  \frac{\partial \mathcal{L}^{\mathcal{M}}}{\partial \boldsymbol{X}^{\mathcal{F}}_{\Gamma_{\mathcal{I}}}} + \frac{\partial \mathcal{L}^{\mathcal{S}}}{\partial \boldsymbol{X}^{\mathcal{S}}_{\Gamma_{\mathcal{I}}}} \cdot (\boldsymbol{H}^{\mathcal{S}})^T.
\end{equation}
$\boldsymbol{H}^{\mathcal{S}}$ is applied directly on the structural sensitivities to ensure that a constant sensitivity field is mapped exactly on the fluid interface mesh. In order to critically evaluate the spatial coupling techniques used in Section \ref{Steady_state_aeroelastic_analysis_of_ONERA_M6_wing}, the interface energy is chosen for the aeroelastic shape sensitivity analysis. Figures \ref{fig:oneraM6_AFSI_mapping_comparison} and \ref{fig:oneraM6_AFSI_mapping_comparison_cross_sections} compare qualitatively and quantitatively the results obtained with each mapping technique. Comparisons show good agreement between the conservative mapping results and the reference. Interestingly, the oscillatory behavior of the NE mapper in the primal problem is reversed in the sensitivity analysis. This can be explained by the fact that, in the adjoint sensitivity analysis all operations are transposed w.r.t the primal problem (see Fig. \ref{fig:AFSI_three_field}). Specially, mapping operations on the interface are reversed as
 \begin{subequations}
 	\label{eq:last_eq}
	\begin{align}
	\boldsymbol{d} & = (\boldsymbol{H}^{\mathcal{F}})^T \cdot \boldsymbol{\Psi}^{\mathcal{S}}_{\Gamma_\mathcal{I}} \\
	\boldsymbol{f}^{\mathcal{S},a}_{\Gamma_\mathcal{I}} & = (\boldsymbol{H}^{\mathcal{S}})^T \cdot \boldsymbol{\Psi}^{\mathcal{M}}_{\Gamma_\mathcal{I}}.
	\end{align}
\end{subequations}
In the literature \citep{Wang2016papaer,DeBoer2008} and also here (Section \ref{Steady_state_aeroelastic_analysis_of_ONERA_M6_wing}), it has been observed that the conservative force transfer with the NE mapper, i.e. $\boldsymbol{f}^{\mathcal{S}}_{\Gamma_\mathcal{I}} = (\boldsymbol{H}^{\mathcal{S}})^T \cdot  \boldsymbol{f}^{\mathcal{F}}_{\Gamma_\mathcal{I}}$, produces spurious oscillations, while the direct mapping, $\boldsymbol{f}^{\mathcal{S}}_{\Gamma_\mathcal{I}} = \boldsymbol{H}^{\mathcal{F}} \cdot  \boldsymbol{f}^{\mathcal{F}}_{\Gamma_\mathcal{I}}$, delivers accurate and oscillation-free traction field on the structure interface. Therefore we can associate the noisy behavior of the NE mapper in the coupled adjoint sensitivity analysis with the fact that transposed NE mapping matrices are used to transfer interface adjoint displacements of fluid and structure (see Eq. \ref{eq:last_eq}). On the other hand, from Figures \ref{fig:oneraM6_AFSI_mapping_comparison} and \ref{fig:ONERA_M6_FSI_traction_compariosn}, it is observed that the mortar method does not introduce noise in neither primal fields nor coupled shape sensitivity field. 

As seen from Figure \ref{fig:oneraM6_AFSI_mapping_comparison_cross_sections}, there are local  differences between the conservative-mapping-based gradients and the reference. They can be explained by the following facts:
\begin{itemize}
	\item Primal fields (displacements and tractions) computed on the non-matching interfaces have quantitative inaccuracies w.r.t those computed on the matching interfaces. As a result, adjoint fields and subsequently sensitivities can not be expected to be the same as the references values.
	\item Shape derivatives of the mapping matrices are neglected in Eq. \ref{eq:last_eq}. \cite{Wang2016papaer} has shown that the quality of spatial mapping results can be deteriorated at curved boundaries like leading edge. This means that spatial mapping operation is sensitive to the shape of interface surface. Therefore, it is reasonable to attribute discrepancies observed in the leading edge region (see Figure \ref{fig:oneraM6_AFSI_mapping_comparison_cross_sections}) to the omission of spatial mapping sensitivities.  
\end{itemize}

\section{Conclusions}
\label{conclusion_sec}
In this paper, adjoint-based shape sensitivity analysis for FSI problems was revisited from the mathematical and, particularly, the black-box implementation point of view. To exploit advantages of existing single-disciplinary solvers, a mixed Lagrangian-Eulerian formulation was used to solve the FSI problem. In a manner consistent with the primal problem, the adjoint FSI problem was partitioned 
using coupling conditions which were realized as auxiliary objective functions for single-disciplinary adjoint solvers. This requires a minimal modification to existing solvers. The proposed partitioned scheme projects the shape sensitivities of multidisciplinary objective functionals to the undeformed configuration, which is considered to be a great advantage. 

The presented scheme requires domain-based adjoint sensitivities of the fluid to be transferred as adjoint forces to the adjoint mesh motion problem. Since adjoint fluid solver might not use the domain-based formulation rather the boundary-based one (the so-called reduced formulation), a reduced adjoint coupling scheme was also developed. Although the complete formulation unconditionally showed accurate coupled shape gradients, the reduced one was suffering from accuracy issues in regions of strong flow gradients and near singularities. 
 
This paper also investigated the performance of two spatial mapping techniques in primal and adjoint shape sensitivity analyses of FSI problems involving non-matching interface meshes. Tests with a representative aeroelastic wing showed that the conservative mortar method, unlike the nearest element interpolation, does not introduce spurious oscillations, neither in the interface traction received by the structure nor in the interface shape sensitivity field.

\section{Replication of results}
The software packages used in this work are open-source and available for download at the URLs given in the list of references. Furthermore, the datasets analyzed during the current study are available in the following link: \url{https://1drv.ms/f/s!AkrOhpK6P2FWgYVD-M7VHT37MZhYlA}. 

\section*{Conflict of interest}
On behalf of all authors, the corresponding author states that there is no conflict of interest. 
\section*{Acknowledgements}
The authors gratefully acknowledge the support of the International Graduate School of Science and Engineering (IGSSE) of the Technische Universität München, Germany, under project 9.10. We are also grateful to the editor and anonymous reviewers who contributed to the improvement of the paper quality by their insightful and constructive comments.

\bibliographystyle{apalike}       

\bibliography{mybibfile}

\end{document}